\documentclass[reqno]{amsart}
\usepackage{amssymb,latexsym}
\usepackage{amsfonts,mathrsfs}
\usepackage[margin=1.4in]{geometry}
\usepackage[all]{xy}

\newtheorem{claim}{Claim}[section]

\newtheorem{theorem}{Theorem}[section]
\newtheorem{lemma}[theorem]{Lemma}
\newtheorem{proposition}[theorem]{Proposition}
\newtheorem{corollary}[theorem]{Corollary}
\newtheorem{definition}[theorem]{Definition}

\newcommand{\ba}{\begin{array}}
\newcommand{\ea}{\end{array}}

\def \qed{\cqfd}

\newcommand{\Vol}{\mathrm{Vol}}
\newcommand{\rank}{\mathrm{rank}}
\newcommand{\Id}{\mathrm{Id}}
\newcommand{\End}{\mathrm{End}}

\def\qed{\vbox{\hrule
\hbox{\vrule\hbox to 5pt{\vbox to 8pt{\vfil}\hfil}\vrule}\hrule}}

\newcommand{\beg}{\begin{eqnarray*}}
\newcommand{\begn}{\begin{eqnarray}}
\newcommand{\en}{\end{eqnarray*}}
\newcommand{\enn}{\end{eqnarray}}

\newcommand{\tr}{\mbox{\rm tr\,}}

\begin{document}
\title{The limit of the Hermitian-Yang-Mills flow  on reflexive sheaves}
\author{Jiayu Li}
\address{Key Laboratory of Wu Wen-Tsun Mathematics\\ Chinese Academy of Sciences\\School of Mathematical Sciences\\
University of Science and Technology of China\\
Hefei, 230026\\ and AMSS, CAS, Beijing, 100080, P.R. China\\} \email{jiayuli@ustc.edu.cn}
\author{Chuanjing Zhang}
\address{School of Mathematical Sciences\\
University of Science and Technology of China\\
Hefei, 230026,P.R. China\\ } \email{chjzhang@mail.ustc.edu.cn}
\author{Xi Zhang}
\address{Key Laboratory of Wu Wen-Tsun Mathematics\\ Chinese Academy of Sciences\\School of Mathematical Sciences\\
University of Science and Technology of China\\
Hefei, 230026,P.R. China\\ } \email{mathzx@ustc.edu.cn}
\thanks{The authors were supported in part by NSF in
China,  No. 11625106, 11571332, 11131007.}

\begin{abstract} In this paper, we study the asymptotic behavior of the Hermitian-Yang-Mills flow on a reflexive sheaf. We prove that the limiting reflexive sheaf is isomorphic to the double dual of
the graded sheaf associated to the Harder-Narasimhan-Seshadri filtration, this answers a question by Bando and Siu.
\end{abstract}

\maketitle

\section{Introduction}
\setcounter{equation}{0}

\hspace{0.4cm}

Let $(M, \omega )$ be a compact K\"ahler manifold and $\mathcal{E}$ a  coherent sheaf on $M$.
 The $\omega$-degree and the $\omega$-slope of $\mathcal{E}$ are defined by
\begin{equation*}
\deg_{\omega}(\mathcal{E}):=\int_X c_1(\mathcal{E})\wedge\displaystyle{\frac{\omega^{n-1}}{(n-1)!}},
\end{equation*}
and
\begin{equation*}
\mu_{\omega}(\mathcal{E}):=\frac{\deg_{\omega}(\mathcal{E})}{\mbox{rank}(\mathcal{E})},
\end{equation*}
where $c_1(\mathcal{E})$ is the first Chern class of $\mathcal{E}$. We say that a torsion free coherent sheaf $\mathcal{E}$  is $\omega$-stable ($\omega$-semi-stable) in the sense of Mumford-Takemoto if for every proper coherent sub-sheaf $\mathcal{F}\hookrightarrow \mathcal{E}$ we have
\begin{equation*}
\mu_{\omega}(\mathcal{F})< (\leq )\mu_{\omega}(\mathcal{E}).
\end{equation*}
We denote by $\Sigma_{\mathcal{E}}$ the set of singularities where $\mathcal{E}$ is not locally free. It is well known that the coherent sheaf $\mathcal{E}$ can be seen as a holomorphic vector bundle on $M\setminus \Sigma_{\mathcal{E}}$.
 A Hermitian metric $H$ on the sheaf $\mathcal{E}$ is called {\it admissible } if it is a Hermitian metric which is  defined on the  holomorphic vector bundle $ \mathcal{E}|_{M\setminus \Sigma_{\mathcal{E}}}$ and satisfies: (1) $|F_{H}|_{H, \omega}$ is square integrable;
(2) $|\Lambda_{\omega } F_{H}|_{H}$ is uniformly bounded,
where $F_{H}$ is the curvature tensor of Chern connection $D_{H}$ with respect to  $H$, and $\Lambda_{\omega }$ denotes the contraction with  the K\"ahler metric $\omega $.
 A Hermitian metric $H$ on the  holomorphic vector bundle $\mathcal{E}|_{M\setminus \Sigma_{\mathcal{E}}}$ is said to be
$\omega$-Hermitian-Einstein  if it
  satisfies the following Einstein condition on $M\setminus \Sigma_{\mathcal{E}}$, i.e.
\begin{equation}\label{HE}
\sqrt{-1}\Lambda_{\omega} F_{H}
=\lambda_{\mathcal{E}, \omega} \Id_{\mathcal{E}},
\end{equation}
where   $\lambda_{\mathcal{E}, \omega}  =\frac{2\pi}{\Vol(M, \omega)} \mu_{\omega} (\mathcal{E})$.

The Donaldson-Uhlenbeck-Yau theorem (\cite{NS,D0,D1,UY}) states that, if $\mathcal{E}$ is locally free on the whole $M$, i.e. $\Sigma_{\mathcal{E}}=\emptyset$,  the $\omega$-stability of $\mathcal{E}$ implies the existence of  $\omega$-Hermitian-Einstein metric on $\mathcal{E}$. This theorem has several interesting and important generalizations and extensions (\cite{LY, Hi, Si, Br1, BS, GP, BG, HL, AG1, Bi, BT, JZ, LN1, LN2, M, Mo1, Mo}, etc.).
In \cite{BS},
 Bando and Siu  introduced the notion of admissible Hermitian metrics on torsion-free sheaves, and proved the Donaldson-Uhlenbeck-Yau theorem on stable reflexive sheaves.
In fact, they obtained  a long time solution $H(t)$ of the Hermitian-Yang-Mills flow on $(M\setminus \Sigma_\mathcal{E}) \times [0, +\infty)$, i.e. $H(t)$ satisfies:
\begin{equation}\label{SSS1}
\left \{\begin{split} &H^{-1}(t)\frac{\partial H(t)}{\partial
t}=-2(\sqrt{-1}\Lambda_{\omega}F_{H(t)}-\lambda_{\mathcal{E}, \omega } \Id_{\mathcal{E}}),\\
&H(0)=\hat{H},\\
\end{split}
\right.
\end{equation}
where $\hat{H}$ is an initial metric which will be described in section $2$ in details. The above Hermitian-Yang-Mills flow was introduced and studied  by Donaldson in \cite{D0,D1}. Bando and Siu have shown that $H(t)$ is admissible for every $t>0$. Furthermore, they proved that: if the reflexive coherent sheaf $\mathcal{E}$ is $\omega$-stable, then along the Hermitian-Yang-Mills flow, $H(t)$ converges to $ H_{\infty}$ subsequently  in weak $L_{2, loc}^{p}$-topology  and $H_{\infty}$ is an  admissible $\omega$-Hermitian-Einstein metric.
 There are also some results on the  existence of  approximate solution of  Hermitian-Einstein equation (\ref{HE}) on a semi-stable holomorphic bundle and a semi-stable Higgs bundle, see references
\cite{Ko2,Ja,Card1,Card2,LZ2} for details. Recently, the authors (\cite{LZZ1}) obtain the existence of  admissible approximate  $\omega$-Hermitian-Einstein structure on an $\omega$-semi-stable reflexive sheaf, i.e. they proved that, if the reflexive  sheaf $\mathcal{E}$ is $\omega$-semi-stable, along the Hermitian-Yang-Mills flow (\ref{SSS1}), we have
\begin{equation}\label{SS01}
\sup _{x\in M\setminus \Sigma_\mathcal{E} } |\sqrt{-1}\Lambda_{\omega }F_{H(t)}-\lambda_{\mathcal{E}, \omega } \Id_{\mathcal{E}}|_{H(t)}(x)\rightarrow 0,
\end{equation}
as $t\rightarrow +\infty $.

For an unstable torsion-free coherent sheaf $\mathcal{E}$, one can associate a
filtration (\cite{Ko2}, \cite{Br}) by  sub-sheaves
\begin{equation}\label{HNS01}
0=\mathcal{E}_{0}\subset \mathcal{E}_{1}\subset \cdots \subset \mathcal{E}_{l}=\mathcal{E} ,
\end{equation}
such that every quotient sheaf $Q_{i}=\mathcal{E}_{i}/\mathcal{E}_{i-1}$ is torsion-free and $\omega$-stable, which is called the
Harder-Narasimhan-Seshadri filtration of the reflexive sheaf $\mathcal{E}$ (abbr.
HNS-filtration). Moreover, $\mu (Q_{i})\geq \mu (Q_{i+1})$ and the
associated graded object
\begin{equation}Gr^{HNS}(\mathcal{E}
)=\oplus_{i=1}^{l}Q_{i}
\end{equation} is uniquely determined by the isomorphism
class of $\mathcal{E}$ and the K\"ahler class $[\omega]$.

 If the reflexive sheaf $\mathcal{E}$ is not stable, Bando and Siu (\cite{BS})  proved that: there exists a subsequence $H(t_{i})$ along the Hermitian-Yang-Mills flow (\ref{SSS1}) such that $\int_{M}|\nabla \Lambda_{\omega}F_{H(t_{i})}|\frac{\omega^n}{n!} \rightarrow 0$ as $t_i\rightarrow +\infty$. By Uhlenbeck's theorem (\cite{U1,UY}), taking suitable complex gauge transformations one can choose a subsequence $t_i\rightarrow +\infty$ so that Chern connections
  $ A(t_{i})\rightarrow A_{\infty}$ weakly in $L_1^{2}$-topology outside a closed subset $\tilde{\Sigma}\subset M$ of Hausdorff codimension at least $4$. Since       $\sqrt{-1}\Lambda_{\omega}F_{A_{\infty}}$ is parallel, we can decompose $E_{\infty}$ according to the
eigenvalues of $\sqrt{-1}\Lambda_{\omega}F_{A_{\infty}}$ on $M\setminus \tilde{\Sigma}$. Then we  obtain a holomorphic
orthogonal decomposition
\begin{equation}
E_{\infty}=\bigoplus_{i=1}^{l}E_{\infty}^{i},
\end{equation}
every $E_{\infty}^{i}$ admits a Hermitian-Einstein metric and can be extended to a reflexive sheaf. In \cite{BS}, Bando and Siu propose an interesting question: whether
\begin{equation}
\bigoplus_{i=1}^{l}E_{\infty}^{i} \cong Gr^{HNS}(\mathcal{E})^{\ast \ast}.
\end{equation}
Atiyah and Bott (\cite{AB}) first raised the same question for Riemann surfaces case, which has been proved by  Daskalopoulos (\cite{Da}). When $\mathcal{E}$ is locally free on the whole $M$, the conjecture was confirmed by Daskalopoulos and Wentworth (\cite{DW1}) for  K\"ahler surfaces case;  by Jacob (\cite{Ja3}) and Sibley (\cite{Sib})  for  higher dimensional case.
The above Atiyah-Bott-Bando-Siu conjecture is also valid for Higgs bundles, see references
\cite{Wi,LZ1,LZ3} for details.  In this paper, we study  the  asymptotic behavior of the Hermitian-Yang-Mills flow (\ref{SSS1}) on a reflexive sheaf $\mathcal{E}$, and give a confirm answer to the above Bando-Siu's question.
We obtain the following theorem.

\medskip

\begin{theorem}\label{thm 1.1}
Let $\mathcal{E}$ be a reflexive sheaf on a compact K\"ahler manifold $(M, \omega )$, and $H(t)$ be the solution of the  Hermitian-Yang-Mills flow (\ref{SSS1}) on  $\mathcal{E}$ with the initial metric $\hat{H}$. We have a family of integrable connections
\begin{equation}
A(t)=g (t) (\hat{A})
\end{equation}
on $\mathcal{E}|_{M\setminus \Sigma_{\mathcal{E}}}$ for $t\in [0, +\infty )$, where $g(t)\in {\mathcal G }^{\mathbb{C}}$ satisfies  $g
^{\ast H_{0}}(t)\circ g(t)=h (t)=H_{0}^{-1}H(t)$,  $\Sigma_{\mathcal{E}}$ is the singularity set of $\mathcal{E}$ and $\hat{A}$ is the Chern connection with respect to the initial metric $\hat{H}$, such that:

(1) For every sequence $\{t_{i}\}$ there exists a subsequence
$\{ t_{j}\}$ such that, $A(t_{j})$ converges, modulo gauge transformations, to a Yang-Mills connection
 $A_{\infty}$ on a Hermitian vector bundle $(E_{\infty}, H_{\infty})$ over $M\setminus \Sigma $ in $C_{loc}^{\infty}$-topology as $t_{j}\rightarrow +\infty$, where $\Sigma\subset M$ is
a closed set of Hausdorff codimension at least $4$. Furthermore, the limiting $(E_{\infty}, \overline{\partial }_{A_{\infty}})$ can be extended to the whole $M$ as a reflexive  sheaf  with a holomorphic orthogonal splitting
  \begin{equation}
(E_{\infty}, H_{\infty}, \overline{\partial}_{A_{\infty}}
)=\bigoplus_{i=1}^{l}(\mathcal{E}_{\infty}^{i} , H_{\infty}^{i}
),
\end{equation}
where $H_{\infty}^{i}$ is an admissible Hermitian-Einstein metric on the reflexive
 sheaf $\mathcal{E}_{\infty}^{i} $.

(2) Moreover, the extended reflexive  sheaf is isomorphic to the double dual of the graded  sheaf associated to the HNS-filtration of $\mathcal{E}$, i.e. we have
\begin{equation}
(E_{\infty}, \overline{\partial }_{A_{\infty}})\simeq Gr^{HNS}(\mathcal{E})^{\ast \ast }.
 \end{equation}

\end{theorem}

\medskip

We now give an overview of our proof. The conclusion in the the part (1) of Theorem \ref{thm 1.1} is stronger than that in Theorem 4 in \cite{BS}, because we prove that the convergence holds not only for every sequence $\{t_{i}\}$ but also in much stronger topology, i.e. in $C_{loc}^{\infty}$-topology. To prove the part (1), we follow Hong-Tian's argument  in \cite{HT}. Even though the global approach is similar, some key estimates require new analytical ideas because the base manifold $M\setminus \Sigma_{\mathcal{E}} $ in our case is not compact.
For examples: to prove that $\int_{M}|D_{H(t)}(\Lambda_{\omega}F_{H(t)})|_{H(t)}^{2}\frac{\omega ^{n}}{n!}\rightarrow 0$ as $t\rightarrow +\infty$ in Proposition \ref{I}; to analyze the limiting behavior of the  Yang-Mills flow  on $\mathcal{E}|_{M\setminus\Sigma_{\mathcal{E}}}$ in Theorem \ref{thmlimit}.

To prove the second part of Theorem \ref{thm 1.1}, we will use
the basic idea in \cite{DW1} for a locally free sheaf in the K\"ahler surface case, but there are two points where we need new arguments for reflexive sheaves case. The first one is to prove that the HN type of the limiting  sheaf is in fact equal to that of $\mathcal{E}$; and the second one is to construct a non-zero holomorphic map from any stable quotient  sheaf in HNS-filtration of $\mathcal{E}$ to the limiting sheaf.

The first one is closely related to the existence of an $L^{p}$-approximate critical Hermitian metric (as defined in \cite{DW1}). When $\mathcal{E}$ is locally free, Sibley (\cite{Sib}) constructs a resolution of the HNS-filtration of $\mathcal{E}$ by subbundles, i.e. there exists a finite sequence of blow-ups with smooth centers such that the pullback bundle $\pi^{\ast }\mathcal{E}$ has a filtration by subbundles, where $\pi : \tilde{M}\rightarrow M $ is the composition of the blow-ups involved in the resolution. The metric $\pi^{\ast}\omega $ is degenerate along the exceptional divisor $\pi^{-1}(\Sigma_{HNS})$, where $\Sigma_{HNS}$ is the singularity set of the HNS-filtration of $\mathcal{E}$, and it can be approximated by a family of K\"ahler metrics $\omega_{\epsilon}$ on $\tilde{M}$.  Since every  quotient  subbundle is $\omega_{\epsilon}$-stable  for small $\epsilon$, one can use Donaldson-Uhlenbeck-Yau theorem  to take the direct sum of the Hermitian-Einstein metrics on  quotient  subbundles in the resolution. By choosing any fixed smooth Hermitian metric $H_{0}$ on $\pi^{\ast}\mathcal{E}$ over a neighborhood of $\pi^{-1}(\Sigma_{HNS})$ such that $|\Lambda_{\omega_{\epsilon}}F_{H_{0}}|_{H_{0}}$ is uniformly bounded, Sibley uses Daskalopoulos and Wentworth's cut-off argument (\cite{DW1}) to obtain a smooth $L^{p}$-approximate critical Hermitian metric on the locally free sheaf $\mathcal{E}$. In our case,  $\mathcal{E}$ is only reflexive, we can not find such  smooth metric $H_{0}$. So we can not use Sibley's result directly, and need new arguments to obtain a smooth $L^{p}$-approximate critical Hermitian metric, see Proposition \ref{claim} and Proposition \ref{prop 4.1} for details. Furthermore, in Lemma \ref{continuous}, we prove  the continuous dependence of the Hermitian-Yang-Mills flow (\ref{SSS1}) on initial metrics, this is fully nontrivial for noncompact base manifolds case. Then we can follow Daskalopoulos and Wentworth's trick (Lemma 4.3 in \cite{DW1}) to prove that the HN type of the limiting  sheaf is in fact equal to that of $\mathcal{E}$.

For the second one, we use Donaldson's idea (\cite{D1}) to  construct a nonzero holomorphic map to the limiting bundle as the limit of the sequence of gauge transformations  defined by the Yang-Mills flow. There are many difficulties to obtain uniform estimates, because we have no uniform $L^{\infty}$-bound on the mean curvature (i.e. $|\sqrt{-1}\Lambda_{\omega}F_{A}|$) of the induced connection for subsheaves.
 Using the resolution of singularities, we can pull back the HNS-filtration to $\tilde{M}$ by subbundles. Evolving the induced Hermitian metric on the subbundle by the Hermitian-Yang-Mills flow with respect to the K\"ahler metric $\omega_{\epsilon}$, by the result in \cite{BS}, we can get a uniform $L^{\infty}$-bound on the mean curvature and a local uniform $C^{0}$-estimate of the evolved Hermitian metrics. Using these estimates and following the argument in Proposition 4.1 in \cite{LZ3}, we can obtain a local uniform $C^{0}$-estimate of a sequence of holomorphic maps and then construct a nonzero holomorphic map to the limiting bundle. It should be pointed out that in Proposition 4.1 in \cite{LZ3}, we need the assumption that the pulling back geometric objects including the complex gauge transformations and induced metrics on the subsheaves  can be extended smoothly on the whole $\tilde{M}$. This assumption may not be satisfied in our case, so we modify the argument in \cite{LZ3} suitably to the case that  the geometric object we consider can be approximated by a sequence of smooth ones, see Proposition \ref{propmap} for details.

This paper is organized as follows. In Section 2, we recall Bando and Siu's  regularization on the reflexive sheaf and some basic estimates for the Hermitian-Yang-Mills flow, and we prove that along the Hermitian-Yang-Mills flow, $\int_{M}|D_{H(t)}\Lambda_{\omega}F_{H(t)}|_{H(t)}^{2}\frac{\omega ^{n}}{n!}\rightarrow 0$ as $t\rightarrow +\infty$. In section 3, we analyze the limiting behavior of the Yang-Mills flow on $(M\setminus\Sigma_{\mathcal{E}}, \mathcal{E}|_{M\setminus\Sigma_{\mathcal{E}}}, \omega)$ and give a proof for the part (1) of Theorem \ref{thm 1.1}.
  In section 4 and section 5, we obtain an $L^{p}$-approximate critical Hermitian metric and prove that the HN type of the limiting  sheaf is in fact equal to that of the initial one. In the last section,  we construct a non-zero holomorphic map between  sheaves and  complete the proof of Theorem \ref{thm 1.1}.

\hspace{0.3cm}

\section{Analytic preliminaries and basic estimates }
\setcounter{equation}{0}

In this section, we first recall Bando and Siu's  regularization on the reflexive sheaf, and then give some basic estimates for the Hermitian-Yang-Mills flow. Let $(M, \omega )$ be a compact K\"ahler manifold of complex dimension $n$, and $\mathcal{E}$  be a reflexive sheaf on $M$.  The singularity set of  $\mathcal{E}$ will be denoted by  $\Sigma_{\mathcal{E}}$.
Bando and Siu (\cite{BS}) proved that: there is a regularization on the reflexive sheaf $\mathcal{E}$, by successively blowing up $\pi_{i}: M_{i}\rightarrow M_{i-1}$ with smooth center $Y_{i-1}\subset M_{i-1}$ finite times such that the pull-back of $\mathcal{E}$ to $M_{k} $ modulo torsion  is locally free and the composition
  \begin{equation}
  \pi = \pi_{1}\circ \cdots \circ \pi_{k}: \tilde{M}\rightarrow M
  \end{equation}
  is biholomorphic outside $\Sigma_{\mathcal{E}} $, where $i=1, \cdots , k$, $M_{0}=M$ and $\tilde{M}=M_{k}$. It is easy to see that the holomorphic vector bundle $E=\pi^{\ast}\mathcal{E}/tor(\pi^{\ast}\mathcal{E})$ is isomorphic to $\mathcal{E}$ on $\tilde{M}\setminus \pi^{-1}(\Sigma_{\mathcal{E}})$, where $tor(\pi^{\ast}\mathcal{E})$ is the torsion sheaf of $\pi^{\ast}\mathcal{E}$.

  It is well known that every $M_{i}$ is  K\"ahler (\cite{GH}). As in \cite{BS}, we fix arbitrary K\"ahler metrics $\eta_{i} $ on $M_{i}$ and set
 \begin{equation}\label{omegai}
 \omega_{1, \epsilon }=\pi_{1}^{\ast }\omega +\epsilon_{1}\eta_{1},  \quad  \quad  \omega_{i, \epsilon }=\pi_{i}^{\ast }\omega_{i-1 , \epsilon } +\epsilon_{i}\eta_{i}
 \end{equation}
for all $1\leq i \leq k$,
  where $0<\epsilon_{i} \leq 1$ and $\epsilon =(\epsilon_{1}, \cdots , \epsilon_{k})$. Bando and Siu (Lemma 3 in \cite{BS}) derived a uniform  Sobolev inequality for $(\tilde{M}, \omega_{\epsilon})$,
   by using Cheng and Li's estimate (\cite{CL}), they obtained the following uniform upper bounds of the heat kernels.

 \medskip

\begin{proposition}\label{Prop 2.1}{\bf(Proposition 2 in \cite{BS})}
Let $(M, \omega )$ be a compact K\"ahler manifold, and $\pi : \tilde{M}\rightarrow M$ be a single blow-up with smooth centre. Fix a K\"ahler metric $\eta $ on $\tilde{M}$ and set $\omega_{\epsilon }=\pi^{\ast}\omega +\epsilon \eta $, where $0<\epsilon \leq 1$.   Let $K_{\epsilon}$ be the heat kernel with respect to the metric $\omega_{\epsilon}$. Then, for any $\tau >0$, there exists a constant $C_{K}(\tau)$ independent of $\epsilon $, such that
\begin{equation}\label{kernel01}0\leq K_{\epsilon}(x , y, t)\leq C_{K}(\tau) (t^{-n}\exp{(-\frac{(d_{\omega_{\epsilon}}(x, y))^{2}}{(4+\tau )t})}+1)\end{equation} for every $x, y\in \tilde{M}$ and $0<t < +\infty$, where $d_{\omega_{\epsilon}}(x, y)$ is the distance between $x$ and $y$ with respect to the metric $\omega_{\epsilon}$. There also exists a constant $C_{G}$ such that
\begin{equation}\label{green} G_{\epsilon}(x , y)\geq -C_{G}\end{equation}
for every $x, y\in \tilde{M}$ and $0<\epsilon \leq 1$, where $G_{\epsilon}$ is the Green function with respect to the metric $\omega_{\epsilon}$.
\end{proposition}

 \medskip

  Given a smooth Hermitian metric $\hat{H}$ on the bundle $E$, we denote the corresponding Chern connection by $D_{\hat{H}}$, and the corresponding curvature form by $F_{\hat{H}}$.
\begin{equation}
|\Lambda_{\omega_{k, \epsilon }}F_{\hat{H}}|_{\hat{H}}\frac{\omega_{k, \epsilon}^{n}}{\eta_{k}^{n}}=n\Big|\frac{F_{\hat{H}}\wedge \omega_{k, \epsilon}^{n-1}}{\eta_{k}^{n}}\Big|_{\hat{H}}\leq \tilde{C}_{0}\sup_{\tilde{M}}|F_{\hat{H}}|_{\hat{H}},
\end{equation}
 where  $\tilde{C}_{0}$ is a uniform constant independent of $\epsilon$.
So there exists a uniform constant $\hat{C}_{0}$ such that
 \begin{equation}\label{initial1}
 \int_{\tilde{M}} |\Lambda_{\omega_{k, \epsilon }}F_{\hat{H}}|_{\hat{H}}\frac{\omega_{k, \epsilon}^{n}}{n!}\leq \hat{C}_{0},
 \end{equation}
 for all $ \epsilon $.

We consider the evolving metric $H_{k, \epsilon } (t)$  along  the Hermitian-Yang-Mills flow (\ref{SSS1}) on the holomorphic bundle $E $ over $\tilde{M}$ with the fixed smooth initial metric $\hat{H}$ and with respect to the K\"ahler metric $\omega_{k, \epsilon}$, i.e. it satisfies
\begin{equation}\label{DDD1}
\left \{\begin{split} &H_{k, \epsilon}^{-1}(t)\frac{\partial H_{k, \epsilon}(t)}{\partial
t}=-2(\sqrt{-1}\Lambda_{\omega_{k, \epsilon}}F_{H_{k, \epsilon}(t)} -\lambda_{k, \epsilon }\Id_{E}),\\
&H_{k, \epsilon}(0)=\hat{H},\\
\end{split}
\right.
\end{equation}
where  $\lambda_{k, \epsilon}  =\frac{2\pi}{\Vol(\tilde{M}, \omega_{k, \epsilon})} \mu_{\omega_{k, \epsilon }} (E)$.
 For simplicity, set:
\begin{equation}
\theta (H, \omega)=\sqrt{-1}\Lambda_{\omega}F_{H}-\lambda_{\omega }\Id_{E} .
\end{equation}
Along  the heat flow (\ref{DDD1}), we have the following estimates (the proof can be found in Siu's lecture notes \cite{Siu02}):
\begin{equation}\label{F1}
(\Delta_{k, \epsilon}-\frac{\partial }{\partial t})\tr (\theta (H_{k, \epsilon}(t), \omega_{k, \epsilon}))=0,
\end{equation}
\begin{equation}\label{F2}
(\Delta_{k, \epsilon}-\frac{\partial }{\partial t} )|\theta (H_{k, \epsilon}(t), \omega_{k, \epsilon})|_{H_{k, \epsilon}(t)}^{2}=2|D_{H_{k, \epsilon}(t)} (\theta (H_{k, \epsilon}(t), \omega_{k, \epsilon}))|^{2}_{H_{k, \epsilon}(t), \omega_{k, \epsilon}},
\end{equation}
\begin{equation}\label{H00}
(\Delta_{k, \epsilon } -\frac{\partial }{\partial t}) |\theta (H_{k, \epsilon}(t), \omega_{k, \epsilon})|_{H_{k, \epsilon}(t)}\geq  0,
\end{equation}
\begin{equation}\label{HH01}
 |\frac{\partial }{\partial t} \ln (\tr (H_{k, \epsilon}^{-1}(t_{0})H_{k, \epsilon}(t))+ \tr (H_{k, \epsilon}^{-1}(t)H_{k, \epsilon}(t_0)))|\leq 2|\theta (H_{k, \epsilon}(t), \omega_{k, \epsilon})|_{H_{k, \epsilon}(t)}.
\end{equation}

Using the maximum principle and the above inequalities, we derive
\begin{equation}\label{H001}
\int_{\tilde{M}}|\theta (H_{k, \epsilon}(t), \omega_{k, \epsilon}) |_{H_{k, \epsilon}(t)}\frac{\omega_{k, \epsilon }^{n}}{n!}\leq \int_{\tilde{M}}|\theta (\hat{H}, \omega_{k, \epsilon}) |_{\hat{H}}\frac{\omega_{k, \epsilon }^{n}}{n!}\leq \hat{C}_{1},
\end{equation}
 \begin{equation}\label{H000}
 |\theta (H_{k, \epsilon}(t), \omega_{k, \epsilon}) |_{H_{k, \epsilon}(t)}(x) \leq \int_{\tilde{M}}K_{k, \epsilon} (x, y, t)|\theta (\hat{H}, \omega_{k, \epsilon}) |_{\hat{H}}\frac{\omega_{k, \epsilon }^{n}}{n!} ,
\end{equation}
and
\begin{equation}\label{H0001}
 |\theta (H_{k, \epsilon}(t+1), \omega_{k, \epsilon}) |_{H_{k, \epsilon}(t+1)}(x) \leq \int_{\tilde{M}}K_{k, \epsilon} (x, y, 1)|\theta (H_{k, \epsilon}(t), \omega_{k, \epsilon}) |_{H_{k, \epsilon}(t)}\frac{\omega_{k, \epsilon }^{n}}{n!} ,
\end{equation}
for all $x\in \tilde{M}$ and $t>0$.

\medskip
%
%
%
%
%
%

After obtaining  local uniform $C^{\infty}$-bounds on $H_{k, \epsilon}(x, t)$, Bando and Siu (\cite{BS}) get the following lemma.

\medskip

\begin{lemma}\label{lemmaBS}{\bf (\cite{BS})}
By choosing a subsequence,  $H_{k, \epsilon}(x, t)$  converges successively  to a long time solution $H(x, t)$ of  the Hermitian-Yang-Mills flow (\ref{SSS1}) on $M\setminus \Sigma_{\mathcal{E}} \times [0, +\infty)$ in $C_{loc}^{\infty}$-topology as $(\epsilon_{1}, \cdots , \epsilon_{k})\rightarrow 0$.  Furthermore,  $H(x, t)$ is admissible and  satisfies:
\begin{equation}\label{H0013}
\int_{M}|\theta (H(t), \omega) |_{H(t)}\frac{\omega^{n}}{n!}\leq \int_{M}|\theta (\hat{H}, \omega) |_{\hat{H}}\frac{\omega^{n}}{n!}\leq \hat{C}_{1},
\end{equation}
\begin{equation}\label{H00013}
 |\theta (H(t+\tilde{t}), \omega ) |_{H(t+\tilde{t})}(x) \leq \int_{M}K_{\omega } (x, y, t)|\theta (H(\tilde{t}), \omega ) |_{H(\tilde{t})}\frac{\omega^{n}}{n!}
\end{equation}
for all $x\in M\setminus \Sigma_{\mathcal{E}}$, $t> 0$ and $\tilde{t}\geq 0$.
\end{lemma}

\medskip

Denote by $D_{\hat{A}}$ the Chern connection  on the holomorphic bundle $\mathcal{E}|_{M\setminus \Sigma_\mathcal{E}}$ with respect to the initial metric $\hat{H}$.
Let
$h(t)=\hat{H}^{-1}H(t)$, using the identities
\begin{equation}\label{id1}
\begin{split}
& \partial _{H(t)}-\partial_{\hat{H}} =h^{-1}(t)\partial_{\hat{H}}h(t) ,\\
& F_{H(t)}-F_{\hat{H}}=\overline{\partial }_{\hat{A}} (h^{-1}(t)\partial_{\hat{H}} h(t)) ,\\
\end{split}
\end{equation}  then we can rewrite (\ref{SSS1}) as
\begin{equation}
\frac{\partial h(t)}{\partial
t}=-2\sqrt{-1}h(t)\Lambda_{\omega}(F_{\hat{H}}+\overline{\partial
}_{\hat{A}}(h^{-1}(t)\partial_{\hat{H}}h(t)))+2\lambda_{\mathcal{E}, \omega} h(t).
\end{equation}

\medskip

Let's consider the Hermitian vector bundle $(\mathcal{E}|_{M\setminus \Sigma_{\mathcal{E}}}, \hat{H})$.
We denote by
$\textbf{A}_{\hat{H}}$ the space of connections of $\mathcal{E}|_{M\setminus \Sigma_{\mathcal{E}}}$
compatible with $\hat{H}$, by $\textbf{A}^{1,1}_{\hat{H}}$
 the space of unitary integrable connections of $\mathcal{E}|_{M\setminus \Sigma_{\mathcal{E}}}$, i.e.
\begin{equation}
\textbf{A}^{1,1}_{\hat{H}}=\{A\in \textbf{A}_{\hat{H}}|F_{A}^{0, 2}= F_{A}^{2, 0}=0\},
\end{equation}
and by $\textbf{G}^{\mathbb{C}}$ (resp. $\textbf{G}$, where $\textbf{G}=\{\sigma \in \textbf{G}^{\mathbb{C}}| \sigma^{\ast \hat{H}}\sigma=\Id\}$) the complex gauge group (resp. unitary gauge group) of the Hermitian
vector bundle $(\mathcal{E}|_{M\setminus \Sigma_{\mathcal{E}}}, \hat{H})$. $\textbf{G}^{\mathbb{C}}$ acts on the space $\textbf{A}_{\hat{H}}$ as follows: let
$\sigma \in \textbf{G}^{\mathbb{C}}$ and $A\in \textbf{A}_{\hat{H}}$,
\begin{equation}\label{id2}
\overline{\partial }_{\sigma(A)}=\sigma \circ \overline{\partial
}_{A}\circ \sigma^{-1}, \quad \partial _{\sigma (A)}=(\sigma^{\ast
\hat{H}})^{-1} \circ \partial _{A}\circ \sigma^{\ast\hat{H}}.
\end{equation}

In \cite{D1}, Donaldson has shown that the Hermitian-Yang-Mills flow (\ref{SSS1}) is formally gauge-equivalent to the
Yang-Mills flow, i.e. we have the following proposition:

\medskip

\begin{proposition}\label{propgauge}
There is a family of complex gauge transformations $\sigma (t)\in \textbf{G}^{\mathbb{C}}$ satisfying
$\sigma ^{\ast \hat{H}}(t)\sigma (t)=h (t)=\hat{H}^{-1}H(t)$, where $H(t)$ is the long time solution of the Hermitian-Yang-Mills flow (\ref{SSS1}) with the initial metric $\hat{H}$, such that $A(t)=\sigma (t)(\hat{A})$ is a long time solution of the Yang-Mills flow with the initial connection $\hat{A}$, i.e. it satisfies:
\begin{equation}\label{YM2}
\left \{\begin{split} &\frac{\partial A(t)}{\partial t}=-D_{A(t)}^{\ast } F_{A(t)},\\
&A(0)=\hat{A}.\\
\end{split}
\right.
\end{equation}
\end{proposition}

\medskip
It is well known that
\begin{equation}
\sigma^{-1} (t)\circ F_{A(t)}\circ \sigma (t)=F_{\hat{A}}+\overline{\partial }_{\hat{A}}(h^{-1}(t)\partial_{\hat{A}}h(t))=F_{H(t)},
\end{equation}
\begin{equation}
\sigma^{-1} (t)\circ D_{A(t)}(\Lambda_{\omega}F_{A(t)})\circ \sigma (t)=D_{H(t)}(\Lambda_{\omega}F_{H(t)}),
\end{equation}
and then
\begin{equation}
|F_{H(t)}|_{H(t)}^{2}=|F_{A(t)}|_{\hat{H}}^{2},
\end{equation}
\begin{equation}
|D_{H(t)}(\Lambda_{\omega}F_{H(t)})|_{H(t)}^{2}=| D_{A(t)}(\Lambda_{\omega}F_{A(t)})|_{\hat{H}}^{2}.
\end{equation}

\medskip

For simplicity, set
\begin{equation}
\theta (A(t), \omega ) =\sqrt{-1}\Lambda_{\omega}F_{A(t)}-\lambda_{\mathcal{E}, \omega }\Id,
\end{equation}
and
\begin{equation}
I(t)=\int_{M}|D_{A(t)}\theta (A(t), \omega)|_{\hat{H}}^{2}\frac{\omega ^{n}}{n!}=\int_{M}|D_{H(t)}\theta (H(t), \omega)|_{H(t)}^{2}\frac{\omega ^{n}}{n!}.
\end{equation}
In the following we will prove that $I(t)\rightarrow 0$ as $t\rightarrow +\infty $. When $\mathcal{E}$ is locally free, i.e. $\Sigma_{\mathcal{E}}=\emptyset $, this was prove by
Donaldson and Kronheimer (\cite{DK}). In the case that $\mathcal{E}$ is only reflexive, we need new arguments because the base manifold $M\setminus \Sigma_{\mathcal{E}}$ is non-compact.

\medskip

\begin{proposition}\label{I}
Let $H(t)$ be the long time solution of the Hermitian-Yang-Mills flow (\ref{SSS1}) with the initial metric $\hat{H}$, then $I(t)\rightarrow 0$ as $t\rightarrow +\infty $.
\end{proposition}

\medskip

{\bf Proof.} As that in the beginning of this section,   there is  a finite sequence of blowing up $\pi_{i}: M_{i}\rightarrow M_{i-1}$ with smooth center, where $i=1, \cdots , k$, such that $E=\pi^\ast\mathcal{E}/tor(\pi^\ast\mathcal{E})$ is locally free on $\tilde{M}$, where $\pi : \tilde{M}\rightarrow M$ is the composition of the sequence of blow-ups. The initial Hermitian metric $\hat{H}$ is a smooth metric on $E$. By induction, we can assume that there is just one blow-up, i.e. $k=1$.  Set $\omega_\epsilon=\pi^\ast\omega+\epsilon\eta$, where $\eta$ is a fixed K\"ahler metric on $\tilde{M}$.
Let $H_{ \epsilon } (t)$  be  the long time solution of the Hermitian-Yang-Mills flow (\ref{SSS1}) on the holomorphic bundle $E $ over $\tilde{M}$ with the fixed smooth initial metric $\hat{H}$ and with respect to the K\"ahler metric $\omega_{ \epsilon}$, i.e. it satisfies
\begin{equation}\label{DDD11}
\left \{\begin{split} &H^{-1}_{ \epsilon}(t)\frac{\partial H_{ \epsilon}(t)}{\partial
t}=-2(\sqrt{-1}\Lambda_{\omega_{ \epsilon}}F_{H_{\epsilon}(t)} -\lambda_{ \epsilon }\Id_{E}),\\
&H_{ \epsilon}(0)=\hat{H}.\\
\end{split}
\right.
\end{equation}
Lemma \ref{lemmaBS} says that $H_{ \epsilon}(x, t)$  converges  to the long time solution $H(x, t)$ of  the Hermitian-Yang-Mills flow (\ref{SSS1}) on $M\setminus \Sigma_\mathcal{E} \times [0, +\infty)$ in $C_{loc}^{\infty}$-topology as $\epsilon \rightarrow 0$. We also denote by $\hat{A}$ the Chern connection on the holomorphic vector bundle $E$ with respect to the smooth metric $\hat{H}$.
Let $A_\epsilon(t)$ be the long time solution of the Yang-Mills flow on the Hermitian vector bundle $(E, \hat{H})$ over the K\"ahler manifold $(\tilde{M}, \omega_{\epsilon})$, i.e.
\begin{equation}\label{YM1}
\left \{\begin{split}
& \frac{\partial A_\epsilon(t)}{\partial t}=-D_{A_\epsilon(t)}^{\ast} F_{A_\epsilon(t)},\\
& A_\epsilon(0)=\hat{A}.\\
\end{split}
\right.
\end{equation}
Set
\begin{equation}
I_{\epsilon}(t)=\int_{\tilde{M}}|D_{A_{\epsilon}(t)}\theta (A_{\epsilon}(t), \omega_{\epsilon})|_{\hat{H}}^{2}\frac{\omega_{\epsilon} ^{n}}{n!}=\int_{\tilde{M}}|D_{H_{\epsilon}(t)}\theta (H_{\epsilon}(t), \omega_{\epsilon})|_{H_{\epsilon}(t)}^{2}\frac{\omega_{\epsilon} ^{n}}{n!}.
\end{equation}

By the uniform bound on the heat kernel (\ref{kernel01}) and (\ref{H000}), there exists a uniform constant $C$ such that
\begin{equation}
\sup_{\tilde{M}}|\theta (A_{\epsilon}(t), \omega_{\epsilon})|_{\hat{H}}\leq \tilde{C}
\end{equation}
for  any $0<t_{0}\leq t $ and $0\leq \epsilon \leq 1$. Direct computations show that
\begin{equation}
\begin{split}
&\frac{dI_\epsilon(t)}{dt}=-2\int_{\tilde{M}}|D_{A_\epsilon(t)}^\ast D_{A_\epsilon(t)}\theta (A_{\epsilon}(t), \omega_{\epsilon})|_{\hat{H}}^2\frac{\omega_\epsilon^n}{n!}\\
&+2 Re\int_{\tilde{M}}\langle[\bar{\partial}_{A_\epsilon(t)}\theta (A_{\epsilon}(t), \omega_{\epsilon})-\partial_{A_\epsilon(t)}\theta (A_{\epsilon}(t), \omega_{\epsilon}), \theta (A_{\epsilon}(t), \omega_{\epsilon})], D_{A_\epsilon(t)}\theta (A_{\epsilon}(t), \omega_{\epsilon})\rangle _{\hat{H}}\frac{\omega_\epsilon^n}{n!}\\
\leq & 8\int_{\tilde{M}}|D_{A_\epsilon(t)}\theta (A_{\epsilon}(t), \omega_{\epsilon})|_{\hat{H}}^2|\theta (A_{\epsilon}(t), \omega_{\epsilon})|_{\hat{H}}\frac{\omega_\epsilon^n}{n!}\\
\leq& C I_{\epsilon}(t),
\end{split}
\end{equation}
where $C$ is a uniform constant. So we know that there exists a uniform constant such that
\begin{equation}I_\epsilon(t)\leq e^{C(t-s)}I_\epsilon(s),\end{equation} for any $0<t_{0}\leq s \leq t $ and $0\leq \epsilon \leq 1$.

Of course the formula (\ref{F2}) yields
\begin{equation}\label{moFF}
\begin{split}
&\int_{\tilde{M}}|\sqrt{-1}\Lambda_{\omega_\epsilon}F_{A_\epsilon(t)}-\lambda \Id|_{\hat{H}}^2\frac{\omega_\epsilon^n}{n!}+2\int_{t_0}^t\int_{\tilde{M}}|D_{A_\epsilon(s)}\theta (A_{\epsilon}(t), \omega_{\epsilon})|^2\frac{\omega_\epsilon^n}{n!}ds\\
=&\int_{\tilde{M}}|\sqrt{-1}\Lambda_{\omega_\epsilon}F_{A_\epsilon(t_0)}-\lambda \Id|_{\hat{H}}^2\frac{\omega_\epsilon^n}{n!}.
\end{split}
\end{equation}
According to Fatou's lemma, we get
\begin{equation}\label{mono1}
\begin{split}
&\int_{M\setminus\Sigma_\mathcal{E}}|\sqrt{-1}\Lambda_\omega F_{A(t)}-\lambda \Id|_{\hat{H}}^2\frac{\omega^n}{n!}+2\int_{t_0}^t\int_{M\setminus\Sigma_\mathcal{E}}\Big|\frac{\partial A(s)}{\partial s}\Big|^2\frac{\omega^n}{n!}ds\\
\leq &\int_{M\setminus\Sigma_\mathcal{E}}|\sqrt{-1}\Lambda_\omega F_{A(t_0)}-\lambda \Id|_{\hat{H}}^2\frac{\omega^n}{n!}.
\end{split}
\end{equation}
This implies that  $\int_{\tilde{M}}|\sqrt{-1}\Lambda_{\omega_\epsilon}F_{A_\epsilon(t)}-\lambda \Id|_{\hat{H}}^2\frac{\omega_\epsilon^n}{n!}$ and $\int_{M\setminus\Sigma_\mathcal{E}}|\sqrt{-1}\Lambda_\omega F_{A(t)}-\lambda \Id|_{\hat{H}}^2\frac{\omega^n}{n!}$ both are monotonically nonincreasing with respect to $t$. Then we must have
\begin{equation}
\int_{\tilde{M}}|\sqrt{-1}\Lambda_{\omega_\epsilon}F_{A_\epsilon(t)}-\lambda \Id|_{\hat{H}}^2\frac{\omega_\epsilon^n}{n!}-\int_{\tilde{M}}|\sqrt{-1}\Lambda_{\omega_\epsilon}F_{A_\epsilon(t+1)}-\lambda \Id|_{\hat{H}}^2\frac{\omega_\epsilon^n}{n!}\rightarrow 0
\end{equation}
and
\begin{equation}\label{moFF2}
\int_{M\setminus\Sigma_\mathcal{E}}|\sqrt{-1}\Lambda_\omega F_{A(t)}-\lambda \Id|_{H_0}^2\frac{\omega^n}{n!}-\int_{M\setminus\Sigma_\mathcal{E}}|\sqrt{-1}\Lambda_\omega F_{A(t+1)}-\lambda \Id|_{H_0}^2\frac{\omega^n}{n!}\rightarrow 0,
\end{equation}
as $t\rightarrow +\infty$.

For any $m\geq t_0>0$, there exists $t_m\in [m, m+1]$, such that $I_\epsilon(t_m)=\int_m^{m+1}I_\epsilon(t)dt$. From the formula (\ref{moFF}), it follows that
\begin{equation}
\begin{split}
I_\epsilon(t)\leq &e^{2C}I_\epsilon(t_m)=e^{2C}\int_{m}^{m+1}I_\epsilon(t)dt\\
=&\frac{e^{2C}}{2}\Big(\int_{\tilde{M}}|\sqrt{-1}\Lambda_\omega F_{A(m)}-\lambda \Id|_{H_0}^2\frac{\omega_\epsilon^n}{n!}-\int_{\tilde{M}}|\sqrt{-1}\Lambda_\omega F_{A(m+1)}-\lambda \Id|_{H_0}^2\frac{\omega_\epsilon^n}{n!}\Big),
\end{split}
\end{equation}
for any $t\in [m+1, m+2]$.
Applying Fatou's lemma again, we derive
\begin{equation}
\begin{split}
I(t)\leq& \lim_{\epsilon\rightarrow 0}I_\epsilon(t)\\
\leq & \lim_{\epsilon\rightarrow 0}\frac{e^{2C}}{2}\Big(\int_{\tilde{M}}|\sqrt{-1}\Lambda_{\omega_\epsilon} F_{A_\epsilon(m)}-\lambda \Id|_{\hat{H}}^2\frac{\omega_\epsilon^n}{n!}-\int_{\tilde{M}}|\sqrt{-1}\Lambda_{\omega_\epsilon} F_{A_\epsilon(m+1)}-\lambda \Id|_{\hat{H}}^2\frac{\omega_\epsilon^n}{n!}\Big)\\
=&\frac{e^{2C}}{2}\Big(\int_{M\setminus\Sigma_\mathcal{E}}|\sqrt{-1}\Lambda_\omega F_{A(m)}-\lambda \Id|_{\hat{H}}^2\frac{\omega^n}{n!}-\int_{M\setminus\Sigma_\mathcal{E}}|\sqrt{-1}\Lambda_\omega F_{A(m+1)}-\lambda \Id|_{\hat{H}}^2\frac{\omega^n}{n!}\Big),\\
\end{split}
\end{equation}
for any $t\in [m+1, m+2]$. This together with (\ref{moFF2}) means that $I(t)\rightarrow 0$, as $t\rightarrow +\infty$.

\hfill $\Box$ \\

\medskip

Now we recall  other Hermitian-Yang-Mills type functionals which are introduced in \cite{DW1}. For any $\verb"a"\in
\verb"u"(R)$, let $\varphi_{\alpha
}(\verb"a")=\sum_{j=1}^{R}|\lambda_{j}|^{\alpha }$, where $\verb"u"(R)$ is the Lie algebra of the unitary group
$U(R)$,
$\sqrt{-1}\lambda_{j}$ are the eigenvalues of $\verb"a"$, and $\alpha \geq 1$ is a real number.
 For a given real number
$N$, define the Hermitian-Yang-Mills type functionals as follows:
\begin{equation}
HYM_{\alpha , N}(A, M, \omega )=\frac{1}{\Vol(M, \omega )}\int_{M}\varphi_{\alpha }(\frac{1}{2\pi}\Lambda_{\omega }F_{A} -\sqrt{-1}N \Id_{E}) \frac{\omega^{n}}{n!}.
\end{equation}
Let $A_\epsilon(t)$ be the long time solution of the Yang-Mills flow (\ref{YM1}) on the Hermitian vector bundle $(E, \hat{H})$ over the K\"ahler manifold $(\tilde{M}, \omega_{\epsilon})$. For any smooth convex ad-invariant function $\varphi $, we have
\begin{equation}\label{2.33}
(\Delta_{\omega_{\epsilon}} -\frac{\partial }{\partial t})\varphi(\frac{1}{2\pi}\Lambda_{\omega_{\epsilon}} F_{A_\epsilon(t)} -\sqrt{-1}N \Id_{E})\geq 0,
\end{equation}
whose proof can be found in \cite{DW1} (Proposition 2.25). From \cite{AB} (Proposition 12.16), we know that $\varphi_{\alpha }$ is a convex function on $\verb"u"(R)$ and it can be approximated by a family of smooth convex ad-invariant functions $\varphi_{\alpha , \rho }$ as $\rho \rightarrow 0$. Integrating  (\ref{2.33}) gives that $t\mapsto HYM_{\alpha , N}(A_{\epsilon}(t), \tilde{M}, \omega_{\epsilon} )$ is
nonincreasing along the Yang-Mills flow, for any $0<\epsilon \leq 1$. Since $H_{ \epsilon}(x, t)$  converges  to the long time solution $H(x, t)$ of  the Hermitian-Yang-Mills flow (\ref{SSS1}) outside $\Sigma_{\mathcal{E}} $ in $C_{loc}^{\infty}$-topology as $\epsilon \rightarrow 0$, and $|\Lambda_{\omega } F_{A_{\epsilon}(t)}|_{\hat{H}}$ is uniformly bounded for any $0<\epsilon \leq 1$ and $0<t_{0}\leq t$, it is easy to see that $HYM_{\alpha , N}(A_{\epsilon}(t), \tilde{M}, \omega_{\epsilon} )\rightarrow HYM_{\alpha , N}(A(t), M, \omega )$ as $\epsilon \rightarrow 0$ and  $t\mapsto HYM_{\alpha , N}(A(t), M, \omega )$ is
also nonincreasing. So we obtain the following lemma.

\begin{lemma}\label{lemmaHYM}
Let $A(t)$ be the long time solution of the Yang-Mills flow (\ref{YM2}) on the Hermitian vector bundle $(\mathcal{E}|_{M\setminus \Sigma_{\mathcal{E}}}, \hat{H})$, then $t\mapsto HYM_{\alpha , N}(A(t), M, \omega )$ is
 nonincreasing.
\end{lemma}

Clearly Fatou's lemma tells us
\begin{equation}\label{interg}
\begin{split}
&4\pi^{2}\int_{M} (2c_{2}(\mathcal{E})-c_{1}(\mathcal{E})\wedge c_{1}(\mathcal{E}))\wedge\frac{\omega^{n-2}}{(n-2)!}\\
=&\lim_{\epsilon \rightarrow 0}4\pi^{2}\int_{\tilde{M}} (2c_{2}(E)-c_{1}(E)\wedge c_{1}(E))\wedge\frac{\omega_{\epsilon}^{n-2}}{(n-2)!}\\
=&\lim_{\epsilon \rightarrow 0}\int_{\tilde{M}}\tr (F_{A_{\epsilon}(t) }\wedge F_{A_{\epsilon}(t)})\wedge \frac{\omega_{\epsilon }^{n-2}}{(n-2)!}\\
=&\lim_{\epsilon \rightarrow 0}\int_{\tilde{M}}(|F_{A_{\epsilon}(t)}|_{\hat{H}, \omega_{\epsilon}}^{2}-|\Lambda_{\omega_{\epsilon}} F_{A_{\epsilon}(t)}|_{\hat{H}}^{2}) \frac{\omega_{\epsilon}^{n}}{n!}\\
\geq&  \int_{M\setminus \Sigma_\mathcal{E}}|F_{A(t)}|_{\hat{H}, \omega}^{2}\frac{\omega^{n}}{n!}-\int_{M\setminus \Sigma_\mathcal{E}}|\sqrt{-1}\Lambda_{\omega} F_{A(t)}|_{\hat{H}}^{2} \frac{\omega^{n}}{n!},\\
\end{split}
\end{equation}
and then it holds that
\begin{equation}\label{mono2}
\begin{split}
&\int_{M\setminus \Sigma_\mathcal{E}}|F_{A(t)}|_{\hat{H}, \omega}^{2}\frac{\omega^{n}}{n!}\\
\leq &\int_{M\setminus \Sigma_\mathcal{E}}|\sqrt{-1}\Lambda_{\omega} F_{A(t)}|_{\hat{H}}^{2} \frac{\omega^{n}}{n!}+4\pi^{2}\int_{M} (2c_{2}(\mathcal{E})-c_{1}(\mathcal{E})\wedge c_{1}(\mathcal{E}))\wedge\frac{\omega^{n-2}}{(n-2)!}\\
\leq &\int_{M\setminus \Sigma_\mathcal{E}}|\sqrt{-1}\Lambda_{\omega}(F_{A(t_0)})-\lambda \Id|_{\hat{H}}^{2} \frac{\omega^{n}}{n!}+\lambda^2 \rank E\int_{M\setminus \Sigma_\mathcal{E}}\frac{\omega^n}{n!}\\
&+4\pi^{2}\int_{M} (2c_{2}(\mathcal{E})-c_{1}(\mathcal{E})\wedge c_{1}(\mathcal{E}))\wedge\frac{\omega^{n-2}}{(n-2)!},
\end{split}
\end{equation}
for all $0< t_0 \leq t$.
For simplicity, in the sequel we set
\begin{equation}
\begin{split}
\mbox{HYM}(A(t_0))=&\int_{M\setminus \Sigma_\mathcal{E}}|\sqrt{-1}\Lambda_{\omega}(F_{A(t_0)})-\lambda \Id|_{\hat{H}}^{2} \frac{\omega^{n}}{n!}\\
&+\lambda^2 \rank E\int_{M\setminus \Sigma_\mathcal{E}}\frac{\omega^n}{n!}+4\pi^{2}\int_{M} (2c_{2}(\mathcal{E})-c_{1}(\mathcal{E})\wedge c_{1}(\mathcal{E}))\wedge\frac{\omega^{n-2}}{(n-2)!}.
\end{split}
\end{equation}
Let $f$ be a smooth function with support in $M\setminus \Sigma_{\mathcal{E}}$, we have
\begin{equation}\label{mono3}
 \begin{split}
& \frac{d}{dt}\Big(\int_{M} f ^{2}|F_{A(t)}|_{\hat{H}, \omega}^{2}\, \frac{\omega^{n}}{n!}\Big)\\
 =&-2\int_{M} f^{2}\Big(\Big|\frac{\partial A (t)}{\partial t
}\Big|^{2}\Big)\frac{\omega^{n}}{n!}+
2Re \int_{M}\Big\langle
\sqrt{-1}\Lambda_{\omega}((\overline{\partial
}-\partial )(f^{2})\wedge (F_{A(t)})) , \frac{d A(t)}{d t}\Big\rangle \, \frac{\omega^{n}}{n!}\\
&- 2Re \int_{M}\Big\langle ( \sqrt{-1}\Lambda_{\omega}F_{A(t)})(\overline{\partial
}-\partial )(f^{2}) , \frac{d A(t)}{d t}\Big\rangle \, \frac{\omega^{n}}{n!}.\\
\end{split}
\end{equation}
Integrating over $[s, \tau ]$ with
respect to $t$ on both sides of  (\ref{mono3}) and using the inequality (\ref{mono2}),
 we deduce the following local energy estimate.

\medskip

\begin{lemma}\label{localenergy}{\bf(Lemma 5 in \cite{HT})}
Let $A(t)$ be the long time solution of the Yang-Mills flow (\ref{YM2}) on the Hermitian vector bundle $(\mathcal{E}|_{M\setminus \Sigma_{\mathcal{E}}}, \hat{H})$. For any $x_0$ with $B_{2R}(x_0)\subset M\setminus\Sigma_{\mathcal{E}}$ and for any two finite numbers $s, \tau\geq t_0>0$, we have
\begin{equation}\label{locale1}
\begin{split}
&\int_{B_R(x_0)}|F_{A(t)}|_{\hat{H}, \omega}^{2}(\cdot, s)\frac{\omega^{n}}{n!}\\
\leq &\int_{B_{2R}(x_0)}|F_{A(t)}|_{\hat{H}, \omega}^{2}(\cdot, \tau)\frac{\omega^{n}}{n!}+2\int_{\min{\{s, \tau\}}}^{\max{\{s, \tau\}}}\int_{B_{2R}(x_0)}\Big(\Big|\frac{\partial A}{\partial t}\Big|^2+2\Big|\frac{\partial \phi}{\partial t}\Big|^2\Big)\frac{\omega^{n}}{n!}dt\\
&+C \Big(\frac{|s-\tau|}{R^2}\mbox{\rm{HYM}}(A(t_0))\int_{\min{\{s, \tau\}}}^{\max{\{s, \tau\}}}\int_{M\setminus\Sigma_\mathcal{E}}\Big(\Big|\frac{\partial A}{\partial t}\Big|^2+2\Big|\frac{\partial \phi}{\partial t}\Big|^2\Big)\frac{\omega^{n}}{n!}dt\Big)^{\frac{1}{2}},\\
\end{split}
\end{equation}
where $C$ is a uniform constant.
\end{lemma}

\section{The limit behaviour of the Yang-Mills flow}
\setcounter{equation}{0}

 In this section, we consider the limit behaviour of the Yang-Mills flow (\ref{YM2}) on the Hermitian bundle $(\mathcal{E}|_{M\setminus \Sigma_{\mathcal{E}}}, \hat{H})$. We first recall the monotonicity inequality and the $\varepsilon$-regularity theorem obtained by Hong and Tian in \cite{HT}.
  For a fixed point $u_0=(x_0,t_0)\in M\times \mathbb{R}_+$, denote
\begin{equation}
\begin{split}
 &T_r(x_0,t_0)=\left\{u=(x,t):t_0-4r^2<t<t_0-r^2, x\in X\right\},
 \\& P_r(u_0)=B_r(x_0)\times[t_0-r^2,t_0+r^2].
\end{split}
\end{equation}
 The fundamental solution of (backward) heat equation with singularity at $(z_0,t_0)\in C^{n}\times \mathbb{R}_+$ is
\begin{equation}
\tilde{G}_{(z_0,t_0)}(z,t)=\frac{1}{(4\pi(t_0-t))^n}\exp{\Big(-\frac{|z-z_0|^2}{4(t_0-t)}\Big)},\ \ \ (t<t_0).
\end{equation}
Denote the exponential map centered at $x_0$ on $(M, \omega)$ by $\exp_{x_0}$, and set
\begin{equation}
G_{u_0}(x, t)=\tilde{G}_{(0, t_0)}(\exp_{x_0}^{-1}(x), t).
\end{equation}

In the following, we denote $d_{x_0}=\min\{dist(x_0, \Sigma_{\mathcal{E}}), i(M)\}$, where $dist(x_0, \Sigma_{\mathcal{E}})$ is the distance from $x_0$ to the closed set $\Sigma_{\mathcal{E}}$, $i(M)$ is the injective radius of $(M, \omega)$. Let $\varphi_{x_0}\in C_0^\infty(B_{d_{x_0}}(x_0))$ be a cut-off function such that $\varphi_{x_0}\equiv 1$ on $B_{d_{x_0}/2}(x_0)$, $\varphi_{x_{0}}\equiv 0$ outside $B_{d_{x_0}}(x_0)$ and $|\nabla\varphi_{x_{0}}|\leq 4/d_{x_0}$. Let $A(t)$ be the long time solution of the Yang-Mills flow (\ref{YM2}) on the Hermitian vector bundle $(\mathcal{E}|_{M\setminus \Sigma_{\mathcal{E}}}, \hat{H})$ with initial value $\hat{A}$. Set
\begin{equation}
\Phi(r; A)= r^2\int_{T_r(x_0, t_0)}\varphi_{x_{0}}^2G_{u_0}|F_{A(t)}|_{\hat{H}, \omega }^{2}\frac{\omega^{n}}{n!}dt.
\end{equation}
The same argument in \cite{HT}, only replacing the energy inequality by the above inequality (\ref{mono2}) concludes the following monotonicity inequalities.

\begin{theorem}\label{mono111}{\bf(Theorem 2 and 2' in \cite{HT})}
Let $A(t)$ be the long time solution of the Yang-Mills flow (\ref{YM2}) with initial connection $\hat{A}$ on $(M\setminus\Sigma_{\mathcal{E}}, \mathcal{E}|_{M\setminus\Sigma_\mathcal{E}}, \omega)$. Then for any fixed $t^{\ast}>0$, $u_{0}=(x_0, t_0)\in M\setminus\Sigma_{\mathcal{E}}\times [t^{\ast}, T]$, and for $r_1$ and $r_2$ with $0< r_1\leq r_2<\min{\{d_{x_0}, \sqrt{t_0-t^{\ast}}/2\}}$, we have
\begin{equation}
\Phi(r_1;A)\leq C\exp(C(r_2-r_1))\Phi(r_2;A)+C(r_2^2-r_1^2)\mbox{\rm{HYM}}(A(t^{\ast})),
\end{equation}
where $C$ is a positive constant which depends only on $dist^{-1}(x_0, \Sigma_{\mathcal{E}})$ and the geometry of $(M, \omega)$. Furthermore, if  $R\leq d_{x_{0}}$ and $f_{x_{0}, R}\in C_0^\infty(B_{R}(x_{0}))$ is a cut-off function satisfying $ 0\leq f_{x_{0}, R}\leq1$, $f_{x_{0},R}\equiv 1$ on $B_{R/2}(x_0)$£¬, $|\nabla f_{x_{0}, R}| \leq 8/R$ on $B_{R}(x_0)\setminus B_{R/2}(x_0)$, then we have
\begin{equation}\label{mono}
\begin{split}
&r_1^2\int_{T_{r_{1}}(x_0,t_0)}|F_{A(t)}|_{\hat{H}, \omega }^{2}f_{x_{0}, R}^2G_{u_0}\frac{\omega^{n}}{n!}\,dt\\
\leq & C\exp(C(r_2-r_1))r_{2}^2\int_{T_{r_{2}}(x_0,t_0)}|F_{A(t)}|_{\hat{H}, \omega }^{2}f_{x_{0}, R}^2G_{u_0}\frac{\omega^{n}}{n!}\,dt+C(r_2^2-r_1^2)\mbox{\rm{HYM}}(A(t^{\ast}))\\
& +CR^{2-2n}\int_{P_R(x_0,t_0)}|F_{A(t)}|_{\hat{H}, \omega }^{2}\frac{\omega^{n}}{n!}\,dt ,\\
\end{split}
\end{equation}
for any $0<r_{1}\leq r_{2}\leq \min{\{R/2, \sqrt{t_0-t^{\ast}}/2\}}$,
where $C$ is a positive constant depending only on the geometry of $(M, \omega )$.
\end{theorem}

\medskip

Using the above monotonicity inequality (\ref{mono}),  Hong and Tian obtain the following $\varepsilon $-regularity theorem.

\medskip

\begin{theorem}\label{mono111}{\bf(Theorem 4' in \cite{HT})}
Let $A(t)$ be the long time solution of the Yang-Mills flow (\ref{YM2}) with initial connection $\hat{A}$ on $(M\setminus\Sigma_{\mathcal{E}}, \mathcal{E}|_{M\setminus\Sigma_{\mathcal{E}}}, \omega)$, and $t^\ast$ be a positive number.  There exist positive constants $\varepsilon_0, \delta_0<1/4$ such that for any $x_{0}\in M\setminus \Sigma_{\mathcal{E}} $, if it holds that
\begin{equation}
  R^{2-2n}\int_{P_{R}(x_{0},t_{0})}|F_{A(t)}|_{\hat{H}, \omega }^{2}\frac{\omega^{n}}{n!}\,dt\leq \varepsilon_0 ,
\end{equation}
where $0<R<\min\{d_{x_{0}}, \frac{\sqrt{t_0-t^{\ast}}}{2}\}$, then for any $ \delta\in(0,\, \delta_0)$, we have
\begin{equation}\label{k1}
\sup\limits_{P_{\delta R }(x_{0},t_{0})}|F_{A(t)}|_{\hat{H}, \omega }^{2} \leq C(\delta R)^{-4},
\end{equation}
where $C$ depends only on the geometry of $(M,\omega)$, $\delta_0^{-1}$ and $\mbox{\rm{HYM}}(A(t^{\ast}))$.
\end{theorem}

\medskip

 Using the above $\varepsilon$-regularity theorem, we can analyze the limiting behavior of the  Yang-Mills flow (\ref{YM2}) on $(M\setminus\Sigma_{\mathcal{E}}, \mathcal{E}|_{M\setminus\Sigma_{\mathcal{E}}}, \omega)$. We will modify Tian's argument (Proposition 3.1.2 in \cite{Tian}) and Hong-Tian's argument (Proposition 6 in \cite{HT}) to be suitable for the non-compact case.

\medskip

\begin{theorem}\label{thmlimit}
 Let $A(t)$ be the long time solution of the Yang-Mills flow (\ref{YM2}) with initial connection $\hat{A}$ on the Hermitian vector bundle $(\mathcal{E}|_{M\setminus\Sigma_{\mathcal{E}}}, \hat{H})$  over $(M\setminus\Sigma_{\mathcal{E}},  \omega)$.  Then for every sequence $t_{k}\rightarrow +\infty$, there exists a subsequence
$\{ t_{j}\}$ such that as $t_{j}\rightarrow +\infty$, $A(t_{j})
$ converges, modulo gauge transformations, to a
solution $A_{\infty}$ of the Yang-Mills
equation  on a Hermitian vector bundle $(E_{\infty} , H_{\infty})$  in $C^{\infty}_{loc}$-topology outside
$\Sigma\subset M$, where $\Sigma $ is a closed set of Hausdorff complex codimension at least $2$ and  $\Sigma_{\mathcal{E}} \subset \Sigma$.
\end{theorem}

\medskip

{\bf Proof. } By Proposition \ref{I}, we know that
$\int_{M}\mid\frac{\partial{A}}{\partial{t}}\mid_{\hat{H}}^{2}
\frac{\omega^{n}}{n!}\rightarrow 0$ as $t\rightarrow +\infty$, and then
\begin{equation}
\int_{t_{k}-a}^{t_{k}+a}\int_{M}\Big|\frac{\partial{A}}{\partial{t}}\Big|_{\hat{H}}^{2}
\frac{\omega^{n}}{n!}dt\rightarrow
0,
\end{equation}
as
$t_{k}\rightarrow +\infty$,
for any $a>0$.
Choosing $r_{0}$ small enough and
assuming that
\begin{equation}\label{ep1}
r_0^{4-2n}\int_{B_{r_{0}}(x_{0})}|F_{A(t_{k})}|^{2}_{\hat{H}, \omega}\frac{\omega^{n}}{n!}<
\varepsilon_{1},
\end{equation}
where $\varepsilon_{1}$ is determined later.
Using the local energy estimate (\ref{locale1}) gives us that $\forall  t_k-r_0^2\leq s \leq t_k+r_0^2$, it holds that
\begin{equation}
\begin{split}
&\int_{B_{r_0/2}(x_0)}|F_{A(s)}|_{\hat{H}, \omega }^{2}(\cdot)\frac{\omega^{n}}{n!}\\
\leq &\int_{B_{r_0}(x_0)}|F_{A(t_{k})}|_{\hat{H}, \omega }^{2}(\cdot)\frac{\omega^{n}}{n!}+2\int_{t_k-r_0^2}^{t_k+r_0^2}\int_{B_{r_0}(x_0)}\Big|\frac{\partial A}{\partial t}\Big|^2\frac{\omega^{n}}{n!}dt\\
&+C \Big(\mbox{HYM}(A(t_0))\int_{t_k-r_0^2}^{t_k+r_0^2}\int_{M}\Big|\frac{\partial A}{\partial t}\Big|^2\frac{\omega^{n}}{n!}dt\Big)^{\frac{1}{2}},
\end{split}
\end{equation}
and
\begin{equation}
\begin{split}
&\int_{t_k-\frac{r_0^2}{4}}^{t_k+\frac{r_0^2}{4}}\int_{B_{r_0/2}(x_0)}|F_{A(s)}|_{\hat{H}, \omega }^{2}(\cdot)\frac{\omega^{n}}{n!} ds\\
\leq &\frac{1}{2}r_0^2\Big(\int_{B_{r_0}(x_0)}|F_{A(t_{k})}|_{\hat{H}, \omega }^{2}(\cdot)\frac{\omega^{n}}{n!}+2\int_{t_k-r_0^2}^{t_k+r_0^2}\int_{B_{r_0}(x_0)}\Big|\frac{\partial A}{\partial t}\Big|^2\frac{\omega^{n}}{n!}dt\\
&+C \Big(\mbox{HYM}(A(t_0))\int_{t_k-r_0^2}^{t_k+r_0^2}\int_{M\setminus\Sigma}\Big|\frac{\partial A}{\partial t}\Big|^2\frac{\omega^{n}}{n!}dt\Big)^{\frac{1}{2}}\Big).
\end{split}
\end{equation}
Then (\ref{ep1}) implies that, for sufficiently large $t_{k}$,
\begin{equation}
\begin{split}
&\Big(\frac{r_0}{2}\Big)^{2-2n}\int_{P_{r_0/2}(x_0, t_k)}|F_{A(t)}|_{\hat{H}, \omega }^{2}(\cdot)\frac{\omega^{n}}{n!}dt\\
\leq &2^{2n-3}r_0^{4-2n}\int_{B_{r_0}(x_0)}|F_{A(t_{k})}|_{\hat{H}, \omega }^{2}(\cdot)\frac{\omega^{n}}{n!}dt\\
&+2^{2n-2}r_0^{4-2n}\int_{t_k-r_0^2}^{t_k+r_0^2}\int_{B_{r_0}(x_0)}\Big|\frac{\partial A}{\partial t}\Big|^2\frac{\omega^{n}}{n!}dt\\
&+2^{2n-3}Cr_0^{4-2n} \Big(\mbox{HYM}(A(t_0))\int_{t_k-r_0^2}^{t_k+r_0^2}\int_{M\setminus\Sigma}\Big|\frac{\partial A}{\partial t}\Big|^2\frac{\omega^{n}}{n!}dt\Big)^{\frac{1}{2}}\\
\leq &2^{2n-1}\varepsilon_1,
\end{split}
\end{equation}
where we choose $2^{2n-1}\varepsilon_{1}\leq \varepsilon_{0}$ and $\varepsilon_{0}$ is
the constant determined in Theorem \ref{mono111}. Therefore, we obtain
\begin{equation}\label{ep2}
\sup_{P_{\delta r_0/2}(x_0, t_k)}|F_{A(t)}|_{\hat{H}, \omega }^{2}\leq C(\delta r_0)^{-4}
\end{equation}
for any $\delta \in(0, \delta_0)$ and sufficiently large $k$, where $C$ is a uniform constant.

Applying (\ref{ep2}), (\ref{locale1}) and Moser's parabolic estimate to the following inequality
\begin{equation}\label{YMes2}
(\Delta -\frac{\partial }{\partial t})|F_{A(t)}|_{\hat{H}, \omega }^{2}\geq -C(|F_{A(t)}|_{\hat{H}, \omega } +|Rm(\omega)|_{\omega})|F_{A(t)}|_{\hat{H}, \omega }^{2},
\end{equation}
we derive for sufficiently large $k$,
\begin{equation}
\begin{split}
&\sup_{(x, t)\in P_{\frac{\delta_{0}r_0}{8}}(x_0, t_k)}|F_{A(t)}|_{\hat{H}, \omega }^{2}(x)\leq C (\delta_{0} r_{0})^{-4-2n}\int_{P_{\frac{\delta_{0}r_0}{4}}(x_0, t_k)}|F_{A(t)}|_{\hat{H}, \omega }^{2}(\cdot)\frac{\omega^{n}}{n!}dt\\
\leq &C (\delta_{0} r_{0})^{-2-2n} \sup_{t_{k}-(\frac{\delta_{0}r_0}{4})^{2}\leq t \leq t_{k}+(\frac{\delta_{0}r_0}{4})^{2}} \int_{B_{\frac{\delta_{0}r_0}{4}}(x_0)}|F_{A(t)}|_{\hat{H}, \omega }^{2}(\cdot)\frac{\omega^{n}}{n!}\\
\leq &C (\delta_{0} r_{0})^{-2-2n}  \int_{B_{\frac{\delta_{0}r_0}{2}}(x_0)}|F_{A(t)}|_{\hat{H}, \omega }^{2}(\cdot)\frac{\omega^{n}}{n!}+C (\delta_{0} r_{0})^{-2-2n} \varepsilon_{1}\\
\leq &C(\delta_{0}^{4-2n}(\delta_{0} r_{0})^{-6}+(\delta_{0} r_{0})^{-2-2n})\varepsilon_{1},
\end{split}
\end{equation}
and then
\begin{equation}\label{ep11}
\begin{split}
&(\delta r_0)^{4-2n}\int_{B_{\delta r_{0}}(x_{0})}|F_{A(t_{k})}|^{2}_{\hat{H}, \omega}\frac{\omega^{n}}{n!}\\
\leq & (\delta r_0)^{4-2n} \Vol (B_{\delta r_{0}}(x_{0}))C(\delta_{0}^{4-2n}(\delta_{0} r_{0})^{-6}+(\delta_{0} r_{0})^{-2-2n})\varepsilon_{1}\\
< & \varepsilon_{1},
\end{split}
\end{equation}
where $\delta \leq r_{0}^{n+1}$ and $r_{0}$ is small enough. Setting
\begin{equation}
r_{i}=(r_{i-1})^{n+2}=r_{0}^{(n+2)^{i}},
\end{equation}
and repeating the above argument, we know that (\ref{ep1}) implies
\begin{equation}\label{ep111}
( r_i)^{4-2n}\int_{B_{ r_{i}}(x_{0})}|F_{A(t_{k})}|^{2}_{\hat{H}, \omega}\frac{\omega^{n}}{n!}
<  \varepsilon_{1},
\end{equation}
for $i\geq 1$ and sufficiently large $t_k$.

We set
\begin{equation}
d_{x}=dist(x, \Sigma_{\mathcal{E}}), \quad  U_d=\{x\in M: d_{x}< d\},
\end{equation}
\begin{equation}
\hat{\Sigma}_{k, j, i}=\{x\in M\setminus U_{r_{j}}: r_{i}^{4-2n}\int_{B_{r_{i}}(x)}|F_{A(t_{k})}|_{\hat{H}, \omega }^{2}(\cdot)\frac{\omega^{n}}{n!}\geq \varepsilon_1 \},
\end{equation}
for any $k\geq 1$ and $i\geq j \geq 1$. By the standard diagonal process, we can choose a subsequence which also is denoted by $\{t_{k}\}$ such that for each $j\leq i$, $\hat{\Sigma}_{k, j, i}$
converges to a closed subset $\Sigma_{j, i}$ as $k\rightarrow +\infty$. From (\ref{ep111}), it is easy to see that $\Sigma_{j, i_{1}}\subset \Sigma_{j, i_{2}}$ for $i_{1}\geq i_{2}$.
Define
\begin{equation}
\Sigma_{j}=\bigcap_{i }\Sigma_{j, i}, \Sigma_{an}=\bigcup_{ j }\Sigma_j, \quad \Sigma=\Sigma_{\mathcal{E}}\bigcup \Sigma_{an}.
\end{equation}

\medskip

\begin{claim}\label{claim1}
$\Sigma$ is closed.
\end{claim}

\medskip

{\bf Proof. } Suppose $x_{0}\in M\setminus \Sigma $ and set $d_0=dist(x_0, \Sigma_{\mathcal{E}})>0$. For any $r_{j}<d_{0}$, we have $x_{0}\notin \Sigma_{j}$, $x_{0}\notin \Sigma_{j, i}$ for $i$ sufficiently large, and $x_{0}\notin \hat{\Sigma}_{k, j, i}$ for $t_{k}$ sufficiently large. Then
it follows that
\begin{equation}
\liminf_{k\rightarrow \infty}r_{i}^{4-2n}\int_{B_{r_{i}}(x_0)}|F_{A(t_{k})}|_{\hat{H}, \omega }^{2}\frac{\omega^{n}}{n!}< \varepsilon_1,
\end{equation}
 for $i$ and $k$ sufficiently large. Together with (\ref{ep2}),
fixing small $r_{i_{0}}$, for any $x\in B_{\frac{1}{8}\delta_{0} r_{i_{0}}}(x_0)\subset M\setminus \Sigma_{\mathcal{E}}$, we get
\begin{equation}\label{sm}
r^{4-2n}\int_{B_{r}(x)}|F_{A(t_{k})}|_{\hat{H}, \omega }^{2}\frac{\omega^{n}}{n!}<\varepsilon_{1},
\end{equation}
when $r$ is small enough and $k$ is large enough. Clearly $(\ref{sm})$ implies that, for  $k$ and $i$ sufficiently large,
$
B_{\frac{1}{8}\delta_{0} r_{i_{0}}}(x_0)\cap \{\Sigma_{\mathcal{E}}\cup \hat{\Sigma}_{k, j, i} \}=\emptyset ,
$
and
\begin{equation}
B_{\frac{1}{16}\delta_{0} r_{i_{0}}}(x_0)\cap \{\Sigma_{\mathcal{E}}\cup \Sigma_{ j, i} \}=\emptyset
\end{equation}
for all $j$.
Then $B_{\frac{1}{16}\delta_{0} r_{i_{0}}}(x_0)\cap \{\Sigma_{\mathcal{E}}\cup \Sigma_{ j} \}=\emptyset$  for all $j$, this means that  $B_{\frac{1}{16}\delta_{0} r_{i_0}}(x_0)\subset M\setminus \Sigma $ and concludes the proof of Claim \ref{claim1}.

\medskip

\begin{claim}\label{claim2}
The Hausdorff codimension of $\Sigma$ is at least 4.
\end{claim}

\medskip

{\bf Proof.} Since the sheaf $\mathcal{E}$ is torsion-free, it is well known that the Hausdorff codimension of $\Sigma_{\mathcal{E}}$ is at least 4 and the $(2n-4)$-dimensional Hausdorff measure is finite (i.e. $H^{2n-4}(\Sigma_{\mathcal{E}})< +\infty$). The definition says that
\begin{equation}
H^{2n-4}_\delta(\Sigma_{\mathcal{E}})=\inf\{\sum_{\alpha}(r_{\alpha})^{2n-4}|\cup_\alpha B_{x_\alpha}(r_\alpha)\supset \Sigma_{\mathcal{E}}, r_\alpha< \delta\},
\end{equation}
\begin{equation}
H^{2n-4}(\Sigma_{\mathcal{E}})=\lim_{\delta\rightarrow 0}H^{2n-4}_\delta(\Sigma_{\mathcal{E}}).
\end{equation}
Because $H^{2n-4}_\delta$ is monotonically nonincreasing with respect to $\delta$, $H^{2n-4}_\delta(\Sigma_{\mathcal{E}})\leq H^{2n-4}(\Sigma_{\mathcal{E}})$ for all $\delta$. Since $\Sigma_{\mathcal{E}}$ is compact,  for an arbitrary $\delta_0 > 0$, there exists a finite cover of $\Sigma_{\mathcal{E}}$, $\{B_{R_\alpha}(x_\alpha )\}$, such that $\sum\limits_\alpha(R_\alpha)^{2n-4}< H^{2n-4}(\Sigma_{\mathcal{E}})+1$ and $R_\alpha < \delta_0$. Then we can find a positive number $\delta'< \delta_0$ such that $M\setminus \cup_\alpha B_{R_\alpha}(x_\alpha )\subset M\setminus U_{\delta'}$. So it follows that $\Sigma\cap (M\setminus \cup_\alpha B_{R_\alpha}(x_\alpha ))\subset \Sigma\cap (M\setminus U_{\delta'})$ and $\Sigma\cap (M\setminus U_{\delta'})$ is closed.
Set
\begin{equation}
\widetilde{\Sigma}_{\delta'}= \Sigma\cap (M\setminus U_{\delta'})= \Sigma_{an}\cap (M\setminus U_{\delta'}).
\end{equation}
Suppose that $r_{\tilde{i}}\leq \frac{1}{8}\delta ' < r_{\tilde{i}-1}$.  Let $r=r_{\tilde{i}} >0$, we can find a finite collection of geodesic balls $\{B_{4r}(x_{\beta })\}$ such that $\{B_{4r}(x_{\beta})\}$ is a cover of $\widetilde{\Sigma}_{\delta'}$, $x_{\beta}\in \widetilde{\Sigma}_{\delta'}$ for all $i$, and $B_{2r}(x_{\beta_{1}})\cap B_{2r}(x_{\beta_{2}})=\emptyset$ for $\beta_{1}\neq \beta_{2}$. For every point $x_{\beta}$, suppose that $x_{\beta}\in \Sigma_{j}$ and take $i$ large enough such that $r_{i}<r$, then for $k$ sufficiently large, there are $y_{\beta}\in \hat{\Sigma}_{k, j, i}$ such that $dist (x_{\beta}, y_{\beta})<r$. It is easy to see that $\{B_{5r}(y_{\beta })\}$ is a finite covering of $\widetilde{\Sigma}_{\delta'}$ and $B_{r}(y_{\beta_{1}})\cap B_{r}(y_{\beta_{2}})=\emptyset$ for $\beta_{1}\neq \beta_{2}$.

Choosing sufficiently large $k$, $y_{\beta}\in \hat{\Sigma}_{k, j, i}$, we know
\begin{equation}
 (r_{i})^{4-2n}\int_{B_{r_{i}}(y_{\beta })}|F_{A(t_{k})}|_{\hat{H}, \omega }^{2}\frac{\omega^{n}}{n!}\geq \varepsilon_{1},
\end{equation}
for every $\beta $.
Repeating the argument in the proof of (\ref{ep111}) yields
\begin{equation}
 (r)^{4-2n}\int_{B_{r}(y_{\beta })}|F_{A(t_{k})}|_{\hat{H}, \omega }^{2}\frac{\omega^{n}}{n!}\geq \varepsilon_{1},
\end{equation}
for every $\beta $.
Summing over $\beta $ and using the inequality (\ref{mono2}), we get
\begin{equation}
\begin{split}
\sum_{\beta} \varepsilon_{1}< &\sum_{\beta} (r)^{4-2n}\int_{B_{r}(y_{\beta })}|F_{A(t_{k})}|_{\hat{H}, \omega }^{2}\frac{\omega^{n}}{n!}\\
< &(r)^{4-2n}\int_{M}|F_{A(t_{k})}|_{\hat{H}, \omega }^{2}\frac{\omega^{n}}{n!}\\
\leq &(r)^{4-2n}\mbox{HYM}(A(t_1)),
\end{split}
\end{equation}
and then
\begin{equation}
\sum_\alpha R_\alpha^{2n-4}+\sum_{\beta }(5r)^{2n-4}< H^{2n-4}(\Sigma_{\mathcal{E}})+1+\frac{1}{\varepsilon_1}5^{2n-4}\mbox{HYM}(A(t_1)).
\end{equation}
It implies that
\begin{equation}
H^{2n-4}_{\delta_0}(\Sigma)< H^{2n-4}(\Sigma_{\mathcal{E}})+1+\frac{1}{\varepsilon_1}5^{2n-4}\mbox{HYM}(A(t_1)).
\end{equation}
Letting $\delta_0\rightarrow 0$, we obtain $H^{2n-4}(\Sigma)$ is finite. This concludes the proof of Claim \ref{claim2}.

\medskip

 Given a compact subset $\Omega \subset M\setminus \Sigma $, we  suppose $\Omega \subset M\setminus U_{d_{0}} $ for some $d_{0}>0$. For any point $x_{0}\in \Omega$, as that in the proof of  Claim \ref{claim2}, we know that, there is $r_{i_{0}}$ such that
 \begin{equation}\label{ep222}
\sup_{B_{\frac{\delta_{0} r_{i_0}}{8}}(x_0)}|F_{A(t_{k})}|_{\hat{H}, \omega }^{2}\leq C(\delta_{0} r_{i_{0}})^{-4}
\end{equation}
 for $t_{k}$ sufficiently large. Since $\Omega$ is compact, we can cover it by a finite union of balls such that every ball satisfies the above estimate (\ref{ep222}). So it follows that $\sup_{\Omega}|F_{A(t_{k})}|_{\hat{H}, \omega }^{2}$ is uniformly bounded. Uhlenbeck's Theorem (Theorem 3.6 in \cite{U2}) implies that there exists a subsequence  of $\{A(t_{k})\}
$,  modulo gauge transformations,
 converging to a connection $A_{\infty}$
weakly in $L_{1,loc}^{2}$-topology outside
$\Sigma$, where $A_{\infty}$ is a solution of the
Yang-Mills equation on a Hermitian vector bundle $(E_{\infty}, H_{\infty})$ which is isometric to $(\mathcal{E}|_{M\setminus\Sigma_{\mathcal{E}}}, \hat{H})$  outsides $\Sigma$. Furthermore, by standard
parabolic regularity techniques and using Hong-Tian's argument (Proposition 6 in \cite{HT}), we know that $A(t_{k})$
converges to $A_{\infty}$ in $C_{loc}^{\infty
}$-topology outside $\Sigma$. This concludes the proof of Theorem \ref{thmlimit}.

\hfill $\Box$ \\

\medskip

 From the estimates (\ref{H0013}) and (\ref{H00013}), we see that $|\theta (A(t), \omega )|_{\hat{H}}$ is uniformly bounded for $t\geq t_{0}>0$. Through the same argument as that in Corollary 2.12 in \cite{DW1} (or Corollary 3.12 in \cite{LZ1}), we have the following  corollary.

\medskip

\begin{corollary}\label{Corollary 2.5.}  Let $A(t_{k})$ be a sequence
of connections along the Yang-Mills flow (\ref{YM2}) with the
limit $A_{\infty } $, then:

(1) $\theta (A(t_{k}) , \omega)\rightarrow \theta (A_{\infty} ,
\omega)$ strongly in $L^{p}$ as $k\rightarrow +\infty$ for all $1\leq p<\infty$,  and consequently
\begin{equation}\lim_{t\rightarrow +\infty }\int_{M}|\theta (A(t) ,
\omega)|_{\hat{H}}^{2}\frac{\omega ^{n}}{n!}=\int_{M}|\theta (A_{\infty} , \omega)|_{H_{\infty}}^{2}\frac{\omega ^{n}}{n!};\end{equation}

(2) $\|\theta (A_{\infty} , \omega)\|_{L^{\infty}}\leq
\|\theta (A(t_{k}) ,
\omega)\|_{L^{\infty}}\leq \|\theta (A(t_{0}) ,
\omega)\|_{L^{\infty}}$ for $0< t_{0}\leq t_{k}$.

\end{corollary}

\medskip

In the sequel, we call $A_{\infty}$ an Uhlenbeck limit of $A(t)$. Since $A_{\infty}$ is a solution of the Yang-Mills equation, i.e. it satisfies
\begin{equation}
D_{A_{\infty}}^{\ast}F_{A_{\infty}}=0,
\end{equation}
by the K\"ahler indentity, we have
\begin{equation}
D_{A_{\infty}}\theta(A_{\infty}, \omega ) =0,
\end{equation}
i.e. $\theta(A_{\infty}, \omega ) $ is parallel. On the other  hand, $(\sqrt{-1}\theta(A_{\infty}, \omega ))^{\ast H_{\infty}
}=\sqrt{-1}\theta(A_{\infty}, \omega ) $, we can decompose $E_{\infty}$ according to the
eigenvalues of $\sqrt{-1}\theta(A_{\infty}, \omega )$ and obtain a holomorphic
orthogonal decomposition:
$
E_{\infty}=\bigoplus_{i=1}^{l}E_{\infty}^{i}
$
on $M\setminus \Sigma$. Let $\lambda_{i}$ be the eigenvalues of $\sqrt{-1}\theta(A_{\infty}, \omega )$, $H_{\infty}^{i}$ be the restrictions of $H_{\infty}$ to $E_{\infty}^{i}$ and $A_{\infty}^{i}=A_{\infty}|_{E^{i}}$, it is easy to see that $A_{\infty}^{i} $ is  a Hermitian-Einstein connection on
$(E_{\infty}^{i} , H_{\infty}^{i})$, i.e. it satisfies
\begin{equation}
\sqrt{-1}\Lambda_{\omega }F_{A_{\infty}^{i}}= \lambda_{i}\Id_{E_{\infty}^{i}}.
\end{equation}
Of course (\ref{mono2}) means that
\begin{equation}
\int_{M\setminus \Sigma } |F_{A_{\infty}}|^{2}_{H_{\infty}} \frac{\omega^{n}}{n!}\leq C <\infty .
\end{equation}
  Since the singularity set $\Sigma $ is of Hausdorff codimension at least $4$, by Theorem 2 in Bando and Siu's paper \cite{BS}, we know that  every $(E_{\infty}^{i}, \overline{\partial }_{A_{\infty}^{i}})$ can be extended to the whole $M$ as a reflexive sheaf (which is also denoted by $(E_{\infty}^{i}, \overline{\partial }_{A_{\infty}^{i}})$ for simplicity),  and $H_{\infty}^{i}$ can be smoothly extended over the place where the sheaf $(E_{\infty}^{i}, \overline{\partial }_{A_{\infty}^{i}})$ is locally free. Therefore, we deduce the following proposition.

  \medskip

\begin{proposition}\label{part1} The limiting $(E_{\infty}, \overline{\partial }_{A_{\infty}})$ can be extended to the whole $M$ as a reflexive sheaf  with a holomorphic orthogonal splitting
  \begin{equation}
(E_{\infty},  \overline{\partial}_{A_{\infty}},H_{\infty})=\bigoplus_{i=1}^{l}(E_{\infty}^{i} ,  \overline{\partial}_{A_{\infty}^{i}} , H_{\infty}^{i}
),
\end{equation}
and $H_{\infty}^{i}$ is an admissible Hermitian-Einstein metric on the reflexive sheaf $(E_{\infty}^{i} , \overline{\partial}_{A_{\infty}^{i}})$ for any $1\leq i\leq l$.
\end{proposition}

\medskip

\vskip 1mm

\section{ $L^{p}$-approximate critical Hermitian metric }
\setcounter{equation}{0}

In this section, we first recall the
Harder-Narasimhan-Seshadri filtration of reflexive sheaves (\cite{Ko2},
v.7.15, 7.17, 7.18; or \cite{Br}, section 7).  Then we  prove the existence of $L^{p}$-approximate critical Hermitian metric.  We will modify Daskalopoulos and Wentworth's  cut-off argument (\cite{DW1}) and Sibley's trick (\cite{Sib}) to be suitable for the reflexive sheaf case.

Let $\mathcal{E}$ be a reflexive sheaf over a compact K\"ahler manifold $(M, \omega)$.
If $\mathcal{E}$ is not $\omega$-stable,  there is a filtration of $\mathcal{E}$ by  coherent
 sub-sheaves
\begin{equation}
0=\mathcal{E}_{0}\subset \mathcal{E}_{1}\subset \cdots \subset \mathcal{E}_{\tilde{k}}=\mathcal{E} ,
\end{equation}
such that the quotients $\mathcal{Q}_{j}=\mathcal{E}_{j}/\mathcal{E}_{j-1}$ are torsion-free, $\omega$-semi-stable and $\mu_{\omega} (\mathcal{Q}_{j})>\mu_{\omega} (\mathcal{Q}_{j+1})$. We call it the
Harder-Narasimhan filtration (abbr. HN-filtration) of $\mathcal{E}$. The
associated graded sheaf $Gr^{HN}(\mathcal{E}
)=\oplus_{j=1}^{\tilde{k}}\mathcal{Q}_{j}$ is uniquely determined by the isomorphism
class of $\mathcal{E}$ and the K\"ahler class $[\omega]$.

\medskip

\begin{definition} For a reflexive sheaf $\mathcal{E}$ of rank $R$ over a compact K\"ahler manifold $(M, \omega )$, construct a nonincreasing
$R$-tuple of numbers
\begin{equation}
\vec{\mu }(\mathcal{E} )=(\mu_{1, \omega} , \cdots , \mu_{R, \omega})
\end{equation}
from the HN-filtration by setting: $\mu_{i, \omega }=\mu_{\omega} (\mathcal{Q}_{j})$, for
$\rank(\mathcal{E}_{j-1})+1\leq i\leq \rank(\mathcal{E}_{j})$. We call $\vec{\mu }(\mathcal{E})$ the Harder-Narasimhan type of $\mathcal{E}$.
\end{definition}

\hspace{0.3cm}

{\bf Remark: } {\it For a pair $\vec{\mu }$, $\vec{\lambda  }$ of
$R$-tuple's satisfying $\mu_1\geq \cdots \geq \mu_R$, $\lambda_1\geq \cdots \geq \lambda_R$, and
$\sum_{i=1}^{R}\mu_{i}=\sum_{i=1}^{R}\lambda_{i}$, we define:
\begin{equation}
\vec{\mu }\leq \vec{\lambda  } \quad \Leftrightarrow \quad
\sum_{i\leq k}\mu_{i}\leq \sum_{i\leq k}\lambda_{i}, \quad for\ all\ k=1, \cdots , R.
\end{equation}

 }

\medskip

Moreover, for every $\omega$-semistable quotient sheaf $\mathcal{Q}_{j}$, there is a further filtration, which is called by the Seshadri filtration, by subsheaves
\begin{equation} 0=\mathcal{E}_{j, 0}\subset
\mathcal{E}_{j, 1}\subset \cdots \subset \mathcal{E}_{j, k_{j}}=\mathcal{Q}_{j} ,\end{equation} such that the quotients $\mathcal{Q}_{j,\alpha}=\mathcal{E}_{j, \alpha}/\mathcal{E}_{j, \alpha-1}$ are
torsion-free and $\omega$-stable, $\mu_{\omega} (\mathcal{Q}_{j, \alpha})=\mu_{\omega}
(\mathcal{Q}_{j})$ for each $\alpha $. We call this double filtration $\{\mathcal{E}_{j, \alpha}\}$  the Harder-Narasimhan-Seshadri filtration (abbr. HNS-filtration) of the sheaf $\mathcal{E}$. The
associated graded sheaf: $ Gr^{HNS}(\mathcal{E} )=\oplus_{j=1}^{\tilde{k}}\oplus_{\alpha=1}^{k_{j}}\mathcal{Q}_{j , \alpha}
$ is uniquely determined by the isomorphism class of $\mathcal{E}$ and the K\"ahler class $[\omega]$.

\medskip

In the following, we denote
the Harder-Narasimhan-Seshadri filtration (or HNS-filtration) of $\mathcal{E}$ simply by:
\begin{equation}
0= \mathcal{E}_0\subset \mathcal{E}_1\subset \cdots \subset\mathcal{E}_{l-1}\subset\mathcal{E}_l=\mathcal{E},
\end{equation}
 where each $\mathcal{E}_i$ is a saturated subsheaf of $\mathcal{E}$. Set
\begin{equation}
\Sigma_{HNS}=\cup_{i=1}^{l}(\Sigma_{\mathcal{Q}_{i}}\cup \Sigma_{\mathcal{E}}),
\end{equation}
and refer to it as the singularity set of the HNS-filtration, where $\mathcal{Q}_{i}=\mathcal{E}_{i}/\mathcal{E}_{i-1}$ for each $1\leq i\leq l$. Since every $\mathcal{Q}_{i}$ is torsion-free, it is well known that $\Sigma_{HNS}$ is a complex analytic subset of complex codimension at least two.

 \medskip

\medskip

By Hironaka's flattening theorem (\cite{Hi2} or \cite{BS}), there is a finite sequence of blowing ups along compact sub-manifolds such that, if we denote by $\pi: \tilde{M}\rightarrow M$ the composition of all the blowing ups, then $\pi^{\ast}\mathcal{E}/tor{(\pi^{\ast
}\mathcal{E})}$ is locally free.

\medskip

\begin{proposition}\label{claim}
Let $E=\pi^{\ast}\mathcal{E}/tor{(\pi^{\ast}\mathcal{E})}$, then we can get a filtration  $\widetilde{\mathcal{F}}=\{\widetilde{\mathcal{E}_{i}}\}_{i=1}^{l}$ of $E$ from the HNS-filtration of $\mathcal{E}$:
\begin{equation}\label{HNS02}
0= \widetilde{\mathcal{E}_0}\subset \widetilde{\mathcal{E}_1}\subset \cdots \subset\widetilde{\mathcal{E}}_{l-1}\subset\widetilde{\mathcal{E}_l}=E,\quad \rm{over}\ \tilde{M}
\end{equation}
such that, for every $1\leq i \leq l$, $\widetilde{\mathcal{E}_i}$ is a reflexive sheaf, $\widetilde{\mathcal{Q}}_i=\widetilde{\mathcal{E}_i}/\widetilde{\mathcal{E}}_{i-1}$ is torsion free and isomorphic to the sheaf $\mathcal{Q}_{i} $  outside $\pi^{-1}(\Sigma_{HNS})$. Furthermore, every quotient sheaf $\widetilde{\mathcal{Q}}_{i}$ in the filtration (\ref{HNS02}) is {\bf $\omega_{\epsilon}$-stable}  for any $0<|\epsilon |\leq \epsilon^{\ast}\ll1$, and $\lim_{\epsilon\rightarrow 0}\deg_{\omega_{\epsilon}}(\widetilde{\mathcal{Q}}_{i})=\deg_{\omega}(\mathcal{Q}_{i})$.
\end{proposition}

\medskip

{\bf Proof.}
Pulling back the following exact sequences:
\[
\begin{aligned}
&0\rightarrow \mathcal{E}_i\rightarrow \mathcal{E}\rightarrow \mathcal{G}_i=\mathcal{E}/\mathcal{E}_i\rightarrow 0, \qquad &\textrm{over} \ M,\\
&0\rightarrow \mathcal{E}_{i-1}\rightarrow \mathcal{E}_i\rightarrow \mathcal{Q}_i=\mathcal{E}_i/\mathcal{E}_{i-1}\rightarrow 0,\qquad &\textrm{over} \ M,
\end{aligned}
\]
we get
\[
\begin{aligned}
&\pi^{\ast}\mathcal{E}_i\xrightarrow{f_i}\pi^{\ast}\mathcal{E}\rightarrow \pi^{\ast}\mathcal{G}_i\rightarrow 0, \qquad &\textrm{over} \ \tilde{M},\\
&\pi^{\ast}\mathcal{E}_{i-1}\xrightarrow{f_{i-1}}\pi^{\ast}\mathcal{E}\rightarrow \pi^{\ast}\mathcal{G}_{i-1}\rightarrow 0, \qquad &\textrm{over} \ \tilde{M},\\
&\pi^{\ast}\mathcal{E}_{i-1}\xrightarrow{g_i}\pi^{\ast}\mathcal{E}_i\rightarrow \pi^{\ast}\mathcal{Q}_i\rightarrow 0, \qquad &\textrm{over} \ \tilde{M}.
\end{aligned}
\]

Set $T_i=tor(Imf_i)=tor(Im\{\pi^{\ast}\mathcal{E}_i\rightarrow\pi^{\ast}\mathcal{E}\})$, $T=tor(\pi^{\ast}\mathcal{E})$. Then we can obtain the following exact sequence:
\[
0\rightarrow Imf_i/T_i\rightarrow \pi^{\ast}\mathcal{E}/T\rightarrow\check{\mathcal{G}}_i\triangleq \frac{\pi^{\ast}\mathcal{E}/T}{Imf_i/T_i}\rightarrow 0, \qquad \textrm{over}\ \tilde{M}.
\]

Set $E=\pi^{\ast}\mathcal{E}/T$, $\widetilde{\mathcal{E}}_i=Sat_E(Imf_i/T_i)$ and then $\widetilde{\mathcal{E}}_i$ is reflexive.
Clearly the definition gives the following exact sequences:
\[
\begin{aligned}
&0\rightarrow \widetilde{\mathcal{E}}_i\rightarrow E\rightarrow \widetilde{\mathcal{G}}_i\triangleq \check{\mathcal{G}}_i/tor(\check{\mathcal{G}}_i)\rightarrow 0,\qquad &\textrm{over} \ \tilde{M},\\
&0\rightarrow \widetilde{\mathcal{E}}_{i-1}\rightarrow E\rightarrow \widetilde{\mathcal{G}}_{i-1}\triangleq \check{\mathcal{G}}_{i-1}/tor(\check{\mathcal{G}}_{i-1})\rightarrow 0,\qquad &\textrm{over} \ \tilde{M}.
\end{aligned}
\]

Consider
\[
\xymatrix{
\pi^{\ast}\mathcal{E}_{i-1} \ar[dr]_{f_{i-1}}\ar[r]^{g_i} &\pi^{\ast}\mathcal{E}_i \ar[d]^{f_i} \ar[r] &\pi^{\ast}(\mathcal{E}_i/\mathcal{E}_{i-1}) \ar[r]& 0\\
&\pi^{\ast}\mathcal{E} \ar[d] \ar[dr] \\
&\pi^{\ast}(\mathcal{E}/\mathcal{E}_{i})\ar[d] &\pi^{\ast}(\mathcal{E}/\mathcal{E}_{i-1}) \ar[dr]\\
&0 &&0}.
\]
Of course $Imf_{i-1}= Im(f_i\circ g_i)$ means that $Imf_{i-1}$ is a subsheaf of $Imf_{i}$. Hence the following commutative diagram holds (all the horizontal sequences are exact):
\[
\xymatrix{
0\ar[r]&T_{i-1}\ar[r]\ar[d]&Imf_{i-1}\ar[r]\ar[d]&Imf_{i-1}/T_{i-1}\ar[r]\ar[d]^{h_i}&0,\\
0\ar[r]&T_{i}\ar[r]&Imf_{i}\ar[r]&Imf_{i}/T_{i}\ar[r]&0,}
\]
where we define the map $h_i$ by the commutation (it is easy to check that $h_i$ is well-defined). Moreover, a simple diagram shows that $h_i$ is injective. Noting that
$\widetilde{\mathcal{E}}_{i-1}=Ker\{E\rightarrow \check{\mathcal{G}}_{i-1}/tor(\check{\mathcal{G}}_{i-1})\}$ and $\widetilde{\mathcal{E}_{i}}=Ker\{E\rightarrow \check{\mathcal{G}}_i/tor(\check{\mathcal{G}}_i)\}$, considering the following sequences:
\[
\begin{aligned}
&E\rightarrow\check{\mathcal{G}}_{i-1}\rightarrow\check{\mathcal{G}}_{i-1}/tor(\check{\mathcal{G}}_{i-1})\rightarrow 0,\qquad &\textrm{over} \ \tilde{M},\\
&E\rightarrow\check{\mathcal{G}}_{i}\rightarrow\check{\mathcal{G}}_{i}/tor(\check{\mathcal{G}}_{i})\rightarrow 0,\qquad &\textrm{over} \ \tilde{M},
\end{aligned}
\]
we can see $\widetilde{\mathcal{E}}_{i-1}$ is a subsheaf of $\widetilde{\mathcal{E}}_{i}$.

Consider the following commutative diagram (all the vertical and horizontal sequences are exact):
\[
\xymatrix{
&&&0\ar[d]\\
&0\ar[d]&0\ar[d]&Ker(p_i)\ar[d]\\
0\ar[r]&\widetilde{\mathcal{E}}_{i-1}\ar[r]\ar[d]&E\ar[r]\ar@{=}[d]&\widetilde{\mathcal{G}}_{i-1}\ar[r]\ar[d]^{p_i}&0\\
0\ar[r]&\widetilde{\mathcal{E}}_{i}\ar[r]\ar[d]&E\ar[r]\ar[d]&\widetilde{\mathcal{G}}_{i}\ar[r]\ar[d]&0\\
&\widetilde{\mathcal{E}}_{i}/\widetilde{\mathcal{E}}_{i-1}\ar[d]&0&0\\
&0}
\]
where we define the map $p_i$ by the commutation and it is easy to check $p_i$ is well-defined (moreover, $p_i$ is surjective). The snake lemma tells us $\widetilde{\mathcal{E}}_{i}/\widetilde{\mathcal{E}}_{i-1}\cong Ker(p_i)$ and then $\widetilde{\mathcal{Q}}_i$ is torsion free (because $\widetilde{\mathcal{G}}_{i-1}$ is torsion free).
By a similar argument as that in Theorem 4.9, Proposition 4.10, Proposition 4.12 in \cite{Sib}, it is easy to see that $\lim_{\epsilon\rightarrow 0}\deg_{\omega_{\epsilon}}(\widetilde{\mathcal{Q}}_{i})=\deg_{\omega}(\mathcal{Q}_{i})$, and $\widetilde{\mathcal{Q}}_{i}$ in the filtration (\ref{HNS02}) is  $\omega_{\epsilon}$-stable  for any $0<|\epsilon |\leq \epsilon^{\ast}\ll1$ and $1\leq i \leq l$.

\hfill $\Box$ \\

\medskip

Let $H$ be a smooth Hermitian metric on the holomorphic bundle
$E$, and $\mathcal{F}=\{F_{i}\}_{i=1}^{l}$ be a
filtration of $E$ by saturated subsheaves: \[0=F_{0}\subset
F_{1}\subset \cdots \subset F_{l-1}\subset F_{l}=E.\]
For each $F_{i}$ and the metric $H$,  we have the associated unitary
projection $\pi_{i}^{H}:E\rightarrow E$ onto $F_{i}$, where
$\pi_{i}^{H}$ is an $L_{1}^{2}$-bounded  Hermitian endomorphism. For
convenience, set $\pi_{0}^{H}=0$. Given real numbers $\mu_{1}
,\cdots , \mu_{l}$ and a filtration $\mathcal{F}$, we define an $L_{1}^{2}$-bounded Hermitian endomorphism of $E$ by
\begin{equation}\Psi (\mathcal{F }, (\mu_{1} ,\cdots , \mu_{l}), H
)=\Sigma_{i=1}^{l}\mu_{i}(\pi_{i}^{H}-\pi_{i-1}^{H}).\end{equation} The
Harder-Narasimhan projection $\Psi_{\omega }^{HN}(E,
 , H)$ is the $L_{1}^{2}$-bounded Hermitian endomorphism
defined above in the particular case where $\mathcal{F}$ is the HN-filtration $\mathcal{F}^{HN}=\{F_{i}^{HN}(E)\}_{i=1}^{l}$ and $\mu_{i}=\mu_{\omega}
(F_{i}^{HN}(E)/F_{i-1}^{HN}(E))$.

\hspace{0.3cm}

\begin{definition}\label{def3.6}
Fix $\delta >0$ and $1\leq p\leq
\infty$. An $L^{p}$-$\delta $-approximate critical Hermitian
metric on a holomorphic bundle $E$ over a compact K\"ahler manifold $(M, \omega )$ is a smooth metric $H$
such that
$$\|\frac{\sqrt{-1}}{2\pi}\Lambda_{\omega }(F_{A_{H}})-\Psi_\omega ^{HN}(E,  H)\|_{L^{p}(\omega )}\leq
\delta ,
$$where $A_{H}$ is the Chern connection determined by
$(\overline{\partial }_{E} , H)$.
\end{definition}

\medskip

\medskip

Now recall the following lemma which was proved by Sibley in \cite{Sib} (Lemma 5.3).

\medskip

\begin{lemma}\label{volume} Let $(M, \omega )$ be a compact K\"ahler manifold of complex dimension $n$, and  $\pi :
\overline{M}\rightarrow M$ be a blow-up along a smooth complex sub-manifold $\Sigma$ of complex co-dimension $k$ where $k\geq 2$. Let $\eta $ be a K\"ahler metric on $\overline{M}$, and consider the family of
K\"ahler metrics $\omega_{\epsilon }=\pi^{\ast
}\omega +\epsilon \eta $, where $0< \epsilon\leq 1$. Then for any $0\leq \gamma < \frac{1}{k-1}$, we have $\frac{\eta ^{n}}{ \omega_{\epsilon}^{n}}\in L^{\gamma}(\overline{M}, \eta )$, and the $L^{\gamma}(\overline{M}, \eta )$-norm of $\frac{\eta ^{n}}{ \omega_{\epsilon}^{n}}$ is uniformly bounded in $\epsilon$, i.e. there is a positive constant $C^{\ast}$ such that
\begin{equation}\label{Int1}
\int_{\overline{M}}\Big( \frac{\eta ^{n}}{ \omega_{\epsilon}^{n}}\Big)^{\gamma } \frac{\eta ^{n}}{n!}\leq C^{\ast}
\end{equation}
for all $\epsilon$.
\end{lemma}

\medskip

Fixing  a number $\tilde{\gamma }<\frac{1}{n-1}$, for $0<\gamma \ll \tilde{\gamma }$, using the H\"older inequality, we have
\begin{equation}\label{3.9}
\begin{split}
&\int_{\tilde{M}}\Big(\frac{\eta_{k}^n}{\omega_{k, \epsilon }^n}\Big)^{\gamma }\frac{\eta_{k}^n}{n!}\\
\leq &\int_{\tilde{M}}\Big(\frac{\eta_{k}^n}{(\pi_{k}^{\ast}\eta_{k-1} +\epsilon_{k}\eta_{k}) ^n}\Big)^{1+\gamma }\Big(\frac{(\pi_{k}^{\ast}\eta_{k-1} +\epsilon_{k}\eta_{k}) ^n}{\omega_{k, \epsilon }^n}\Big)^{\gamma }\frac{(\pi_{k}^{\ast}\eta_{k-1} +\epsilon_{k}\eta_{k}) ^n}{n!}\\
\leq &\Big( \int_{\tilde{M}}\Big(\frac{\eta_{k}^n}{\pi_{k}^{\ast}\eta_{k-1} +\epsilon_{k}\eta_{k}}\Big)^{\gamma }\frac{\eta_{k}^n}{n!}  \Big) ^{\frac{1+\gamma }{1+\tilde{\gamma}}}\Big(\int_{\tilde{M}}\Big(\frac{(\pi_{k}^{\ast}\eta_{k-1} +\epsilon_{k}\eta_{k}) ^n}{\omega_{k, \epsilon }^n}\Big)^{\frac{\gamma(1+\tilde{\gamma})}{\tilde{\gamma}-\gamma}}\frac{(\pi_{k}^{\ast}\eta_{k-1} +\epsilon_{k}\eta_{k}) ^n}{n!}\Big)^{\frac{(\tilde{\gamma }-\gamma )}{1+\tilde{\gamma }}}.\\
\end{split}
\end{equation}
Taking limit $\epsilon_{k}\to 0$ in (\ref{3.9}) and using (\ref{Int1}) yield
\begin{equation}\label{3.10}
\begin{split}
&\int_{\tilde{M}}\Big(\frac{\eta_{k}^n}{\omega_{k-1, \epsilon }^n}\Big)^{\gamma }\frac{\eta_{k}^n}{n!}\\
= &\lim_{\epsilon_{k} \to 0}\int_{\tilde{M}}\Big(\frac{\eta_{k}^n}{\omega_{k, \epsilon }^n}\Big)^{\gamma }\frac{\eta_{k}^n}{n!}\\
\leq &\lim_{\epsilon_{k} \to 0}\Big( C^{\ast}\Big) ^{\frac{1+\gamma }{1+\tilde{\gamma}}}\Big(\int_{\tilde{M}}\Big(\frac{(\pi_{k}^{\ast}\eta_{k-1} +\epsilon_{k}\eta_{k}) ^n}{\omega_{k, \epsilon }^n}\Big)^{\frac{\gamma(1+\tilde{\gamma})}{\tilde{\gamma}-\gamma}}    \frac{(\pi_{k}^{\ast}\eta_{k-1} +\epsilon_{k}\eta_{k}) ^n}{n!}\Big)^{\frac{(\tilde{\gamma }-\gamma )}{1+\tilde{\gamma }}}\\
= & \Big( C^{\ast}\Big) ^{\frac{1+\gamma }{1+\tilde{\gamma}}}\Big(\int_{M_{k-1}}\Big(\frac{(\eta_{k-1} ) ^n}{\omega_{k-1, \epsilon }^n}\Big)^{\frac{\gamma(1+\tilde{\gamma})}{\tilde{\gamma}-\gamma}}    \frac{(\eta_{k-1} ) ^n}{n!}\Big)^{\frac{(\tilde{\gamma }-\gamma )}{1+\tilde{\gamma }}}.\\
\end{split}
\end{equation}
Repeating the argument in (\ref{3.9}) and taking limit $\epsilon_{i}\to 0$ successively, we know that there exists a positive constant $a^{\ast}\ll \frac{1}{n-1} $  such that
\begin{equation}\label{3.101}
\lim_{\epsilon_{1}\to 0 } \cdots \lim_{\epsilon_{k} \to 0}\int_{\tilde{M}}\Big(\frac{\eta_{k}^n}{\omega_{k, \epsilon }^n}\Big)^{\gamma }\frac{\eta_{k}^n}{n!}\leq \tilde{C}^{\ast }
\end{equation}
for all $0<\gamma \leq  a^{\ast }$, where $\tilde{C}^{\ast }$ is a  uniform constant.

\medskip

\medskip

 \begin{proposition}\label{prop 4.1} Let $\mathcal{E}$ be a
reflexive sheaf on a smooth compact K\"ahler manifold $(M , \omega )$, and
 $\mathcal{F}^{HNS}=\{\mathcal{E}_{i}\}_{i=1}^{l}$ be the HNS-filtration of $\mathcal{E}$ by saturated subsheaves:
\begin{equation}\label{HNS021}
0= \mathcal{E}_{0}\subset \mathcal{E}_{1}\subset \cdots \subset\mathcal{E}_{l-1}\subset \mathcal{E}_{l}=\mathcal{E},
\end{equation}
where every quotient $\mathcal{Q}_{i}=\mathcal{E}_{i}/\mathcal{E}_{i-1}$ is torsion-free and $\omega $-stable. Let $\pi: \tilde{M}\rightarrow M$ be the composition of  a finite sequence of blowing ups along compact sub-manifolds such that  $E=\pi^{\ast}\mathcal{E}/tor{(\pi^{\ast
}\mathcal{E})}$ is locally free. Then there exists a positive constant $\tilde{a}^{\ast}$ such that,  for any $\delta >0$ and any $1\leq p < 1+\tilde{a}^{\ast}$, there is a smooth metric $H$ on $E$ such that
  \begin{equation}\label{a3}
\|\frac{\sqrt{-1}}{2\pi }\Lambda_{\omega
}(F_{(\mathcal{E},
H)})-\Psi(\mathcal{F}, (\mu_{ 1, \omega} , \cdots , \mu_{l, \omega }) ,
H)\|_{L^{p}(M, \omega)}\leq \delta ,
\end{equation}
where $ \mu_{i, \omega }$ is the $\omega$-slope of $\mathcal{E}_{i}$ and $(\mathcal{E},
H)$ is the Chern connection on $\mathcal{E}$ with respect to the metric $H$.
\end{proposition}

\medskip

{\bf Proof.}  Consider the filtration  $\widetilde{\mathcal{F}}=\{\widetilde{\mathcal{E}_{i}}\}_{i=1}^{l}$  of $E$ over $\tilde{M}$ which is constructed in Proposition \ref{claim}. Every quotient sheaf $\widetilde{\mathcal{Q}}_{i}$ in the filtration (\ref{HNS02}) is {\bf $\omega_{\epsilon}$-stable} for all $0<|\epsilon | \leq \epsilon^{\ast}$. Following Sibley's argument (\cite{Sib}), we can construct a resolution of the filtration (\ref{HNS02}) such that the pullback bundle
 has a filtration by subbundles, which away from the exceptional divisor  is precisely the filtration (\ref{HNS02}).  Using Daskalopoulos and Wentworth's  cut-off argument (\cite{DW1}) and Sibley's trick (\cite{Sib}), we can obtain an $L^{p}$-approximate critical Hermitian metric on the holomorphic vector bundle $E$ over the K\"ahler manifold $(\tilde{M}, \omega_{\epsilon})$ (Theorem 5.12 in \cite{Sib}).
Given $\delta '>0$ and any $1\leq p' < +\infty$,  for every small $\epsilon '$, there exists a smooth Hermitian metric $H_{\epsilon '}$ on $E$ such that
 \begin{equation}\label{a31}
\|\frac{\sqrt{-1}}{2\pi }\Lambda_{\omega_{\epsilon '}
}(F_{(E, \overline{\partial }_{E} ,
H_{\epsilon '})})-\Psi(\widetilde{\mathcal{F}}, (\mu_{ 1, \omega_{\epsilon '}} , \cdots , \mu_{l, \omega_{\epsilon '} }) ,
H_{\epsilon '})\|_{L^{p'}(\tilde{M}\setminus \pi^{-1}(\Sigma_{HNS}), \omega_{\epsilon '})}\leq \delta '.
\end{equation}

We choose a smooth metric $H_{\epsilon '}$ satisfying (\ref{a31}) for some $\epsilon '$ and $p'$ which will be chosen later. For simplicity, denote $\Theta_{1} =\frac{\sqrt{-1}}{2\pi}(F_{(E, \overline{\partial }_{E} ,
H_{\epsilon '})})$.  A straightforward computation shows that
 \begin{equation}\label{ZZ1}
 \begin{split}
&\|\frac{\sqrt{-1}}{2\pi}\Lambda_{\omega_{\epsilon}
}(F_{(E, \overline{\partial }_{E} ,
H_{\epsilon '})})-\Psi(\widetilde{\mathcal{F}}, (\mu_{ 1, \omega_{\epsilon }} , \cdots , \mu_{l, \omega_{\epsilon } }) ,
H_{\epsilon '})\|_{L^{p}(\tilde{M}, \omega_{\epsilon})}\\
&\leq \|\Lambda_{\omega_{\epsilon}
}(\Theta_{1}-\frac{\omega_{\epsilon '}}{n}\Psi(\widetilde{\mathcal{F}}, (\mu_{ 1, \omega_{\epsilon '}} , \cdots , \mu_{l, \omega_{\epsilon '} }) ,
H_{\epsilon '}))\|_{L^{p}(\tilde{M}, \omega_{\epsilon})}\\
&\quad +\| \frac{1}{n}\Lambda_{\omega_{\epsilon}
}(\omega_{\epsilon '}-\omega_{\epsilon })\Psi(\widetilde{\mathcal{F}}, (\mu_{ 1, \omega_{\epsilon '}} , \cdots , \mu_{l, \omega_{\epsilon '} }) ,
H_{\epsilon '}) \|_{L^{p}(\tilde{M}, \omega_{\epsilon})}\\
&\quad +\|\Psi(\widetilde{\mathcal{F}}, (\mu_{ 1, \omega_{\epsilon }} , \cdots , \mu_{l, \omega_{\epsilon } }) ,
H_{\epsilon '})- \Psi(\widetilde{\mathcal{F}}, (\mu_{ 1, \omega_{\epsilon '}} , \cdots , \mu_{l, \omega_{\epsilon '} }) ,
H_{\epsilon '}) \|_{L^{p}(\tilde{M}, \omega_{\epsilon})}.
\end{split}
\end{equation}

Clearly the Chern-Weil theory implies
\begin{equation}
\begin{split}
&\int_{\tilde{M}} (2c_{2}(E)-c_{1}(E)\wedge c_{1}(E))\wedge\frac{\omega_{\epsilon '}^{n-2}}{(n-2)!}\\
=&\int_{\tilde{M}}|\Theta_{1}|_{H_{\epsilon '}}^{2}-|\Lambda_{\omega_{\epsilon '}} \Theta_{1}|_{H_{\epsilon '}}^{2} \frac{\omega_{\epsilon '}^{n}}{n!},\\
\end{split}
\end{equation}
Setting $\Theta_{2}=\Theta_{1}-\frac{\omega_{\epsilon '}}{n}\Psi(\widetilde{\mathcal{F}}, (\mu_{ 1, \omega_{\epsilon '}} , \cdots , \mu_{l, \omega_{\epsilon '} }) ,
H_{\epsilon '})$, we know that $\|\Theta_{2}\|_{L^{2 }(\tilde{M}, \omega_{\epsilon '})}$ is bounded uniformly.
In the sequel, we always assume that $1\leq p < 1+\frac{a^{\ast}}{2+a^{\ast}}$.
From the definition of $\Lambda_{\omega}$, it follows that
\begin{equation}\label{LLL1}
\begin{split}
&\|\Lambda_{\omega_{\epsilon }}\Theta_{2}\|_{L^{p }(\tilde{M}, \omega_{\epsilon})}= \Big\|\frac{\omega_{\epsilon '}^{n}}{\omega_{\epsilon}^{n}}\Big(\Lambda_{\omega_{\epsilon ' }}\Theta_{2} +n \frac{\Theta_{2} \wedge (\omega_{\epsilon}^{n-1}- \omega_{\epsilon ' }^{n-1})}{\omega_{\epsilon '}^{n}}\Big)\Big\|_{L^{p }(\tilde{M}, \omega_{\epsilon})}\\
\leq&  \Big\|\frac{\omega_{\epsilon '}^{n}}{\omega_{\epsilon}^{n}}(\Lambda_{\omega_{\epsilon ' }}\Theta_{2})\Big\|_{L^{p }(\tilde{M}, \omega_{\epsilon})}\\
&+ n\sum_{\beta =1}^{k}|\epsilon_{\beta}-\epsilon_{\beta}'|\Big\|\frac{\Theta_{2} \wedge \eta_{\beta}  \wedge (\sum_{i=0}^{n-2}\omega_{\epsilon}^{n-i-2}\wedge \omega_{\epsilon ' }^{i})}{\omega_{\epsilon }^{n}}\Big\|_{L^{p }(\tilde{M}, \omega_{\epsilon})}.\\
\end{split}
\end{equation}
Direct calculations show that
\begin{equation}\label{LLL2}
\begin{split}
&\Big\|\frac{\omega_{\epsilon '}^{n}}{\omega_{\epsilon}^{n}}(\Lambda_{\omega_{\epsilon ' }}\Theta_{2} )\Big\|_{L^{p }(\tilde{M}, \omega_{\epsilon})}^{p }
= \int_{\tilde{M}} |\Lambda_{\omega_{\epsilon ' }}\Theta_{2}|^{p } \Big(\frac{\omega_{\epsilon '}^{n}}{\omega_{\epsilon}^{n}}\Big)^{p -1} \frac{\omega_{\epsilon '}^{n}}{n!}\\
&\leq  \Big(\int_{\tilde{M}} |\Lambda_{\omega_{\epsilon ' }}\Theta_{2}|^{p \cdot a } \frac{\omega_{\epsilon '}^{n}}{n!}\Big)^{\frac{1}{a}}
\Big(\int_{\tilde{M}}  \Big(\frac{\omega_{\epsilon '}^{n}}{\omega_{\epsilon}^{n}}\Big)^{(p -1)\cdot b} \frac{\omega_{\epsilon '}^{n}}{n!}\Big)^{\frac{1}{b}},
\end{split}
\end{equation}
where $p \cdot a=p '$, since $\frac{p }{1-\frac{1}{a^{\ast}}(p -1)}<\tilde{p }$ and $\frac{1}{a}+\frac{1}{b}=1$, we have $(p -1)b <a^{\ast}$.
Let $U$ be a neighborhood of the singularity set $\Sigma_{\mathcal{E}}$. Since $\pi^{
\ast}\omega $ is degenerate only along $\pi^{-1}(\Sigma_{\mathcal{E}})$, there must exist a positive constant $C_{U}$ depending only on $U$ such that  $\pi^{\ast}\omega \geq C_{U}\eta_{\beta } $ on $\tilde{M}\setminus U$ for all $1\leq \beta \leq k$. On the other hand, we can suppose that $\pi^{\ast}\omega \leq C_{M}\eta_{k}$ on $\tilde{M}$ for some positive constant $C_{M}$. Now suppose that $\epsilon <\epsilon '$, then we get
\begin{equation}\label{LLL3}
\begin{split}
&\int_{\tilde{M} \setminus U}\Big|\frac{\Theta_{2} \wedge \eta_{\beta}  \wedge (\sum_{i=0}^{n-2}\omega_{\epsilon}^{n-i-2}\wedge \omega_{\epsilon ' }^{i})}{\omega_{\epsilon }^{n}}\Big|^{
p } \frac{\omega_{\epsilon}^{n}}{n!}\\&= \int_{\tilde{M} \setminus U}\Big|\frac{\Theta_{2} \wedge \eta_{\beta}  \wedge (\sum_{i=0}^{n-2}\omega_{\epsilon}^{n-i-2}\wedge \omega_{\epsilon ' }^{i})}{\omega_{\epsilon ' }^{n}}\Big|^{
p } \Big(\frac{\omega_{\epsilon '}^{n}}{\omega_{\epsilon}^{n}}\Big)^{(p -1)} \frac{\omega_{\epsilon '}^{n}}{n!}\\
&\leq C(n)C_{U}^{-p } \int_{\tilde{M} \setminus U} |\Theta_{2}|_{\omega_{\epsilon_1}}^{p }\Big(\frac{\omega_{\epsilon '}^{n}}{\omega_{\epsilon}^{n}}\Big)^{(p -1)} \frac{\omega_{\epsilon '}^{n}}{n!}\\
&\leq C(n)C_{U}^{-(n+1)p +n }(C_{M}+ |\epsilon ' |)^{n(p -1)} \int_{\tilde{M} \setminus U} |\Theta_{2}|_{\omega_{\epsilon '}}^{p } \frac{\omega_{\epsilon '}^{n}}{n!},\\
\end{split}
\end{equation}
where $C(n)$ is a uniform constant.
On the other hand, we know
\begin{equation}\label{LLL4}
\begin{split}
&\int_{ U}|\epsilon_{\beta}-\epsilon_{\beta}'|^{p }\Big|\frac{\Theta_{2} \wedge \eta_{\beta}  \wedge (\sum_{i=0}^{n-2}\omega_{\epsilon}^{n-i-2}\wedge \omega_{\epsilon ' }^{i})}{\omega_{\epsilon }^{n}}\Big|^{
p } \frac{\omega_{\epsilon}^{n}}{n!}\\&= \int_{ U}|\epsilon_{\beta}-\epsilon_{\beta}'|^{p }\Big|\frac{\Theta_{2} \wedge \eta_{\beta}  \wedge (\sum_{i=0}^{n-2}\omega_{\epsilon}^{n-i-2}\wedge \omega_{\epsilon ' }^{i})}{\omega_{\epsilon ' }^{n}}\Big|^{
p } \Big(\frac{\omega_{\epsilon '}^{n}}{\omega_{\epsilon}^{n}}\Big)^{(p -1)} \frac{\omega_{\epsilon '}^{n}}{n!}\\
&\leq C(n) \Big(\int_{U} |\Theta_{2}|_{\omega_{\epsilon ' }}^{2 } \frac{\omega_{\epsilon '}^{n}}{n!}\Big)^{\frac{p }{2}}
\Big(\int_{U}  \Big(\frac{\omega_{\epsilon '}^{n}}{\omega_{\epsilon}^{n}}\Big)^{\frac{2p -2}{2-p }} \frac{\omega_{\epsilon '}^{n}}{n!}\Big)^{\frac{2-p }{2}}\\
&\leq  C(n) (\Vol(U, \omega_{\epsilon '} ))^{\frac{2-p }{2}(1-\frac{1}{\tilde{b}})}\Big(\int_{U} |\Theta_{2}|_{\omega_{\epsilon ' }}^{2 } \frac{\omega_{\epsilon '}^{n}}{n!}\Big)^{\frac{p}{2}}
\Big(\int_{U}  \Big(\frac{\omega_{\epsilon '}^{n}}{\omega_{\epsilon}^{n}}\Big)^{\frac{2p -2}{2-p }\cdot \tilde{b}} \frac{\omega_{\epsilon '}^{n}}{n!}\Big)^{\frac{2-p }{2\tilde{b}}},\\
\end{split}
\end{equation}
where $\tilde{b}= \frac{1}{2}(\frac{2p -2}{2-p } +a^{\ast})\cdot (\frac{2p -2}{2-p })^{-1}$, and note that the condition on $p$ gives us $\frac{2p -2}{2-p }\cdot \tilde{b} <a^{\ast}$.

Combining  (\ref{LLL1}), (\ref{LLL2}), (\ref{LLL3}) and (\ref{LLL4}), we derive
\begin{equation}\label{Int2}
\begin{split}
&\|\Lambda_{\omega_{\epsilon }}\Theta_{2}\|_{L^{p }(\tilde{M}, \omega_{\epsilon})} \leq  \sum_{\beta =1}^{k}\hat{C}_{1}|\epsilon_{\beta}-\epsilon_{\beta}'|C_{U}^{-(n+1)+\frac{n}{p}}\|\Theta_{2}\|_{L^{p }(\tilde{M}, \omega_{\epsilon '})}\\
&+ \hat{C}_{2}\|\Lambda_{\omega_{\epsilon ' }}\Theta_{2}\|_{L^{p' }(\tilde{M}, \omega_{\epsilon_1})}\Big(\int_{\tilde{M}}  \Big(\frac{\omega_{\epsilon '}^{n}}{\omega_{\epsilon}^{n}}\Big)^{(p -1)\cdot b} \frac{\omega_{\epsilon '}^{n}}{n!}\Big)^{\frac{1}{p b}}\\
& +\hat{C}_{2}(\Vol(U, \omega_{\epsilon_1}))^{C(p, a^{\ast})}\|\Theta_{2}\|_{L^{2 }(\tilde{M}, \omega_{\epsilon '})}\Big(\int_{U}  \Big(\frac{\omega_{\epsilon '}^{n}}{\omega_{\epsilon}^{n}}\Big)^{\frac{2p -2}{2-p }\cdot \tilde{b}} \frac{\omega_{\epsilon '}^{n}}{n!}\Big)^{\frac{2-p }{2p \tilde{b}}},\\
\end{split}
\end{equation}
where $C(p, a^{\ast})=\frac{1-p}{2}-\frac{2(p-1)(2-p)}{2p-2+a^{\ast }(2-p)}$. This together with (\ref{3.101}) implies
\begin{equation}\label{Int21}
\begin{split}
&\lim_{\epsilon_{1}\to 0 } \cdots \lim_{\epsilon_{k} \to 0}\|\Lambda_{\omega_{\epsilon }}\Theta_{2}\|_{L^{p }(\tilde{M}, \omega_{\epsilon})} \leq  \sum_{\beta =1}^{k}\hat{C}_{1}|\epsilon_{\beta}'|C_{U}^{-(n+1)+\frac{n}{p}}\|\Theta_{2}\|_{L^{p }(\tilde{M}, \omega_{\epsilon '})}\\
&+ \hat{C}^{\ast}\hat{C}_{2}\Big(\|\Lambda_{\omega_{\epsilon ' }}\Theta_{2}\|_{L^{p' }(\tilde{M}, \omega_{\epsilon_1})}+(\Vol(U, \omega_{\epsilon_1}))^{C(p, a^{\ast})}\|\Theta_{2}\|_{L^{2 }(\tilde{M}, \omega_{\epsilon '})}\Big),\\
\end{split}
\end{equation}
where $\hat{C}^{\ast}$ is a uniform constant.
We may choose $U$ such that $\Vol(U, \omega_{1})$ small enough first, and then $\tilde{\delta }$ and $\epsilon '$ both sufficiently small so that
\begin{equation}
\lim_{\epsilon_{1}\to 0 } \cdots \lim_{\epsilon_{k} \to 0}\|\Lambda_{\omega_{\epsilon }}\Theta_{2}\|_{L^{p }(\tilde{M}, \omega_{\epsilon})} \leq \frac{\delta }{3}.
\end{equation}

 By (\ref{3.101}) and the fact that $\mu_{i, \omega_{\epsilon }} \rightarrow \mu_{i, \omega }$ as $\epsilon \rightarrow 0$,  we may choose $\epsilon '$ small enough so that the second and third terms in (\ref{ZZ1}) are both smaller than $\frac{\delta }{3}$, hence it follows that
 \begin{equation}\label{a310}
\|\frac{\sqrt{-1}}{2\pi }\Lambda_{\omega
}(F_{(\mathcal{E},
H_{\epsilon '})})-\Psi(\mathcal{F}, (\mu_{ 1, \omega} , \cdots , \mu_{l, \omega }) ,
H_{\epsilon '})\|_{L^{p}(M, \omega)}\leq \delta .
\end{equation}

\hfill $\Box$ \\

\medskip

\section{ The HN type of the Uhlenbeck limit}
\setcounter{equation}{0}

\medskip

Let $\mathcal{E}$ be a
reflexive sheaf on a smooth K\"ahler manifold $(M , \omega )$,  $H(t)$ be a solution of the Hermitian-Yang-Mills flow (\ref{SSS1}) on $\mathcal{E}|_{M\setminus \Sigma_{\mathcal{E}}}$ with the initial metric $\hat{H}$, and $A(t)$ be the related Yang-Mills flow (\ref{YM2}) on the Hermitian vector bundle $(\mathcal{E}|_{M\setminus \Sigma_{\mathcal{E}}}, \hat{H})$.
Let $A_{\infty} $ be an Uhlenbeck
limit. From Theorem \ref{thmlimit}, we know that $A_{\infty} $ is a smooth Yang-Mills connection on the Hermitian
bundle $(E_{\infty} , H_{\infty})$ over $M\setminus (\Sigma_{\mathcal{E}}\cup \Sigma_{an})$, and $\theta (A_{\infty}, \omega )$ is parallel, then the constant eigenvalues vector
$\vec{\lambda}_{\infty}=(\lambda_{1} , \cdots , \lambda_{R})$ of
$\frac{\sqrt{-1}}{2\pi }\Lambda_{\omega}F_{A_{\infty}}$ is just the HN type
of the extended Uhlenbeck limit  sheaf $\mathcal{E}_{\infty}= (E_{\infty} , \overline{\partial}_{A_{\infty}} )$. Denote by $\vec{\mu}_{0}=(\mu_{1} ,
\cdots , \mu_{R})$ the HN type  of the reflexive
sheaf $\mathcal{E}$. In this section, we will show  that the HN type of the limiting  sheaf for the Hermitian-Yang-Mills flow (\ref{SSS1}) is in fact equal to the HN type of the reflexive sheaf $\mathcal{E}$, i.e. $\vec{\lambda }_{\infty}=\vec{\mu}_{0}$.

\begin{lemma}\label{lemma3.5}
Let $A(t)$ be the long time solution of the Yang-Mills flow (\ref{YM2})  on
a complex vector bundle $\mathcal{E}|_{M\setminus \Sigma_{\mathcal{E}}}$ of rank $R$ with a Hermitian metric
$\hat{H}$. Let $S$ be a coherent  subsheaf of $\mathcal{E}$. Suppose there is  a sequence $\{A_j\}$, modulo gauge transformations, such that $\sqrt{-1}\Lambda_{\omega
}(F_{A_{j}})\rightarrow \mathcal{B}$ in
$L^{1}$ as $j\rightarrow +\infty $, where $\mathcal{B}\in
L^{1}(\sqrt{-1}\verb"u"(E))$, and the eigenvalues $\lambda_{1 }\geq
\cdots \geq \lambda_{R}$ of $\frac{1}{2\pi}\mathcal{B}$ are constant almost
everywhere. Then: $\deg_{\omega} (S)\leq \sum_{i\leq \rank(S)}\lambda_{i}$.
\end{lemma}

{\bf Proof.} Because $\deg_{\omega} (S)\leq \deg_{\omega}(Sat_{E}(S))$, we may assume that
$S$ is saturated.
 As before, let $\pi : \tilde{M}\rightarrow M $ be the composition of  a finite sequence of blowups resolving the sheaf $\mathcal{E}$, i.e. such that $E=\pi^\ast\mathcal{E}/tor(\pi^\ast\mathcal{E})$ is locally free. Considering the exact sequence
\begin{equation}
0\rightarrow S\rightarrow \mathcal{E}\rightarrow Q\rightarrow 0, \quad \rm{over}\ M,
\end{equation}
 we get the following exact sequences
\begin{equation}
\pi^\ast S\xrightarrow{f_1}\pi^\ast\mathcal{E}\rightarrow \pi^\ast Q\rightarrow 0,  \quad \rm{over}\ \tilde{M},
\end{equation}
and
\begin{equation}
0\rightarrow Im(f_1)/tor(Im(f_1))\rightarrow E=\pi^\ast\mathcal{E}/tor(\pi^\ast\mathcal{E})\rightarrow \check{Q}\rightarrow 0,  \quad \rm{over}\ \tilde{M},
\end{equation}
where $\check{Q}=\frac{E}{Im(f_1)/tor(Im(f_1))}$. Setting $\hat{S}=Sat_E(Im(f_1)/tor(Im(f_1)))$, since $\pi $ is biholomorphic outside $\Sigma_{\mathcal{E}}$ and codim$(\Sigma_{\mathcal{E}})\geq 3$, we have $(\pi)_\ast \hat{S}=S$ on $M\setminus\Sigma_{\mathcal{E}}$.

Let  $H_{k, \epsilon } (t)$  be the long time solution of the Hermitian-Yang-Mills flow (\ref{DDD1}) on the holomorphic bundle $E $ over $\tilde{M}$ with the fixed smooth initial metric $\hat{H}$ and with respect to the K\"ahler metric $\omega_{k, \epsilon}$. Clearly Lemma \ref{lemmaBS} and Proposition \ref{propgauge} say that $H_{k, \epsilon}(x, t)$  converges successively to the long time solution $H(x, t)$ of  the Hermitian-Yang-Mills flow (\ref{SSS1}) as $\epsilon\rightarrow 0$, $A(t)=\sigma (t) (D_{\hat{H}})$ and $\sigma (t)^{\ast }\circ \sigma (t)=\hat{H}^{-1}H(t)$.
Let $\pi^{H(t)}$ (resp. $\pi ^{H_{k, \epsilon }(t)}$) denote the orthogonal projection
onto $S$ (resp. $\hat{S}$) with respect to the Hermitian metric $H(t)$ (resp. $H_{k, \epsilon }(t)$).
Using the Gauss-Codazzi equation and Fatou's lemma, we derive
\begin{equation}\label{deg}
\begin{split}
\deg_{\omega}(S)= &\int_M c_1(\det S)\wedge\frac{\omega^{n-1}}{(n-1)!}\\
= &\int_M c_1(\pi_\ast(\det\hat{S}))\wedge\frac{\omega^{n-1}}{(n-1)!}\\
= &\int_{\tilde{M}} c_1(\det\hat{S})\wedge\frac{\pi^\ast\omega^{n-1}}{(n-1)!}\\
= &\lim_{\epsilon\rightarrow 0}\int_{\tilde{M}}c_1(\hat{S})\wedge\frac{\omega_\epsilon^{n-1}}{(n-1)!}\\
= &\lim_{\epsilon\rightarrow 0}\frac{1}{2\pi}\int_{\tilde{M}}(\tr(\sqrt{-1}\Lambda_{\omega_\epsilon}(F_{H_\epsilon(t)})\pi^{H_\epsilon(t)})-|\bar{\partial}_{E}\pi^{H_\epsilon(t)}|^2)\frac{\omega_{\epsilon }^{n}}{n!}\\
\leq & \frac{1}{2\pi}\int_{M}(\tr(\sqrt{-1}\Lambda_{\omega }(F_{H(t)})\pi^{H(t)})-|\bar{\partial}_{\mathcal{E}}\pi^{H(t)}|^2)\frac{\omega^{n}}{n!}\\
\leq & \frac{1}{2\pi}\int_{M}\tr(\sqrt{-1}\Lambda_{\omega}(F_{H(t)})\pi^{H(t)})\frac{\omega^{n}}{n!}\\
 = & \frac{1}{2\pi}\int_{M}\tr(\sqrt{-1}\Lambda_{\omega}(\sigma^{-1}(t)\circ F_{A(t)}\circ \sigma (t))\pi^{H(t)})\frac{\omega^{n}}{n!}\\
 = & \frac{1}{2\pi}\int_{M}\tr((\sqrt{-1}\Lambda_{\omega} F_{A(t)}-\mathcal{B}) (\sigma (t)\circ \pi^{H(t)}\circ \sigma^{-1}(t)))\frac{\omega^{n}}{n!}\\
&   +\frac{1}{2\pi}\int_{M}\tr(\mathcal{B} (\sigma (t)\circ \pi^{H(t)}\circ \sigma^{-1}(t)))\frac{\omega^{n}}{n!}\\
\end{split}
\end{equation}
for $t>0$.
By a result from linear algebra (Lemma 2.20 in \cite{DW1}), we obtain $\frac{1}{2\pi}\tr(\mathcal{B}(\sigma (t)\circ \pi^{H(t)}\circ \sigma^{-1}(t)))\leq \sum_{i\leq \rank(S)}\lambda_{i}$. So it holds that
$\deg_{\omega} (S)\leq \sum_{i\leq \rank(S)}\lambda_{i} +\frac{1}{2\pi}\|\sqrt{-1}\Lambda_{\omega}(F_{A(t_{j})})-\mathcal{B}\|_{L^{1}}$. Letting $j\rightarrow \infty$ concludes the proof of the lemma.

\hfill $\Box$ \\

\medskip

Combining (\ref{F1}), Lemma \ref{lemmaBS} and Corollary \ref{Corollary 2.5.}, we know
 \begin{equation}
 \int_{\tilde{M}} \tr (\theta (H_{k, \epsilon }(t), \omega_{k, \epsilon}))\frac{\omega_{\epsilon}^{n}}{n!}=\int_{\tilde{M}} \tr (\theta (\hat{H},  \omega_{k, \epsilon}))\frac{\omega_{\epsilon}^{n}}{n!},
 \end{equation}
and then
 \begin{equation}
 \deg_{\omega}(\mathcal{E})=\int_{M} \tr (\frac{\sqrt{-1}}{2\pi }\Lambda_{\omega}F_{\hat{H}})\frac{\omega^{n}}{n!}=\int_{M} \tr (\frac{\sqrt{-1}}{2\pi }\Lambda_{\omega}F_{A_{\infty}})\frac{\omega^{n}}{n!}=\deg_{\omega}(E_{\infty}, \overline{\partial}_{A_{\infty}}),
 \end{equation}
i.e.
\begin{equation}
\sum_{\alpha =1}^{R}\mu_{\alpha}=\sum_{\alpha =1}^{R}\lambda_{\alpha}.
\end{equation}
Let $\{\mathcal{E}_{i}\}_{i=1}^{l}$ be the HNS-filtration of the reflexive sheaf $\mathcal{E}$. Applying Lemma \ref{lemma3.5} yields:
\begin{equation}
\sum_{\alpha \leq \rank \mathcal{E}_{i}} \mu_{\alpha }=\deg_{\omega}(\mathcal{E}_{i})\leq \sum_{\alpha \leq \rank \mathcal{E}_{i}} \lambda_{\alpha }
\end{equation}
for all $i$. Of course Lemma 2.3 in \cite{DW1} means \begin{equation}\label{P5.5}\vec{\mu }_{0}\leq
\vec{\lambda}_{\infty}.\end{equation}

\hspace{0.3cm}

For further consideration, we show the continuous dependence of the Hermitian-Yang-Mills flow (\ref{SSS1}) on initial metrics.

\medskip

\begin{lemma}\label{continuous} Let $\hat{H}_{1}$ and $\hat{H}_{2}$ be two smooth metrics on the holomorphic bundle $E$ over $\tilde{M}$, and
 $\tilde{\delta }=\sup_{x\in \tilde{M}}(\tr \hat{H}_{1}^{-1}\hat{H}_{2}+\tr \hat{H}_{2}^{-1}\hat{H}_{1}-2\, \rank(E))$. If
 $H_{i}(t)$ is the long time solution of the Hermitian-Yang-Mills flow (\ref{SSS1}) on $\mathcal{E}|_{M\setminus \Sigma_{\mathcal{E}}}$ with the initial metric $\hat{H}_{i}$ respectively for $i=1, 2$, then for any $t>0$,
\begin{equation}\label{CON1}
\int_{M}|\sqrt{-1}\Lambda_{\omega}F_{H_2(t)}-\sqrt{-1}\Lambda_{\omega}F_{H_1(t)}|_{H_1(t)}^2\frac{\omega ^n}{n!}dt \leq f_{t}(\tilde{\delta }),
\end{equation}
where  $f_{t}:R^{+}\rightarrow R^{+}$ is a continuous function satisfying $f_{t} (x) \rightarrow 0$ as $x\rightarrow 0$.

\end{lemma}

\medskip

{\bf Proof.}
Let $H_{(i, \epsilon)} (t)$ be the long time solution of the Hermitian-Yang-Mills flow (\ref{SSS1}) on the holomorphic bundle $E $ over $\tilde{M}$ with the  smooth initial metric $\hat{H}_{i}$ and with respect to the K\"ahler metric $\omega_{\epsilon}$, where $i=1, 2$. Set
\begin{equation}
\tilde{h}_{\epsilon }(t)=H_{(1, \epsilon)}^{-1}(t)H_{(2, \epsilon)}(t).
\end{equation}
It is easy to check that
\begin{equation}
\tr(\sqrt{-1}\tilde{h}_{\epsilon }(t)(\Lambda_{\omega_\epsilon}F_{H_{(2, \epsilon)}(t)}-\Lambda_{\omega_\epsilon}F_{H_{(1, \epsilon)}(t)}))=-\frac{1}{2}
\Delta_{\omega_\epsilon}\tr\tilde{h}_{\epsilon }(t)+\tr(-\sqrt{-1}\Lambda_{\omega_\epsilon}\bar{\partial}\tilde{h}_{\epsilon }(t)
\tilde{h}_{\epsilon }^{-1}(t)\partial_{H_1(t)}\tilde{h}_{\epsilon }(t)),
\end{equation}
\begin{equation}
\tr(\sqrt{-1}\tilde{h}_{\epsilon }^{-1}(t)(\Lambda_{\omega_\epsilon}F_{H_{(1, \epsilon)}(t)}-\Lambda_{\omega_\epsilon}F_{H_{(2, \epsilon)}(t)}))=-\frac{1}{2}
\Delta_{\omega_\epsilon}\tr\tilde{h}_{\epsilon }^{-1}(t)+\tr(-\sqrt{-1}\Lambda_{\omega_\epsilon}\bar{\partial}\tilde{h}_{\epsilon }^{-1}(t)
\tilde{h}_{\epsilon }(t)\partial_{H_{(2, \epsilon)}(t)}\tilde{h}_{\epsilon }^{-1}(t)),
\end{equation}
\begin{equation}\label{A}
\begin{split}
&\tr(\sqrt{-1}(\tilde{h}_{\epsilon }(t)-\tilde{h}_{\epsilon }^{-1}(t))(\Lambda_{\omega_\epsilon}F_{H_{(2, \epsilon)}(t)}-\Lambda_{\omega_\epsilon}F_{H_{(1, \epsilon)}(t)}))\\
=&-\frac{1}{2}\Delta_{\omega_\epsilon}(\tr\tilde{h}_{\epsilon }(t)+\tr\tilde{h}_{\epsilon }^{-1}(t))+\tr(-\sqrt{-1}\Lambda_{\omega_\epsilon}
\bar{\partial}\tilde{h}_{\epsilon }(t)\tilde{h}_{\epsilon }^{-1}(t)\partial_{H_{(1, \epsilon)}(t)}\tilde{h}_{\epsilon }(t))\\
&+\tr(-\sqrt{-1}\Lambda_{\omega_\epsilon}
\bar{\partial}\tilde{h}_{\epsilon }^{-1}(t)\tilde{h}_{\epsilon }(t)\partial_{H_{(2, \epsilon)}(t)}\tilde{h}_{\epsilon }^{-1}(t)),
\end{split}
\end{equation}
and
\begin{equation}\label{C0mono}
(\Delta_{\omega_\epsilon}-\frac{\partial}{\partial t})(\tr\tilde{h}_{\epsilon }(t)+\tr\tilde{h}_{\epsilon }^{-1}(t)-2\, \rank(E))\geq 0.
\end{equation}
The inequality (\ref{C0mono}) together with the maximum principle gives us
\begin{equation}
\sup_{x\in \tilde{M}}(\tr\tilde{h}_{\epsilon }(t)+\tr\tilde{h}_{\epsilon }^{-1}(t)-2\, \rank(E))\leq \sup_{x\in \tilde{M}}(\tr\tilde{h}_{\epsilon }(0)+\tr\tilde{h}_{\epsilon }^{-1}(0)-2\, \rank(E))=\tilde{\delta}.
\end{equation}
In the following, we  assume that $\tilde{\delta}$ is small enough. Suppose $\lambda_i$ is the eigenvalue of $\tilde{h}_{\epsilon }(t)$ for $1\leq i \leq n$, then
\begin{equation}
1-\sqrt{\tilde{\delta}(\tilde{\delta}+4n)}\leq \lambda_i \leq 1+\sqrt{\tilde{\delta}(\tilde{\delta}+4n)},
\end{equation}
\begin{equation}\label{eigen1}
1-\sqrt{\tilde{\delta}(\tilde{\delta}+4n)}\leq \frac{1}{\lambda_i} \leq 1+\sqrt{\tilde{\delta}(\tilde{\delta}+4n)},
\end{equation}
and
\begin{equation}\label{eigen2}
-2\sqrt{\tilde{\delta}(\tilde{\delta}+4n)}\Id\leq \tilde{h}_{\epsilon }(t)-\tilde{h}_{\epsilon }^{-1}(t)
\leq 2\sqrt{\tilde{\delta}(\tilde{\delta}+4n)}\Id.
\end{equation}
A direct computation shows that
\begin{equation}
\begin{split}
&|\tr(\sqrt{-1}(\tilde{h}_{\epsilon }(t)-\tilde{h}_{\epsilon }^{-1}(t))(\Lambda_{\omega_\epsilon}F_{H_{(2, \epsilon)}(t)}-\Lambda_{\omega_\epsilon}F_{H_{(1, \epsilon)}(t)}))|\\
\leq&2\sqrt{n}\sqrt{\tilde{\delta}(\tilde{\delta}+4n)}|\Lambda_{\omega_\epsilon}F_{H_{(2, \epsilon)}(t)}-\Lambda_{\omega_\epsilon}F_{H_{(1, \epsilon)}(t)}|_{H_{(1, \epsilon)}(t)}\\
\leq&\sqrt{\tilde{\delta}(\tilde{\delta}+4n)}\Big(\sqrt{n}\Big(1+\sqrt{\tilde{\delta}(\tilde{\delta}+4n)}\Big)|\Lambda_{\omega_\epsilon}F_{H_{(2, \epsilon)}(t)}|_{H_{(2, \epsilon)}(t)}
+|\Lambda_{\omega_\epsilon}F_{H_{(1, \epsilon)}(t)}|_{H_{(1, \epsilon)}(t)}\Big).
\end{split}
\end{equation}
Set $\tilde{T}_{\epsilon}(t)=\tilde{h}_{\epsilon }^{-1}(t)\partial_{H_{(1, \epsilon)}(t)}\tilde{h}_{\epsilon }(t)$. Clearly (\ref{eigen1}) implies that
\begin{equation}
|\tilde{T}_{\epsilon}(t)|^2_{H_1(t), \omega_\epsilon}\leq \sqrt{n}\Big(1+\sqrt{\tilde{\delta}(\tilde{\delta}+4n)}\Big)\tr(-\sqrt{-1}\Lambda_{\omega_\epsilon}\bar{\partial}\tilde{h}_{\epsilon }(t)\tilde{h}_{\epsilon }^{-1}(t)\partial_{H_{(1, \epsilon)}(t)}\tilde{h}_{\epsilon }(t)).
\end{equation}
From (\ref{H001}), it follows that there exists a uniform constant $C$ such that
\begin{equation}\label{H00011}
\int_{\tilde{M}}|\Lambda_{\omega_\epsilon}F_{H_{(2, \epsilon)}(t)}|_{H_{(2, \epsilon)}(t)}\frac{\omega_\epsilon^n}{n!}+\int_{\tilde{M}}|\Lambda_{\omega_\epsilon}F_{H_{(1, \epsilon)}(t)}|_{H_{(1, \epsilon)}(t)}\frac{\omega_\epsilon^n}{n!}\leq C,
\end{equation}
for any $t\geq 0$ and $0<\epsilon \leq 1$.
Combining (\ref{A}) and (\ref{eigen2}), we get
\begin{equation}\label{C6}
\begin{split}
\int_{\tilde{M}}|\tilde{T}_{\epsilon}(t)|_{H_1(t), \omega_\epsilon}^2\frac{\omega_\epsilon^n}{n!}\leq & Cn\Big(1+\sqrt{\tilde{\delta}(\tilde{\delta}+4n)}\Big)^2\sqrt{\tilde{\delta}(\tilde{\delta}+4n)}\\
\leq& C\sqrt{\tilde{\delta}}.
\end{split}
\end{equation}
By straightforward calculations, we deduce
\begin{equation}\label{C}
\begin{split}
\frac{\partial}{\partial t}\tilde{T}_{\epsilon}(t)=& -\tilde{h}_{\epsilon }^{-1}(t)\frac{\partial \tilde{h}_{\epsilon }(t)}{\partial t}\tilde{h}_{\epsilon }^{-1}(t)\partial_{H_{(1, \epsilon)}(t)}\tilde{h}_{\epsilon }(t)+\tilde{h}_{\epsilon }^{-1}(t)\partial_{H_{(1, \epsilon)}(t)}\Big(\tilde{h}_{\epsilon }(t)\tilde{h}_{\epsilon }^{-1}(t)\frac{\partial \tilde{h}_{\epsilon }(t)}{\partial t}\Big)\\
&+\tilde{h}_{\epsilon }^{-1}(t)\frac{\partial}{\partial t}(\partial_{H_{(1, \epsilon)}(t)})\tilde{h}_{\epsilon }(t)\\
=&\partial_{H_{(2, \epsilon)}(t)}\Big(\tilde{h}_{\epsilon }^{-1}(t)\frac{\partial \tilde{h}_{\epsilon }(t)}{\partial t}\Big)+\tilde{h}_{\epsilon }^{-1}(t)\frac{\partial}{\partial t}\Big(H_{(1, \epsilon)}^{-1}(t)\partial H_{(1, \epsilon)}(t)\Big)\tilde{h}_{\epsilon }(t)\\
&-\tilde{h}_{\epsilon }^{-1}(t)\tilde{h}_{\epsilon }(t)\frac{\partial}{\partial t}\Big(H_{(1, \epsilon)}^{-1}(t)\partial H_{(1, \epsilon)}(t)\Big)\\
=&\partial_{H_{(2, \epsilon)}(t)}\Big(\tilde{h}_{\epsilon }^{-1}(t)\frac{\partial \tilde{h}_{\epsilon }(t)}{\partial t}\Big)+\partial_{H_{(2, \epsilon)}(t)}\Big(\tilde{h}_{\epsilon }^{-1}(t)\Big(H_{(1, \epsilon)}^{-1}(t)\frac{\partial H_{(1, \epsilon)}(t)}{\partial t}\Big)\tilde{h}_{\epsilon }(t)\Big)\\
&-\partial_{H_{(1, \epsilon)}(t)}\Big(H_{(1, \epsilon)}^{-1}(t)\frac{\partial H_{(1, \epsilon)}(t)}{\partial t}\Big)\\
=&2\partial_{H_{(1, \epsilon)}(t)}(\sqrt{-1}\Lambda_{\omega_\epsilon}F_{H_{(1, \epsilon)}(t)}-\sqrt{-1}\Lambda_{\omega_\epsilon}F_{H_{(2, \epsilon)}(t)})\\
&-2\tilde{T}_{\epsilon}\sqrt{-1}\Lambda_{\omega_\epsilon}F_{H_{(2, \epsilon)}(t)}+2\sqrt{-1}\Lambda_{\omega_\epsilon}F_{H_{(2, \epsilon)}(t)}\tilde{T}_{\epsilon}(t)\\
=&-2\partial_{H_{(1, \epsilon)}(t)}(\sqrt{-1}\Lambda_{\omega_\epsilon}\bar{\partial}_{E}\tilde{T}_{\epsilon}(t))-2\tilde{T}_{\epsilon}(t)\sqrt{-1}\Lambda_{\omega_\epsilon}F_{H_{(2, \epsilon)}(t)}
+2\sqrt{-1}\Lambda_{\omega_\epsilon}F_{H_{(2, \epsilon)}(t)}\tilde{T}_{\epsilon}(t),
\end{split}
\end{equation}
and
\begin{equation}
\begin{split}
\frac{\partial}{\partial t}|\tilde{T}_{\epsilon}(t)|_{H_{(1, \epsilon)}(t)}^2=&\frac{\partial}{\partial t}\sqrt{-1}\Lambda_{\omega_\epsilon}\tr(\tilde{T}_{\epsilon}(t)H_{(1, \epsilon)}^{-1}(t)\wedge \overline{\tilde{T}_{\epsilon}}^{\mathrm{t}}(t)H_{(1, \epsilon)}^{-1}(t))\\
=&2Re \sqrt{-1}\Lambda_{\omega_\epsilon}\tr\Big(\frac{\partial}{\partial t}\tilde{T}_{\epsilon}(t)H_{(1, \epsilon)}^{-1}(t)\wedge \overline{\tilde{T}_{\epsilon}}^{\mathrm{t}}(t)H_{(1, \epsilon)}(t)\Big)\\
&+\sqrt{-1}\Lambda_{\omega_{\epsilon}}\tr(\tilde{T}_{\epsilon}(t)(-H_{(1, \epsilon)}^{-1}(t)\frac{\partial H_{(1, \epsilon)}(t)}{\partial t}H_{(1, \epsilon)}^{-1}(t))\wedge\overline{\tilde{T}_{\epsilon}}^{\mathrm{t}}(t)H_{(1, \epsilon)}(t))\\
&+\sqrt{-1}\Lambda_{\omega_{\epsilon}}\tr\Big(\tilde{T}_{\epsilon}(t)H_{(1, \epsilon)}^{-1}(t)\wedge\overline{\tilde{T}_{(1, \epsilon)}}^{\mathrm{t}}(t)\frac{\partial H_{(1, \epsilon)}(t)}{\partial t}\Big)\\
\leq&-4Re\sqrt{-1}\Lambda_{\omega_\epsilon}\partial\tr(\sqrt{-1}\Lambda_{\omega_\epsilon}\bar{\partial}_E\tilde{T}_{\epsilon}(t)H_{(1, \epsilon)}^{-1}(t)\wedge\overline{\tilde{T}}^\mathrm{t}(t)H_{(1, \epsilon)}(t))\\
&+4Re\sqrt{-1}\Lambda_{\omega_\epsilon}\tr(\sqrt{-1}\Lambda_{\omega_\epsilon}\bar{\partial}_E\tilde{T}_{\epsilon}(t)H_{(1, \epsilon)}^{-1}(t)\wedge\overline{\bar{\partial}_E\tilde{T}_{\epsilon}}^\mathrm{t}(t)H_{(1, \epsilon)}(t))\\
&+8|\sqrt{-1}\Lambda_{\omega_\epsilon}F_{H_{(2, \epsilon)}(t)}|_{H_{(1, \epsilon)}(t)}|\tilde{T}_{\epsilon}(t)|_{H_{(1, \epsilon)}(t)}^2+4|\sqrt{-1}\Lambda_{\omega_\epsilon}F_{H_{(1, \epsilon)}(t)}|_{H_{(1, \epsilon)}(t)}|\tilde{T}_{\epsilon}(t)|_{H_{(1, \epsilon)}(t)}^2.
\end{split}
\end{equation}
Recall that
\begin{equation}
F_{H_{(2, \epsilon)}(t)}-F_{H_{(1, \epsilon)}(t)}=\bar{\partial}_{E}(\tilde{h}_{\epsilon}^{-1}(t)\partial_{H_{(1, \epsilon)}(t)}\tilde{h}_{\epsilon}(t)),
\end{equation}
we have
\begin{equation}
\begin{split}
&\frac{\partial}{\partial t}\int_{\tilde{M}}|\tilde{T}_{\epsilon}(t)|_{H_{(1, \epsilon)}(t), \omega_\epsilon}^2\frac{\omega_\epsilon^n}{n!}\\
\leq& -4\int_{\tilde{M}}|\sqrt{-1}\Lambda_{\omega_\epsilon}F_{H_{(2, \epsilon)}(t)}-\sqrt{-1}\Lambda_{\omega_\epsilon}F_{H_{(1, \epsilon)}(t)}|_{H_{(1, \epsilon)}(t)}^2\frac{\omega_\epsilon^n}{n!}\\
&+\int_{\tilde{M}}(8|\sqrt{-1}\Lambda_{\omega_\epsilon}F_{H_{(2, \epsilon)}(t)}|_{H_{(1, \epsilon)}(t)}+4|\sqrt{-1}\Lambda_{\omega_\epsilon}F_{H_{(1, \epsilon)}(t)}|_{H_{(1, \epsilon)}(t)})|\tilde{T}_{\epsilon}(t)|_{H_{(1, \epsilon)}(t), \omega_\epsilon}^2\frac{\omega_\epsilon^n}{n!}.
\end{split}
\end{equation}
Then it holds that
\begin{equation}\label{B}
\begin{split}
&\frac{\partial}{\partial t}\Big(t^n\int_{\tilde{M}}|\tilde{T}_{\epsilon}(t)|_{H_{(1, \epsilon)}(t), \omega_\epsilon}^2\frac{\omega_\epsilon^n}{n!}\Big)\\
=&nt^{n-1}\int_{\tilde{M}}|\tilde{T}_{\epsilon}(t)|_{H_{(1, \epsilon)}(t), \omega_\epsilon}^2\frac{\omega_\epsilon^n}{n!}+t^n\frac{\partial}{\partial t}\int_{\tilde{M}}|\tilde{T}_{\epsilon}(t)|_{H_{(1, \epsilon)}(t), \omega_\epsilon}^2\frac{\omega_\epsilon^n}{n!}\\
\leq&nt^{n-1}\int_{\tilde{M}}|\tilde{T}_{\epsilon}(t)|_{H_{(1, \epsilon)}(t), \omega_\epsilon}^2\frac{\omega_\epsilon^n}{n!}-4t^n\int_{\tilde{M}}|\sqrt{-1}\Lambda_{\omega_\epsilon}F_{H_{(2, \epsilon)}(t)}-\sqrt{-1}\Lambda_{\omega_\epsilon}F_{H_{(1, \epsilon)}(t)}|_{H_{(1, \epsilon)}(t)}^2\frac{\omega_\epsilon^n}{n!}\\
&+t^n\int_{\tilde{M}}(8|\sqrt{-1}\Lambda_{\omega_\epsilon}F_{H_{(2, \epsilon)}(t)}|_{H_{(1, \epsilon)}(t)}+4|\sqrt{-1}\Lambda_{\omega_\epsilon}F_{H_{(1, \epsilon)}(t)}|_{H_{(1, \epsilon)}(t)})|\tilde{T}_{\epsilon}(t)|_{H_{(1, \epsilon)}(t), \omega_\epsilon}^2\frac{\omega_\epsilon^n}{n!}.
\end{split}
\end{equation}
The inequalities (\ref{H000}), (\ref{H00011}) and (\ref{kernel01}) tell us that there must exist uniform constants $\hat{C}_1$ and $\hat{C}_1$ such that
\begin{equation}\label{C2}
\sup_{\tilde{M}}t^n|\sqrt{-1}\Lambda_{\omega_\epsilon}F_{H_{(1, \epsilon)}(t)}|_{H_{(1, \epsilon)}(t)}\leq \hat{C}_1,
\end{equation}
\begin{equation}\label{C3}
\sup_{\tilde{M}}t^n|\sqrt{-1}\Lambda_{\omega_\epsilon}F_{H_{(2, \epsilon)}(t)}|_{H_1(t)}\leq \hat{C}_2,
\end{equation}
for $t\in [0, 1]$,
and
\begin{equation}\label{C4}
\sup_{\tilde{M}}|\sqrt{-1}\Lambda_{\omega_\epsilon}F_{H_{(1, \epsilon)}(t)}|_{H_{(1, \epsilon)}(t)}\leq \hat{C}_1,
\end{equation}
\begin{equation}\label{C5}
\sup_{\tilde{M}}|\sqrt{-1}\Lambda_{\omega_\epsilon}F_{H_{(2, \epsilon)}(t)}|_{H_{(1, \epsilon)}(t)}\leq \hat{C}_2,
\end{equation}
for $t\geq 1$.
The above inequalities together with (\ref{B}) yield
\begin{equation}\label{C1}
\begin{split}
&\int_0^1t^n\int_{\tilde{M}}|\sqrt{-1}\Lambda_{\omega_\epsilon}F_{H_{(2, \epsilon)}(t)}-\sqrt{-1}\Lambda_{\omega_\epsilon}F_{H_{(1, \epsilon)}(t)}|_{H_{(1, \epsilon)}(t)}^2\frac{\omega_\epsilon^n}{n!}dt\\
\leq&-\frac{1}{4}\int_{\tilde{M}}|\tilde{T}_{\epsilon}(1)|_{H_{(1, \epsilon)}(1), \omega_\epsilon}^2\frac{\omega_\epsilon^n}{n!}+\frac{n}{4}\int_0^1t^{n-1}\int_{\tilde{M}}|\tilde{T}_{\epsilon}(t)|_{H_{(1, \epsilon)}(t), \omega_\epsilon}^2\frac{\omega_\epsilon^n}{n!}dt\\
&+\frac{1}{4}\int_0^1t^n\int_{\tilde{M}}(8|\sqrt{-1}\Lambda_{\omega_\epsilon}F_{H_{(2, \epsilon)}(t)}|_{H_{(1, \epsilon)}(t)}+4|\sqrt{-1}\Lambda_{\omega_\epsilon}F_{H_{(1, \epsilon)}(t)}|_{H_{(1, \epsilon)}(t)})|\tilde{T}_{\epsilon}(t)|_{H_{(1, \epsilon)}(t), \omega_\epsilon}^2\frac{\omega_\epsilon^n}{n!}dt\\
\leq&\hat{C}\sqrt{\tilde{\delta}},
\end{split}
\end{equation}
where $\hat{C}$ is a uniform constant.
For simplicity, set
\begin{equation}
\tilde{\Phi}_\epsilon(t)=\sqrt{-1}\Lambda_{\omega_\epsilon}(F_{H_{(2, \epsilon)}(t)}-F_{H_{(1, \epsilon)}(t)})=\sqrt{-1}\Lambda_{\omega_\epsilon}\bar{\partial}_E(\tilde{h}_{\epsilon }^{-1}(t)\partial_{H_{(1, \epsilon)}(t)}\tilde{h}_{\epsilon }(t))
\end{equation}
Using (\ref{C}), we obtain
\begin{equation}
\begin{split}
\frac{\partial}{\partial t}\tilde{\Phi}_\epsilon(t)
=&\sqrt{-1}\Lambda_{\omega_\epsilon}\bar{\partial}_E\big(\frac{\partial}{\partial t}\tilde{T}_{\epsilon}(t)\big)\\
=&\sqrt{-1}\Lambda_{\omega_\epsilon}\bar{\partial}_E(-2\partial_{H_{(1, \epsilon)}(t)}(\sqrt{-1}\Lambda_{\omega_\epsilon}\bar{\partial}_E\tilde{T}_{\epsilon}(t))-2\tilde{T}_{\epsilon}(t)\sqrt{-1}\Lambda_{\omega_\epsilon}F_{H_{(2, \epsilon)}(t)}\\
&+2\sqrt{-1}\Lambda_{\omega_\epsilon}F_{H_{(2, \epsilon)}(t)}\tilde{T}_{\epsilon}(t))\\
=&-2\sqrt{-1}\Lambda_{\omega_\epsilon}\bar{\partial}_E(\partial_{H_{(1, \epsilon)}(t)}(\tilde{\Phi}_\epsilon(t))+\tilde{T}_{\epsilon}(t)\sqrt{-1}\Lambda_{\omega_\epsilon}F_{H_{(2, \epsilon)}(t)}-\sqrt{-1}\Lambda_{\omega_\epsilon}F_{H_{(2, \epsilon)}(t)}\tilde{T}_{\epsilon}(t)).
\end{split}
\end{equation}
On the other hand, we have
\begin{equation}
\begin{split}
\Delta_{\omega_\epsilon}|\tilde{\Phi}_\epsilon(t)|_{H_{(1, \epsilon)}(t)}^2=&4Re\langle-\sqrt{-1}\Lambda_{\omega_\epsilon}\bar{\partial}_E\partial_{H_{(1, \epsilon)}(t)}\tilde{\Phi}_\epsilon(t), \tilde{\Phi}_\epsilon(t)\rangle_{H_{(1, \epsilon)}(t)}\\
&+2|\partial_{H_{(1, \epsilon)}(t)}\tilde{\Phi}_\epsilon(t)|_{H_{(1, \epsilon)}(t)}^2+2|\bar{\partial}_E\tilde{\Phi}_\epsilon(t)|_{H_{(1, \epsilon)}(t)}^2\\
&-2\sqrt{-1}\Lambda_{\omega_\epsilon}(\tilde{\Phi}_\epsilon(t)H_{(1, \epsilon)}^{-1}(t)\overline{[F_{H_{(1, \epsilon)}(t)}, \tilde{\Phi}_\epsilon(t)]}^\mathrm{t}H_{(1, \epsilon)}(t)).
\end{split}
\end{equation}
and then
\begin{equation}
\begin{split}
&(\Delta_{\omega_\epsilon}-\frac{\partial}{\partial t})|\tilde{\Phi}_\epsilon(t)|_{H_{(1, \epsilon)}(t)}^2\\
=&2|\partial_{H_{(1, \epsilon)}(t)}\tilde{\Phi}_\epsilon(t)|_{H_{(1, \epsilon)}(t)}^2+2|\bar{\partial}_E\tilde{\Phi}_\epsilon(t)|_{H_{(1, \epsilon)}(t)}^2\\
&-4Re \tr(\sqrt{-1}\Lambda_{\omega_\epsilon}\bar{\partial}_E([\sqrt{-1}\Lambda_{\omega_\epsilon}F_{H_{(2, \epsilon)}(t)}, \tilde{T}_{\epsilon}(t)])H_{(1, \epsilon)}^{-1}(t)\overline{\tilde{\Phi}_\epsilon}^\mathrm{t}(t)H_{(1, \epsilon)}(t))\\
&+4\langle\sqrt{-1}\Lambda_{\omega_\epsilon}[F_{H_{(1, \epsilon)}(t)}, \tilde{\Phi}_\epsilon(t)], \tilde{\Phi}_\epsilon(t)\rangle_{H_{(1, \epsilon)}(t)}.
\end{split}
\end{equation}
Integrating this over $\tilde{M}$ gives
\begin{equation}\label{C9}
\begin{split}
-\frac{\partial}{\partial t}\int_{\tilde{M}}|\tilde{\Phi}_\epsilon(t)|_{H_{(1, \epsilon)}(t)}^2\frac{\omega_\epsilon^n}{n!}=&2\int_{\tilde{M}}(|\partial_{H_{(1, \epsilon)}(t)}\tilde{\Phi}_\epsilon(t)|_{H_{(1, \epsilon)}(t)}^2
+|\bar{\partial}_E\tilde{\Phi}_\epsilon(t)|_{H_{(1, \epsilon)}(t)}^2)\frac{\omega_\epsilon^n}{n!}\\
&-4Re \int_{\tilde{M}}\langle[\sqrt{-1}\Lambda_{\omega_\epsilon}F_{H_{(2, \epsilon)}(t)}, \tilde{T}_{\epsilon}(t)], \partial_{H_{(1, \epsilon)}(t)}\tilde{\Phi}_\epsilon(t)\rangle_{H_{(1, \epsilon)}(t)}\frac{\omega_\epsilon^n}{n!}\\
&+4\int_{\tilde{M}}\langle\sqrt{-1}\Lambda_{\omega_\epsilon}[F_{H_{(1, \epsilon)}(t)}, \tilde{\Phi}_\epsilon(t)], \tilde{\Phi}_\epsilon(t)\rangle_{H_{(1, \epsilon)}(t)}\frac{\omega_\epsilon^n}{n!}\\
\geq&-8\int_{\tilde{M}}|\sqrt{-1}\Lambda_{\omega_\epsilon}F_{H_{(2, \epsilon)}(t)}|_{H_{(1, \epsilon)}(t)}^2|\tilde{T}_{\epsilon}(t)|_{H_{(1, \epsilon)}(t)}^2\frac{\omega_\epsilon^n}{n!}\\
&-8\int_{\tilde{M}}|\sqrt{-1}\Lambda_{\omega_\epsilon}F_{H_{(1, \epsilon)}(t)}|_{H_{(1, \epsilon)}(t)}|\tilde{\Phi}_\epsilon(t)|_{H_{(1, \epsilon)}(t)}^2\frac{\omega_\epsilon^n}{n!},
\end{split}
\end{equation}
and then
\begin{equation}
\begin{split}
-\frac{\partial}{\partial t}\Big(t^{2n}\int_{\tilde{M}}|\tilde{\Phi}_\epsilon(t)|_{H_{(1, \epsilon)}(t)}^2\frac{\omega_\epsilon^n}{n!}\Big)
\geq&-2nt^{2n-1}\int_{\tilde{M}}|\tilde{\Phi}_\epsilon(t)|_{H_{(1, \epsilon)}(t)}^2\frac{\omega_\epsilon^n}{n!}\\
&-8t^{2n}\int_{\tilde{M}}|\sqrt{-1}\Lambda_{\omega_\epsilon}F_{H_{(2, \epsilon)}(t)}|_{H_{(1, \epsilon)}(t)}^2|\tilde{T}_{\epsilon}(t)|_{H_{(1, \epsilon)}(t)}^2\frac{\omega_\epsilon^n}{n!}\\
&-8t^{2n}\int_{\tilde{M}}|\sqrt{-1}\Lambda_{\omega_\epsilon}F_{H_{(1, \epsilon)}(t)}|_{H_{(1, \epsilon)}(t)}|\tilde{\Phi}_\epsilon(t)|_{H_{(1, \epsilon)}(t)}^2\frac{\omega_\epsilon^n}{n!}.
\end{split}
\end{equation}
 By (\ref{C1}), (\ref{C2}), (\ref{C3}) and (\ref{C6}), we immediately get that
\begin{equation}\label{CON2}
\begin{split}
t^{2n}\int_{\tilde{M}}|\tilde{\Phi}_\epsilon(t)|_{H_{(1, \epsilon)}(t)}^2\frac{\omega_\epsilon^n}{n!}\leq &2n\int_0^1t^{2n-1}\int_{\tilde{M}}|\tilde{\Phi}_\epsilon(t)|_{H_{(1, \epsilon)}(t)}^2\frac{\omega_\epsilon^n}{n!}dt\\
&+8\int_0^1t^{2n}\int_{\tilde{M}}|\sqrt{-1}\Lambda_{\omega_\epsilon}F_{H_{(2, \epsilon)}(t)}|_{H_{(1, \epsilon)}(t)}^2|\tilde{T}_{\epsilon}(t)|_{H_{(1, \epsilon)}(t)}^2\frac{\omega_\epsilon^n}{n!}dt\\
&+8\int_0^1t^{2n}\int_{\tilde{M}}|\sqrt{-1}\Lambda_{\omega_\epsilon}F_{H_{(1, \epsilon)}(t)}|_{H_{(1, \epsilon)}(t)}|\tilde{\Phi}_\epsilon(t)|_{H_{(1, \epsilon)}(t)}^2\frac{\omega_\epsilon^n}{n!}dt\\
\leq&2n\int_0^1t^{n}\int_{\tilde{M}}|\tilde{\Phi}_\epsilon(t)|_{H_{(1, \epsilon)}(t)}^2\frac{\omega_\epsilon^n}{n!}dt
+8\hat{C}_2^2\int_0^1\int_{\tilde{M}}|\tilde{T}_{\epsilon}(t)|_{H_{(1, \epsilon)}(t)}^2\frac{\omega_\epsilon^n}{n!}dt\\\\
&+8\hat{C}_1\int_0^1t^n\int_{\tilde{M}}|\tilde{\Phi}_\epsilon(t)|_{H_{(1, \epsilon)}(t)}^2\frac{\omega_\epsilon^n}{n!}dt\\
\leq &\hat{C}_{4}\sqrt{\tilde{\delta}}
\end{split}
\end{equation}
for any $t\in [0, 1]$, where $\hat{C}_{4}$ is a uniform constant. For  any $t>1$, (\ref{C9}), (\ref{C4}), (\ref{C5}), (\ref{C6}) and (\ref{CON2}) imply that
\begin{equation}\label{CON3}
\int_{\tilde{M}}|\tilde{\Phi}_\epsilon(t)|_{H_{(1, \epsilon)}(t)}^2\frac{\omega_\epsilon^n}{n!}\leq \sqrt{\tilde{\delta}}\exp{(\hat{C}_{5}t)},
\end{equation}
where $\hat{C}_{5}$ is a uniform constant. Since $H_{i, \epsilon}(x, t)$  converges  to the long time solution $H_{i}(x, t)$  outside $\Sigma_{\mathcal{E}} $ in $C_{loc}^{\infty}$-topology as $\epsilon \rightarrow 0$,  (\ref{CON2}) and (\ref{CON3}) mean the inequality (\ref{CON1}).

\hfill $\Box$ \\

\medskip

From Lemma \ref{lemmaHYM}, we see that $t\mapsto HYM_{\alpha , N}(A(t), M, \omega )$ is
nonincreasing along the Yang-Mills flow (\ref{YM2}).
Note that Corollary \ref{Corollary 2.5.} says we can choose a sequence
$t_{j}\rightarrow +\infty$, such that
\begin{equation}
HYM_{\alpha , N}(A(t_{j}), M, \omega )\rightarrow HYM_{\alpha , N}(A_{\infty}, M, \omega ).
\end{equation}
 Then we have
\begin{equation}\label{Prop.5.3}\lim_{t\rightarrow \infty }HYM_{\alpha , N}(A(t), M, \omega )=HYM_{\alpha , N}(A_{\infty}, M, \omega )\end{equation} for any $\alpha \geq 1 $ and any $N$.
In the following we assume that $\Vol(M, \omega )=1$, and set $HYM_{\alpha , N}(\vec{\mu })=HYM_{\alpha
}(\vec{\mu }+N)= \varphi_{\alpha }(\sqrt{-1}(\vec{\mu
}+N))$, where $\vec{\mu }+N=diag (\mu_{1}+N , \cdots ,
\mu_{R}+N)$. Using Proposition \ref{prop 4.1}, Lemma \ref{continuous}, and following the argument in Theorem 4.1 in \cite{DW1}, we can obtain $\vec{\lambda }_{\infty}=\vec{\mu}_{0}$. We give a proof briefly for readers' convenience.

\hspace{0.3cm}

\begin{theorem}\label{thm key} Let $\mathcal{E}$ be a
reflexive sheaf on a smooth K\"ahler manifold $(M , \omega )$,  $H(t)$ be a solution of the Hermitian-Yang-Mills flow (\ref{SSS1}) on $\mathcal{E}|_{M\setminus \Sigma_{\mathcal{E}}}$ with the initial metric $\hat{H}$, and $A(t) $ be  the related Yang-Mills flow (\ref{YM2}) on the Hermitian vector bundle $(\mathcal{E}|_{M\setminus \Sigma_{\mathcal{E}}}, \hat{H})$. Assume that $A_{\infty }$ is an Uhlenbeck limit of $A(t)$, and $(E_{\infty}, H_{\infty})$ is the corresponding Hermitian vector bundle defined on $M\setminus (\Sigma_{\mathcal{E}}\cup \Sigma_{an})$.  Then there is a constant $p_{0}>1$ such that
\begin{equation}\label{CC1}
HYM_{\alpha , N} (A_{\infty },  M, \omega )=\lim_{t\rightarrow \infty  }HYM_{\alpha , N} (A(t), M, \omega)=HYM_{\alpha , N} (\vec{\mu}_{0}),
\end{equation}
for all $1\leq \alpha < p_{0}$ and all $N\in \mathbb{R}$; and the HN type of the reflexive sheaf $(E_{\infty}, \overline{\partial}_{A_{\infty}})$ is the same as that of $\mathcal{E}$, i.e. $\vec{\lambda }_{\infty}=\vec{\mu}_{0}$.
\end{theorem}

\medskip

{\bf Proof }  As before, let $\pi : \tilde{M}\rightarrow M $ be the composition of  a finite sequence of blowups resolving the sheaf $\mathcal{E}$, i.e. such that $E=\pi^\ast\mathcal{E}/tor(\pi^\ast\mathcal{E})$  is locally free. Firstly, since the norm $(\int_{M}
\varphi_{\alpha }(\verb"a") dvol )^{\frac{1}{\alpha }}$ is equivalent to the
$L^{\alpha }$-norm on $\verb"u"(E)$,  we have
\begin{equation*}
\begin{split}
& |(HYM_{\alpha , N}((\mathcal{E},
H), M, \omega))^{\frac{1}{\alpha }}-(HYM_{\alpha , N}(\vec{\mu }_{0}))^{^{\frac{1}{\alpha }}}|\\
&\leq \Big(\int_{M} |(\varphi_{\alpha }(\frac{\sqrt{-1}}{2\pi }\Lambda_{\omega
}(F_{(\mathcal{E},
H)})+N \Id_{E}))^{\frac{1}{\alpha}}-(\varphi_{\alpha }(\vec{\mu}+N))^{\frac{1}{\alpha}}|^{\alpha }\frac{\omega ^{n}}{n!}\Big)^{\frac{1}{\alpha }}\\
&\leq  \Big(\int_{M} \varphi_{\alpha }(\frac{\sqrt{-1}}{2\pi }\Lambda_{\omega
}(F_{(\mathcal{E},
H)})- \Psi(\textit{F}, (\mu_{ 1} , \cdots , \mu_{l}) ,
H))\frac{\omega ^{n}}{n!} \Big)^{\frac{1}{\alpha}}\\
&\leq
C(\alpha )\|\frac{\sqrt{-1}}{2\pi }\Lambda_{\omega
}(F_{(\mathcal{E},
H)})-\Psi(\mathcal{F}, (\mu_{ 1, \omega} , \cdots , \mu_{l, \omega }) ,
H)\|_{L^{\alpha }(M, \omega)}.
\end{split}
\end{equation*}
This together with Proposition \ref{prop 4.1} gives us that for any $\delta >0$ and any $1\leq \alpha  < p_{0}$ there is $H$ on the bundle $E$ such that
\begin{equation}\label{BB2}
HYM_{\alpha , N}((\mathcal{E},
H), M, \omega)\leq  HYM_{\alpha , N}(\vec{\mu }_{0}) +\delta .
\end{equation}

 For fixed $\alpha $
and  $N$, since the image of the degree map on line bundles is discrete, we can  define $\delta_{0}>0$ such that
\begin{equation}
2\delta_{0}+HYM_{\alpha , N}(\vec{\mu}_{0})=\min\{HYM_{\alpha ,
N}(\vec{\mu}): HYM_{\alpha , N}(\vec{\mu})>HYM_{\alpha ,
N}(\vec{\mu}_{0})\},
\end{equation}
where $\vec{\mu }$ runs over all possible HN types of torsion-free sheaves
on $M$ with the rank of $\mathcal{E}$.

\medskip

Let $\hat{H}$ be a smooth Hermitian metric on the holomorphic vector bundle $E$, $H(t)$ be the solution of the Hermitian-Yang-Mills flow (\ref{SSS1}) on $\mathcal{E}|_{M\setminus \Sigma_{\mathcal{E}}}$ with the initial metric $\hat{H}$ and
$A^{\hat{H}}(t) $ be the solution of the related Yang-Mills
flow (\ref{YM2}) on the Hermitian vector bundle $(\mathcal{E}|_{M\setminus \Sigma_{\mathcal{E}}}, \hat{H})$ with the initial
connection $\hat{A}=(\mathcal{E}, \hat{H})$. Let $A_{\infty }^{\hat{H}} $ be an Uhlenbeck limit along
the Yang-Mills flow (\ref{YM2}).
Assume  that $\hat{H}_{0}$ satisfies:
\begin{equation}\label{5.12}
HYM_{\alpha , N}((\mathcal{E}, \hat{H}_{0}), M, \omega )\leq  HYM_{\alpha , N}(\vec{\mu }_{0}) +\delta_{0} .
\end{equation}
Combining (\ref{Prop.5.3}), Lemma \ref{lemmaHYM} and
(\ref{P5.5}), we obtain:
\begin{equation*}
HYM_{\alpha , N}(\vec{\mu}_{0})\leq HYM_{\alpha , N}(A_{\infty }^{\hat{H}_{0}} ,
M, \omega )\leq
HYM_{\alpha , N}(\vec{\mu}_{0})+\delta_{0} .
\end{equation*}
Hence we must have $HYM_{\alpha , N}(A_{\infty}^{\hat{H}_{0}} , M, \omega) =
HYM_{\alpha , N}(\vec{\mu}_{0})$. This shows that the result holds
if the metric $\hat{H}_{0}$ satisfies (\ref{5.12}).

For any fixed $\delta \leq \frac{\delta_{0}}{2}$,  we denote by $\textbf{H}_{\delta }$  the set of smooth Hermitian
metrics on $E$ satisfying that,  there is
$T\geq 0$ such that
\begin{equation}\label{ineq1}
HYM_{\alpha , N}(A^{\hat{H}}(t) , M, \omega) < HYM_{\alpha ,
N}(\vec{\mu}_{0})+\delta ,
\end{equation}
for all $t\geq T$.  From  (\ref{BB2}) and the discussion above, we see
$\textbf{H}_{\delta }$ is nonempty. In Lemma \ref{continuous}, we have proved the continuous dependence of the Hermitian-Yang-Mills flow (\ref{SSS1}) on initial metrics, this implies the openness of $\textbf{H}_{\delta }$. By Lemma \ref{lemmaBS} and (\ref{interg}), $\|\Lambda_{\omega}F_{A^{\hat{H}}(t)}\|_{L^{\infty}}$ and $\|F_{A^{\hat{H}}(t)}\|_{L^{2}}$ are uniformly bounded along the Yang-Mills flow (\ref{YM2}) for $t\geq t_{0} >0$. On the other hand, the Uhlenbeck compactness theorem (Theorem 5.2 in \cite{UY}) is also valid for the non-compact case, i.e. on the non-compact K\"ahler manifold $(M\setminus \Sigma_{\mathcal{E}}, \omega )$. So we can follow the argument in Lemma 4.3 in \cite{DW1} to show that $\textbf{H}_{\delta }$ is closed. The proof is exactly the same, we omit it.  Since the
space of smooth metrics on $E$ is connected, we conclude that every
metric is in $\textbf{H}_{\delta }$. Then it follows that $HYM_{\alpha , N}(A_{\infty}^{\hat{H}} , M, \omega)=\lim_{t\rightarrow +\infty }HYM_{\alpha , N}(A^{\hat{H}}(t) ,
M, \omega ) =HYM_{\alpha , N}(\vec{\mu}_{0})$ for any metric
$\hat{H}$ on $E$. With Proposition 2.24 in \cite{DW1}, we know $\vec{\mu}_{0} = \vec{\lambda }_{\infty }$. This concludes the proof of Theorem \ref{thm key}.

\hfill $\Box$ \\

\medskip

 Let $H(t)$ be the long time solution of the Hermitian-Yang-Mills flow (\ref{SSS1}) with the initial metric $\hat{H}$, and $A(t)$ be the solution of the related Yang-Mills flow (\ref{YM2}) with the initial connection $\hat{A}$.
As that in Proposition \ref{propgauge}, we have $A(t)=\sigma (t)(\hat{A})$, where $\sigma (t)$ is a family of complex gauge transformations  satisfying
$\sigma ^{\ast \hat{H}}(t)\sigma (t)=h (t)=\hat{H}^{-1}H(t)$.
Consider  the following  HN-filtration of $\mathcal{E}$ by saturated sheaves
\begin{equation}
0=\mathcal{E}_{0}\subset \mathcal{E}_{1}\subset \cdots \subset \mathcal{E}_{\tilde{k}}=\mathcal{E}.
\end{equation}
  Let $\pi _{\alpha}^{H(t)}$ be the orthogonal projection onto $\mathcal{E}_{\alpha}$ with respect to the Hermitian metric $H(t)$, and $\pi _{\alpha}^{(t)}=\sigma (t)\circ \pi _{\alpha}^{H(t)} \circ \sigma ^{-1}(t)$.   It is easy to check that:
 $(\Id -\pi_{\alpha }^{(t)}) \overline{\partial }_{A(t)}\pi_{\alpha }^{(t)} =0$;
$ (\pi_{\alpha }^{(t)})^{2}=\pi_{\alpha }^{(t)}=(\pi_{\alpha }^{(t)})^{\ast \hat{H}}$, $|\overline{\partial }_{A(t)}\pi_{\alpha }^{(t)}|_{\hat{H}}=|\overline{\partial }_{\hat{A}}\pi_{\alpha }^{H(t)}|_{H(t)}$. From (\ref{deg}), it can be seen that $\pi_{\alpha }^{(t)} \in L_{1}^{2}(\End (\mathcal{E}))$.
  Using  Theorem \ref{thm key} and following  the same argument in \cite{DW1} (Proposition 4.5), we deduce the following lemma.

\medskip

\begin{lemma}\label{lemma3.13} Let $\mathcal{E}$ be a reflexive sheaf on a smooth K\"ahler manifold $(M , \omega )$, and satisfy the same assumptions as that in Theorem \ref{thm key}. Assume that $A_{\infty }$ is an Uhlenbeck limit of $A(t)$, and $(E_{\infty}, H_{\infty})$ is the corresponding Hermitian vector bundle defined on $M\setminus (\Sigma_{\mathcal{E}}\cup \Sigma_{an})$.

(1) Let $\{\pi_{\alpha}^{\infty}\}$ be the HN-filtration of the reflexive sheaf $(E_{\infty}, \overline{\partial}_{A_{\infty}})$, then there is a sequence of  $\{ \pi_{\alpha }^{(t_j)} \}$ which converges to  $\{ \pi_{\alpha }^{\infty } \}$  strongly in $L^{p}\cap L_{1, loc}^{2}$ outside $\Sigma_{\mathcal{E}}\cup\Sigma_{an}$ as $j$ tends to $+\infty$.

(2) Assume the sheaf $\mathcal{E}$ is semi-stable and $\{\mathcal{E}_{\alpha }\}$ is the Seshadri filtration of $\mathcal{E}$, then $\{ \pi_{\alpha }^{(t_j)} \}$ converges to a filtration $\{ \pi_{\alpha }^{\infty } \}$  strongly in $L^{p}\cap L_{1, loc}^{2}$ outside $\Sigma_{\mathcal{E}}\cup\Sigma_{an}$ as $j$ tends to $+\infty$
, the rank and degree of $ \pi_{\alpha }^{\infty } $ is equal to the rank and degree of $ \pi_{\alpha }^{t_j } $ for all $\alpha $ and $j$.

\end{lemma}

\medskip

\section{Proof of theorem 1.1. }
\setcounter{equation}{0}

 In this section, we will prove the part (2) of Theorem 1.1
inductively on the length of the HNS-filtration. The inductive
hypotheses are following:

\hspace{0.3cm}

{\bf Inductive hypotheses:} {\it  Let $\mathcal{Q}$ be a torsion-free sheaf on a compact K\"ahler manifold $(M, \omega )$,  $\mathcal{S}$ be a saturated sub-sheaf of $\mathcal{Q}$.

(1) There is a sequence of connections $A_{ j}^{Q} \in \textbf{A}^{1,1}_{Q, H_{0}} $ on the Hermitian bundle $ (\mathcal{Q}|_{M\setminus \tilde{\Sigma}}, H_{0})$ such that $A_{ j}^{Q}
\rightarrow A_{ \infty}^{Q_{\infty }}$ in
$C^{\infty }_{loc}$-topology off $\tilde{\Sigma}$ as $j\rightarrow +\infty$, where $\tilde{\Sigma}$ is a complex analytic subset of $M$ with complex codimension at least $2$ and satisfies $\Sigma_{\mathcal{Q}}\cup \Sigma_{\mathcal{S}}\subset \tilde{\Sigma }$.

(2) $A_{ j}^{Q} =g_{j}(A_{H_0}^{Q} )$
 for some complex gauge transformations
$g_{j}\in  \textbf{G}^{\mathbb{C}}(Q)$ and $\|\sqrt{-1}\Lambda _{\omega }(F_{A_{H_j}^{Q}})\|_{L^{1}(\omega)}$ is uniformly bounded  in $j$,  where $A_{H_j}^{Q}$ is the Chern connection on $\mathcal{Q}$ with respect to the metric $H_{j}=H_{0}g_{j}^{\ast H_{0}} g_{j}$.

(3)  There exists a sequence of blow-ups with smooth center:
$
\pi_{i}: M_{i}\rightarrow M_{i-1},$ $ i=1, \cdots , r,
$
and an exact sequence of holomorphic vector bundles
\begin{equation}
0\rightarrow \widetilde{\mathcal{S}}\rightarrow \widetilde{\mathcal{Q}}\rightarrow \widetilde{\mathcal{W}}\rightarrow 0
\end{equation}
over $\tilde{M}$, such that the composition $\pi =\pi_{r}\circ \cdots \circ \pi_{1}: \tilde{M}\rightarrow M$ is biholomorphic outside $\tilde{\Sigma} $, $\widetilde{\mathcal{S}} $ and $\widetilde{\mathcal{Q}}$ are isomorphic to $\mathcal{S}$ and $\mathcal{Q}$ outside $\tilde{\Sigma}$ respectively, where $M_{0}=M$, $\tilde{M}=M_{r}$.

(4) Set $\epsilon =(\epsilon_{1}, \cdots , \epsilon_{r})$ and define K\"ahler metrics $\omega_{\epsilon}$ on $M_{r}$ as that in (\ref{omegai}). For every $j$, there exists a sequence of metrics $H_{j, \epsilon}$ on $\widetilde{\mathcal{Q}}$ such that $H_{j, \epsilon}\rightarrow H_{j}$ in $C_{loc}^{\infty}$-topology outside $\tilde{\Sigma}$ as $\epsilon \rightarrow 0$,  $\|\sqrt{-1}\Lambda _{\omega_{\epsilon} }(F_{A_{H_{j , \epsilon }}^{\tilde{Q}}})\|_{L^{1}(\omega_{\epsilon })}$ is uniformly bounded, and $\sup_{M_{r}}(\tr (H_{1, \epsilon}^{-1}H_{j, \epsilon})+\tr (H_{j, \epsilon}^{-1}H_{1, \epsilon}))<C_{j}$, where $C_{j}$ is a constant  independent of $\epsilon$. Furthermore, $\|\sqrt{-1}\Lambda _{\omega_{\epsilon} }(F_{A_{H_{1 , \epsilon }}^{\tilde{S}}})\|_{L^{1}(\omega_{\epsilon })}$ is uniformly bounded, where $A_{H_{1 , \epsilon }}^{\tilde{S}}$ is the induced Chern connection on $\widetilde{\mathcal{S}}$.

(5) Two torsion-free sheaves $\mathcal{Q}$ and $(Q_{\infty }, \overline{\partial }_{A_{\infty}^{Q_{\infty}}})$ have the same HN type.
 }

\hspace{0.3cm}

Now we construct non-zero holomorphic maps from subsheaves in the HNS-filtration of $\mathcal{E}$ to the limiting  reflexive sheaf $(E_{\infty}, \overline{\partial }_{\infty})$.  We get a nonzero holomorphic map which we need by limiting a sequence of holomorphic maps. The key problem is to obtain local uniform $C^{0}$-estimate of this sequence of holomorphic maps. We will follow the argument in Proposition 4.1  in \cite{LZ3}   to handle this problem. There is a difference in the assumption for our case, so we write a proof briefly of the following proposition for readers' convenience.

\medskip

\begin{proposition}\label{propmap}
Let $\mathcal{Q}$ be a torsion-free sheaf on a compact K\"ahler manifold $(M, \omega )$,  $\mathcal{S}$ be a saturated sub-sheaf of $\mathcal{Q}$. Assume that the conditions (1), (2), (3), (4) in the above inductive hypotheses are satisfied. Let $i_{0}: \mathcal{S}\rightarrow \mathcal{Q}$ be the holomorphic inclusion, then there is a subsequence of $g_{j} \circ i_{0}$, up to rescale,  converges to a non-zero holomorphic map $f_{\infty }: \mathcal{S}\rightarrow (Q_{\infty }, \overline{\partial}_{A_{\infty}^{Q}})$ in $C^{\infty}_{loc}$-topology off $\tilde{\Sigma}$ as $j\rightarrow +\infty$.
\end{proposition}

\medskip

{\bf Proof. }By induction, we can assume that $\pi : \tilde{M}\rightarrow M$ is a single blow-up with smooth centre. Fix a K\"ahler metric $\eta $ on $\tilde{M}$ and set $\omega_{\epsilon }=\pi^{\ast}\omega +\epsilon \eta $
  for $0<\epsilon \leq 1$.
 On the blow-up $\tilde{M}$, let $H_{j, \epsilon } (t)$ and $H^{S}_{1, \epsilon }(t)$ be the solutions of the following Hermitian-Yang-Mills flow on holomorphic bundles $\widetilde{\mathcal{Q}}$ and $\widetilde{\mathcal{S}}$ with the fixed initial metrics $H_{j, \epsilon }$ and $H^{S}_{1, \epsilon }$  and with respect to the metric $\omega_{\epsilon}$, i.e. they satisfy the following heat equation
\begin{equation}\label{0DDD1}
H^{-1}(t)\frac{\partial H(t)}{\partial
t}=-2\sqrt{-1}\Lambda_{\omega_{\epsilon}}F_{H(t)},
\end{equation}
where $H_{j, \epsilon}$ is defined in condition (4) among the inductive hypotheses.
A direct computation yields
\begin{equation}\label{H00}
(\Delta_{\epsilon } -\frac{\partial }{\partial t}) |\Lambda _{\omega_{\epsilon }}(F_{H_{j, \epsilon}(t)}) |_{H_{j, \epsilon}(t)}\geq  0,
\end{equation}
\begin{equation}\label{H00S}
(\Delta_{\epsilon } -\frac{\partial }{\partial t}) |\Lambda _{\omega_{\epsilon }}(F_{H^{S}_{1, \epsilon}(t)}) |_{H^{S}_{1, \epsilon}(t)}\geq  0,
\end{equation}
and
\begin{equation}
(\Delta_{\epsilon }-\frac{\partial }{\partial t}) |i_{0}|_{H^{S}_{1, \epsilon}(t), H_{j, \epsilon } (t)}^{2}\geq 0.
\end{equation}
 The maximum principle implies that
 \begin{equation}\label{0H000}
 |\Lambda _{\omega_{\epsilon }}(F_{H_{j, \epsilon}(t)}) |_{H_{j, \epsilon}(t)}(x) \leq \int_{\tilde{M}}K_{\epsilon} (t-t_{0}, x, y)|\Lambda _{\omega_{\epsilon }}(F_{H_{j, \epsilon}(t_{0})}) |_{H_{j, \epsilon}(t_{0})}\frac{\omega_{\epsilon }^{n}}{n!} ,
\end{equation}
\begin{equation}\label{H000S}
 |\Lambda _{\omega_{\epsilon }}(F_{H^{S}_{1, \epsilon}(t)}) |_{H^{S}_{1, \epsilon}(t)}(x) \leq \int_{\tilde{M}}K_{\epsilon} (t-t_{0}, x, y)|\Lambda _{\omega_{\epsilon }}(F_{H^{S}_{1, \epsilon}(t_{0})}) |_{H^{S}_{1, \epsilon}(t_{0})}\frac{\omega_{\epsilon }^{n}}{n!} ,
\end{equation}
and
 \begin{equation}\label{H0001}
 |i_{0}|^{2}_{H^{S}_{1, \epsilon}(t), H_{j, \epsilon } (t_{0}+t)}(x) \leq \int_{\tilde{M}}K_{\epsilon} (t-t_{0}, x, y)|i_{0}|^{2}_{H^{S}_{1, \epsilon}(t_{0}), H_{j, \epsilon } (t_{0})}\frac{\omega_{\epsilon }^{n}}{n!} ,
\end{equation}
for any $t>t_{0}\geq 0$.
By \cite{BS} (Lemma 4), the heat kernels $K_{\epsilon} (t, x, y) $ have a uniform bound for $0<\epsilon \leq 1$.  Following Bando and Siu's argument (\cite{BS}), we could choose a subsequence of $H_{j, \epsilon } (t)$ (and the same for $H^{S}_{1, \epsilon}(t)$) which converges to $H_{j}(t)$ (resp. $H^{S}_{1}(t)$) a solution of the Hermitian-Yang-Mills flow (\ref{0DDD1}) on $\mathcal{Q}$ (resp. $\mathcal{S}$) over $M\setminus \tilde{\Sigma }$ as $\epsilon $ tends to $0$. Combining (\ref{0H000}), (\ref{H000S}), (\ref{H0001}) and the condition (4), we derive
\begin{equation}\label{FFF1}
2( |\Lambda _{\omega}(F_{H_{j}(t)}) |_{H_{j}(t)}+  |\Lambda _{\omega}(F_{H^{S}_{1}(t)}) |_{H^{S}_{1}(t)}) (x)\leq C_{F}
\end{equation}
and
 \begin{equation}\label{H0005}
 |i_{0}|^{2}_{H^{S}_{1}(t_{0}+t), H_{j } (t)}(x) \leq \int_{M}K (t-t_{0}, x, y)|i_{0}|_{H^{S}_{1}(t_{0}), H_{j} (t_{0})}^{2}\frac{\omega ^{n}}{n!} ,
\end{equation}
for all $x$ outside $\tilde{\Sigma}$ and $t\geq t'>t_{0}\geq 0$, where $K(t, x, y)$ is the heat kernel of $(M, \omega )$ and $C_{F}$ is a uniform constant which is independent of $j$.

From (\ref{FFF1}), it follows that
\begin{equation}
\Big|\frac{\partial }{\partial t} \ln  |i_{0}|^{2}_{H^{S}_{1}(t), H_{j } (t)}(x)\Big|\leq 2( |\Lambda _{\omega}(F_{H_{j}(t)}) |_{H_{j}(t)}+  |\Lambda _{\omega}(F_{H^{S}_{1}(t)}) |_{H^{S}_{1}(t)}) (x)\leq C_{F} ,
\end{equation}
for all $x\in M\setminus \tilde{\Sigma } $ and $t\geq t'>0$. Then
\begin{equation}
e^{-C_{F}\delta }\leq \frac{|i_{0}|^{2}_{H^{S}_{1}(t'+\delta ), H_{j } (t'+\delta)}}{|i_{0}|^{2}_{H^{S}_{1}(t'), H_{j } (t')}}(x)\leq e^{C_{F}\delta },
\end{equation}
and
\begin{equation}\label{key11}
\begin{split}
&|i_{0}|^{2}_{H^{S}_{1}(t'), H_{j } (t')}(x)\leq e^{C_{F}\delta } |i_{0}|^{2}_{H^{S}_{1}(t'+\delta ), H_{j } (t'+\delta)}(x)\\
& \leq e^{C_{F}\delta } \int_{M}K (\delta , x, y)|i_{0}|^{2}_{H^{S}_{1}(t'), H_{j } (t')}\frac{\omega ^{n}}{n!}\\
& \leq C_{K}e^{C_{F}\delta } (1+\delta^{-n})\int_{M}|i_{0}|^{2}_{H^{S}_{1}(t'), H_{j } (t')}\frac{\omega ^{n}}{n!},\\
\end{split}
\end{equation}
for all $x\in M\setminus \tilde{\Sigma } $ and $\delta >0$.

 Denote $h_{j, \epsilon }(t)=H_{j, \epsilon }^{-1}H_{j, \epsilon }(t)$, and then the heat equation (\ref{0DDD1}) yields
\begin{equation}\label{H01}
 (\Delta_{\epsilon } -\frac{\partial }{\partial t}) \ln (\tr (h_{j, \epsilon }(t))+ \tr (h_{j, \epsilon }^{-1}(t)))
\geq  -2 |\Lambda _{\omega_{\epsilon }}(F_{H_{j, \epsilon }}) |_{H_{j, \epsilon }}.
\end{equation}
Integrating  the above inequality and using the condition (4), we have
\begin{equation}\label{H0002}
 \int_{\tilde{M}} \ln (\tr (h_{j, \epsilon }(t))+ \tr (h_{j, \epsilon }^{-1}(t)))\frac{\omega_{\epsilon }^{n}}{n!}-\ln 2\rank(\widetilde{\mathcal{Q}})\Vol (\tilde{M}, \omega_{\epsilon })\leq  tC_{h}
\end{equation}
and then
\begin{equation}\label{H02}
 \int_{M} \ln (\tr (h_{j }(t))+ \tr (h_{j }^{-1}(t)))\frac{\omega^{n}}{n!}-\ln 2\rank(\mathcal{Q})\Vol (M, \omega)
\leq  t C_{h}.
\end{equation}
On the other hand, it holds that
\begin{equation}\label{H03}
\begin{split}
&\Delta \ln (\tr (h_{j }(t))+ \tr (h_{j }^{-1}(t)))\\
&\geq -2 |\Lambda _{\omega }(F_{H_{j }(t)}) |_{H_{j }(t)} -2 |\Lambda _{\omega }(F_{H_{j}}) |_{H_{j}}\\
\end{split}
\end{equation}
on $M\setminus\tilde{\Sigma}$, for all $t>0$. Here, we should note that $|\Lambda _{\omega }(F_{H_{j}}) |_{H_{j}}=|\Lambda _{\omega }(F_{A_{j}^{Q}})|_{H_{0}}$.

For any compact subset $\Omega \subset M \setminus \tilde{\Sigma }$, the condition (1) implies that $|\Lambda _{\omega }(F_{H_{j}}) |_{H_{j}}=|\Lambda _{\omega }(F_{A_{j}^{Q}})|_{H_{0}}$ is uniformly bounded on $\Omega$. By (\ref{H02}), (\ref{H03}), (\ref{FFF1}) and the Moser's iteration, there must exist a uniform constant $C_{\Omega, Q}$ such that, for all $j$,
\begin{equation}\label{mc0}
 \sup_{x\in \Omega } \ln (\tr (h_{j }(1))+ \tr (h_{j }^{-1}(1)))\leq C_{\Omega , Q}.
\end{equation}

\medskip

Define the holomorphic map $\tilde{\eta}_{j}: (\mathcal{S}|_{M\setminus \tilde{\Sigma}}, \overline{\partial }_{A_{0}^{S}})\rightarrow (\mathcal{Q}|_{M\setminus \tilde{\Sigma}}, \overline{\partial }_{A_{j}^{Q}})$ by $\tilde{\eta}_{j}=g_{j} \circ i_{0}$, where $A_{0}^{S}$ is the induced connection on $\mathcal{S}$ by the connection $A^{Q}_{H_{0}}$ . It is easy to check that
\begin{equation}
|\tilde{\eta}_{j}|_{H_{1}^{S}, H_{0}}=|i_{0}|_{H_{1}^{S}, H_{j}},
\end{equation}
where $H_{1}^{S}$ is the induced metric on $\mathcal{S}$ by the metric $H_{1}$. Set
\begin{equation}
f_{j}= \Big(\int_{M}|i_{0}|^{2}_{H^{S}_{1}(1), H_{j } (1)}\frac{\omega ^{n}}{n!}\Big)^{-\frac{1}{2}} \tilde{\eta}_{j}.
\end{equation}

Clearly (\ref{key11}) means that there is a constant $C_{a}$ such that
\begin{equation}\label{key110}
\sup_{x\in M\setminus \tilde{\Sigma }}\Big(\int_{M}|i_{0}|^{2}_{H^{S}_{1}(1), H_{j } (1)}\frac{\omega ^{n}}{n!}\Big)^{-1}|i_{0}|^{2}_{H^{S}_{1}(1), H_{j } (1)}(x)\leq C_{a},
\end{equation}
 for all $j$.
 Using (\ref{key110}) and (\ref{mc0}), we obtain a local uniform $C^{0}$-estimate on $f_{j}$, i.e. for any compact subset $\Omega \subset M \setminus \tilde{\Sigma }$, there is a constant $C_{\omega , f}$ such that
\begin{equation}\label{C0111}
\sup_{x\in \Omega }|f_{j}|_{H^{S}_{1}, H_{0}} (x) \leq C_{\omega , f}
\end{equation}
for all $j$. By the above local uniform $C^{0}$-bound of $f_{j}$  and  the assumption that $A_{j}\rightarrow A_{\infty}$ in $C^{\infty}_{loc}$-topology outside $\tilde{\Sigma}$ as $j\rightarrow +\infty$, the elliptic theory implies that there exists a subsequence of $f_{j}$ (for simplicity, also denoted by $f_{j}$) such that $f_{j}$ converges to a holomorphic map $ f_{\infty}: \mathcal{S}\rightarrow (Q_{\infty }, \overline{\partial}_{A_{\infty}^{Q}})$ in $C^{\infty}_{loc}$-topology outside $\tilde{\Sigma}$ as $j\rightarrow +\infty$.
Now we only need to prove that $f_{\infty}$ is non-zero. Since  $\tilde{\Sigma}$ is of Hausdorff complex codimension at least $2$, for any small $\delta >0$, we can choose a compact subset $\Omega_{\delta} \subset M \setminus \tilde{\Sigma}$ such that
\begin{equation}\label{dd}
\int_{M\setminus \Omega_{\delta}} 1 \frac{\omega^{n}}{n!}\leq \delta .
\end{equation}
Of course the local uniform estimate (\ref{mc0}) gives us that there is a positive constant $C_{\delta }$ such that
\begin{equation}
 C_{\delta} |i_{0}|^{2}_{H^{S}_{1}(1), H_{j }(1) }(x)\leq   |i_{0}|^{2}_{H^{S}_{1}, H_{j } }(x)\leq C_{\delta}^{-1}|i_{0}|^{2}_{H^{S}_{1}(1), H_{j } (1)}(x),
\end{equation}
for all $x\in \Omega_{\delta }$ and $j$. Then
\begin{equation}
\begin{split}
\int_{M\setminus \Omega_{\delta}} |f_{\infty}|_{H^{S}_{1}, H_{0}} \frac{\omega^{n}}{n!}&= \lim_{j\rightarrow +\infty }\int_{M\setminus \Omega_{\delta}} |f_{j}|_{H^{S}_{1}, H_{0}} \frac{\omega^{n}}{n!}\\
&=\lim_{j\rightarrow +\infty}\Big(\int_{M}|i_{0}|^{2}_{H^{S}_{1}(1), H_{j } (1)}\frac{\omega ^{n}}{n!}\Big)^{-1}\int_{M\setminus \Omega_{\delta}} |i_{0}|_{H^{S}_{1}, H_{j}} \frac{\omega^{n}}{n!}\\
&\geq \lim_{j\rightarrow +\infty} C_{\delta}\Big(\int_{M}|i_{0}|^{2}_{H^{S}_{1}(1), H_{j } (1)}\frac{\omega ^{n}}{n!}\Big)^{-1}\int_{M\setminus \Omega_{\delta}} |i_{0}|_{H^{S}_{1}(1), H_{j}(1)} \frac{\omega^{n}}{n!}\\
&\geq C_{\delta }(1-\delta C_{a}) >0.\\
\end{split}
\end{equation}
Therefore $f_{\infty}$ is a non-zero holomorphic map. This concludes the proof of Proposition \ref{propmap}.

\hfill $\Box$ \\

\medskip

{\bf A proof of Theorem \ref{thm 1.1}} Let $\{\mathcal{E}_{\alpha }\}_{\alpha =1}^{l}$ be
the Harder-Narasimhan-Seshadri filtration  of $\mathcal{E}$,  $Gr^{HNS}(\mathcal{E})=\oplus_{\alpha =1}^{l}\mathcal{Q}_{\alpha }$ be the associated graded object, where $\mathcal{Q}_{\alpha }=\mathcal{E}_{\alpha }/\mathcal{E}_{\alpha -1}$ is torsion-free for each $1\leq \alpha \leq l$.
 We refer to $\Sigma_{HNS}$ as the singularity set of the HNS-filtration, it is a complex analytic subset of $M$ with complex codimension at least $2$.

According to Hironaka's flattening theorem (\cite{Hi2}), there is a finite sequence of blowing ups $\{\pi_{i}\}_{i=1}^{k}$ along compact sub-manifolds such that if we denote by $\pi: \tilde{M}\rightarrow M$ the composition of all the blowing ups, then $E=\pi^{\ast}\mathcal{E}/tor{(\pi^{\ast
}\mathcal{E})}$ is locally free. By Proposition \ref{claim},
 we can get a filtration  $\widetilde{\mathcal{F}}=\{\widetilde{\mathcal{E}_{i}}\}_{i=1}^{l}$ of $E$ :
\begin{equation}\label{0HNS02}
0= \widetilde{\mathcal{E}_0}\subset \widetilde{\mathcal{E}_1}\subset \cdots \subset\widetilde{\mathcal{E}}_{l-1}\subset\widetilde{\mathcal{E}_l}=E,
\end{equation}
such that, for every $1\leq \alpha \leq l$, $\widetilde{\mathcal{E}_i}$ is a reflexive sheaf, $\widetilde{\mathcal{Q}}_{\alpha}=\widetilde{\mathcal{E}}_{\alpha}/\widetilde{\mathcal{E}}_{\alpha -1}$ is torsion free and isomorphic to the sheaf $\mathcal{Q}_{\alpha} $  outside $\pi^{-1}(\Sigma_{HNS})$. By Sibley's result on the resolution of filtration (Proposition 4.3 in \cite{Sib}), there is a finite sequence of blowing ups along complex submanifolds whose composition $\hat{\sigma }: \hat{M}\rightarrow \tilde{M}$ enjoys the following properties.
There is a filtration
\begin{equation}\label{HNS022}
0= E_0\subset E_1\subset \cdots \subset E_{l-1}\subset E_l=\hat{E}=\hat{\sigma }^*E
\end{equation}
by subbundles. If we write $Im\hat{\sigma}^*\widetilde{\mathcal{E}}_{\alpha}$ for the image of $\hat{\sigma}^*\widetilde{\mathcal{E}}_{\alpha}\rightarrow \hat{\sigma}^*E$, then $E_i=Sat_{\hat{\sigma}^*E}(Im\hat{\sigma}^*\widetilde{\mathcal{E}}_{\alpha})$. If $Q_{\alpha}=E_{\alpha}/E_{\alpha-1}$, then we have $\hat{\sigma}_{\ast}E_{\alpha}=\widetilde{\mathcal{E}}_{\alpha}$ and $\widetilde{\mathcal{Q}}_{\alpha}^{\ast\ast}=(\hat{\sigma}_{\ast}Q_{\alpha})^{\ast\ast}$.
Now set $\hat{\pi}=\hat{\sigma} \circ \pi : \hat{M}\rightarrow M$ , we know $(\hat{\pi}_{\ast}E_{\alpha})^{\ast\ast}=\mathcal{E}_{\alpha}$ and $(\mathcal{Q}_{\alpha})^{\ast\ast}=(\hat{\pi}_{\ast}Q_{\alpha})^{\ast\ast}$. It is easy to see that $\hat{\pi}$ is biholomorphic outside $\hat{\pi}^{-1}(\Sigma_{HNS})$, $E_{\alpha} $ and $Q_{\alpha }$ are isomorphic to $\mathcal{E}_{\alpha }$ and $\mathcal{Q}_{\alpha }$ outside $\tilde{\Sigma}$ respectively.

\medskip

Let $H(t)$ be the long time solution of the Hermitian-Yang-Mills flow (\ref{SSS1}) on the holomorphic vector bundle $\mathcal{E}|_{\Sigma_{\mathcal{E}}}$ with the initial metric $\hat{H}$, and $A(t)$ be the solution of the related Yang-Mills flow (\ref{YM2}) on the Hermitian vector bundle $(\mathcal{E}|_{\Sigma_{\mathcal{E}}}, \hat{H})$ with the initial connection $\hat{A}$.
We have $A(t)=\sigma (t)(\hat{A})$, where $\sigma(t)$ satisfies
$\sigma ^{\ast \hat{H}}(t)\sigma (t)=h (t)=\hat{H}^{-1}H(t)$. Note that Lemma \ref{lemmaBS} says there is a sequence of heat flows $H_{ \epsilon}( t)$ on the holomorphic vector bundle  $E$ which converges successively  to  $H( t)$  in $C_{loc}^{\infty}$-topology outside $\Sigma_{\mathcal{E}}$ as $(\epsilon_{1}, \cdots , \epsilon_{k})\rightarrow 0$. In the sequel, we denote by $\overline{H}_{ \epsilon}( t)=\hat{\sigma}^{\ast}H_{ \epsilon}( t)$ the pull back metric on the bundle $\hat{E}$.

Theorem \ref{thmlimit} and Proposition \ref{part1} imply the part (1) of Theorem \ref{thm 1.1}. So we only need to prove the part (2) of Theorem 1.1.
We assume there is a sequence of connections $A(t_{j})
$ which converges to $A_{\infty}$   in $C^{\infty}_{loc}$-topology outside
$\Sigma$ as $j\rightarrow +\infty$.
Let $\mathcal{S}=\mathcal{E}_{1}$ be the first $\omega$-stable  sub-sheaf corresponding to the above HNS-filtration,   $\mathcal{Q}=\mathcal{E}$, and $g_{j}=\sigma (t_{j})$.
Using the formulas (\ref{HH01}), (\ref{deg}), Lemma \ref{lemmaBS}, Theorem \ref{thm key}, and considering the metrics $\overline{H}_{ \epsilon}( t)$, one can check easily that  the conditions (1), (2), (3), (4) in the above inductive hypotheses are satisfied. Based on Theorem \ref{thmlimit}, we suppose that there exists a sequence of isomorphisms
\begin{equation}
a_{j}: (\mathcal{E}|_{M\setminus \Sigma }, \hat{H})\rightarrow (E_{\infty}, H_{\infty})
\end{equation}
such that $(a_{j}^{-1})^{\ast}(A(t_{j}))\rightarrow A_{\infty}$ in $C_{loc}^{\infty}$-topology outside $\Sigma$ as $j\rightarrow +\infty$. Let $i_{0}: \mathcal{E}_{1}\rightarrow \mathcal{E}$ be the holomorphic inclusion, by Proposition \ref{propmap}, then there is a subsequence of $f_{j}=a_{j}\circ g_{j} \circ i_{0}$, up to rescale, converging to a non-zero holomorphic map $f_{\infty }: \mathcal{E}_{1}\rightarrow (E_{\infty}, \overline{\partial }_{A_{\infty}})$ outside $\Sigma_{HNS}\cup \Sigma_{an}$ as $j\rightarrow +\infty$.
Applying Hartog's
theorem, we can extend $f_{\infty}$  to the whole $M$ as a  sheaf homomorphism.

Let $\pi _{1}^{H(t)}: \mathcal{E}\rightarrow \mathcal{E}$ be the orthogonal projection onto $\mathcal{E}_{1}$ with respect to the Hermitian metric $H(t)$, and $\pi _{1}^{(t)}=\sigma (t)\circ \pi _{1}^{H(t)} \circ \sigma ^{-1}(t)$. Set $\tilde{\pi}_{1}^{j}=a_{j}\circ \pi_{1}^{(t_j)} \circ a_{j}^{-1}$.
From Lemma \ref{lemma3.13}, we know that $ \tilde{\pi}_{1}^{j} \rightarrow \pi_{1 }^{\infty } $  strongly in $L^{p}\cap L_{1, loc}^{2}$ outside $\Sigma_{HNS}\cup\Sigma_{an}$ as $j\rightarrow +\infty$, and $\pi_{1}^{\infty}$ determines a subsheaf $E_{1}^{\infty}$ of $(E_{\infty}, \overline{\partial}_{A_{\infty}})$, with $\rank (E_{1}^{\infty})=\rank (\mathcal{E}_{1})$ and $\mu_{\omega} (E_{1}^{\infty})=\mu_{\omega}(\mathcal{E}_{1})$. Because $\tilde{\pi}_{1}^{j}\circ f_{j}=f_{j}$ for all $j$, we see that in the limit $\pi_{1}^{\infty}\circ f_{\infty}=f_{\infty}$, and then
\begin{equation}f_{\infty}: \mathcal{E}_{1}\rightarrow E_{1}^{\infty}.\end{equation}
Moreover, Theorem \ref{thm key} tells us that $(E_{\infty}, \overline{\partial}_{A_{\infty}})$ and $\mathcal{E}$ have the same HN type, and then the  subsheaf $E_{1}^{\infty}$ is $\omega$-semistable.
Recalling that  $\mathcal{E}_{1}$  is $\omega$-stable, with the result in \cite{Ko2} (V.7.11; 7.12), we observe that the non-zero holomorphic map
$f_{\infty }$ must be injective, then \begin{equation}\mathcal{E}_{1}\simeq
E_{1}^{\infty}=f_{\infty }(\mathcal{E}_{1}),\end{equation} and $E_{1}^{\infty }$ is an $\omega$-stable subsheaf of $(E_{\infty}, \overline{\partial}_{A_{\infty}})$.

Let $\{e_{\alpha }\}$ be a local frame of $\mathcal{E}_{1}$, and $H_{j, \alpha \bar{\beta}}=\langle f_{j}(e_{\alpha }), f_{j}(e_{\beta})\rangle_{\hat{H}}$. We derive
\begin{equation}
\tilde{\pi}_{1}^{j} (X)=\langle X, f_{j}(e_{\beta })\rangle_{\hat{H}}H_{j}^{\alpha , \bar{\beta} }f_{j}(e_{\alpha})
\end{equation}
for any $X \in \Gamma(E)$, where $(H_{j}^{\alpha , \bar{\beta} })$ is the inverse of the matrix $(H_{j, \alpha \bar{\beta}})$.
Because $f_{j}\rightarrow
f_{\infty}$ in $C_{loc}^{\infty }$-topology as $j\rightarrow +\infty$, and $f_{\infty}$ is injective, we can prove that $\tilde{\pi}_{1}^{j}\rightarrow
\pi_{1}^{\infty}$ in $C^{\infty }_{loc }$-topology off $\Sigma_{an}\cup \Sigma_{HNS}$ as $j\rightarrow +\infty$.

\medskip

Consider  the orthogonal holomorphic decomposition $(E_{\infty },
\overline{\partial }_{A_{\infty }})=E_{1}^{\infty }\oplus Q_{\infty}$, where $Q_{\infty}=(E_{1}^{\infty})^{\bot }$. Let $\tilde{\pi}_{1}: E_{\infty} \rightarrow E_{\infty}$ be the projection onto $E_{1}^{\infty}$ with respect to the metric $H_{\infty}$.
Using
 Lemma 5.12 in \cite{Da}, we can choose a sequence of unitary gauge
transformations $\tilde{u}_{j}$ such that
$\tilde{\pi}_{1}^{(j)}=\tilde{u}_{j}\tilde{\pi}_{1}\tilde{u}_{j}^{-1}$  and $\tilde{u}_{j}\rightarrow \Id_{E}$
in $C_{loc}^{\infty}$-topology on $M\setminus (\Sigma_{HNS}\cup \Sigma_{an})$ as $j\rightarrow +\infty$.
It is easy to check that $\tilde{u}_{j}(Q_{\infty})=\tilde{u}_{j}((E_{1}^{\infty})^{\bot })=(\tilde{\pi}_{1}^{(j)} (E_{\infty}))^{\bot}$, and the unitary gauge
transformation $\tilde{u}_{0}: E_{\infty}\rightarrow E_{\infty}$ satisfies $a_{0}^{-1}\circ \tilde{u}_{0}(Q_{\infty})= a_{0}^{-1}\circ \tilde{u}_{0}((E_{1}^{\infty})^{\bot })=\mathcal{E}_{1}^{\bot \hat{H}}$.

Set $\mathcal{Q}=\mathcal{E}/\mathcal{E}_{1}$, then we have $Gr^{HNS}(\mathcal{E})=S\oplus Gr^{HNS}(\mathcal{Q} )$.
Denote by $p^{\ast \hat{H}}=\Id-\pi_{1}^{\hat{H}}: \mathcal{Q} \rightarrow \mathcal{E}_{1}^{\bot \hat{H}}$ the induced  bundle isomorphisms  on $M\setminus \Sigma_{HNS}$, and consider
 the induced
connections on $\mathcal{Q}$
 \begin{equation}D_{A_{j}^{Q}}=(p^{\ast \hat{H}})^{-1}\circ a_{0}^{-1}\circ\tilde{u}_{0}\circ
 \tilde{\pi}_{1}^{\bot}\circ \tilde{u}_{j}^{-1}\circ a_{j}\circ D_{A_{j}}\circ a_{j}^{-1}\circ
\tilde{u}_{j} \circ \tilde{\pi}_{1}^{\bot} \circ \tilde{u}_{0}^{-1}\circ a_{0}\circ p^{\ast \hat{H}},\end{equation}
and the complex gauge transformation
\begin{equation}h_{j}=(p^{\ast \hat{H}})^{-1}\circ a_{0}^{-1}\circ\tilde{u}_{0}\circ \tilde{\pi}_{1}^{\bot}\circ
\tilde{u}_{j}^{-1}\circ a_{j}\circ g_{j} \circ p^{\ast \hat{H}} \in \textbf{G}^{\mathbb{C}}(\mathcal{Q}).\end{equation}
Then it holds that
\begin{equation}
\overline{\partial}_{A_{j}^{Q}}=h_{j}\circ \overline{\partial}_{A_{0}^{Q}} \circ h_{j}^{-1},
\end{equation}
and
\begin{equation}
\partial_{A_{j}^{Q}}=(h_{j}^{\ast})^{-1}\circ \partial_{A_{0}^{Q}}\circ h_{j}^{\ast},
\end{equation}
 where we have used the facts $(\pi_{1}^{(t_j)})^{\bot}\circ g_{j}=(\pi_{1}^{(t_j)})^{\bot} \circ g_{j}\circ (\pi_{1}^{\hat{H}})^{\bot}$ and $h_{j}^{-1}=(p^{\ast \hat{H}})^{-1}\circ (\pi_{1}^{(0)})^{\bot } \circ g_{j}^{-1}\circ a_{j}^{-1}\circ
\tilde{u}_{j} \circ \tilde{u}_{0}^{-1}\circ a_{0}\circ p^{\ast \hat{H}}$.
By the
definition, it is easy to check that $((p^{\ast \hat{H}})^{-1} \circ a_{0}^{-1}\circ \tilde{u}_{0})^{\ast}(A_{j}^{Q})\rightarrow A_{\infty
}^{Q_{\infty}}$ in
$C^{\infty}_{loc }$-topology as $j\rightarrow +\infty$, and  $h_{j}^{\ast } h_{j}=(\hat{H}^{Q})^{-1}H^{Q}(t_{j})$, where $H^{Q}(t)$ denotes the induced metric on the quotient $\mathcal{Q}$ by $H(t)$. Combining (\ref{deg}) and Lemma \ref{lemmaBS}, we get that $\|\sqrt{-1}\Lambda _{\omega }(F_{A_{H(t)}^{Q}})\|_{L^{1}(\omega)}$ is uniformly bounded for $t\geq t_{0} >0$.  So inductive hypotheses (1) and (2) are satisfied.

Let $\mathcal{S}=\mathcal{Q}_{2}=\mathcal{E}_{2}/\mathcal{E}_{1}$, then (\ref{HNS022}) implies the inductive hypothesis (3). Considering the induced metric $\overline{H}_{\epsilon}^{Q}(t)$ on the quotient $\hat{E}/E_{1}$ by $\sigma^{\ast}(H_{\epsilon }(t))$, from the formulas (\ref{HH01}), (\ref{deg}) and Lemma \ref{lemmaBS}, we see that the inductive hypothesis (4) is valid.
Using Theorem \ref{thm key} and Lemma \ref{lemma3.13}, one can check easily that the inductive hypothesis (5) is also valid. Repeating the above argument,  we obtain an isomorphism
\begin{equation}
f: (E_{\infty}, \overline{\partial }_{A_{\infty}})\rightarrow Gr^{HNS
}(\mathcal{E})= \oplus_{\alpha =1}^{l} \mathcal{Q}_{\alpha }
\end{equation}
on $M\setminus (\Sigma_{HNS}\cup \Sigma_{an})$.  By the uniqueness of reflexive extension in \cite{Siu}, we know that $f$ can be extended to a sheaf isomorphism  on the whole $M$.
 This completes the
proof of Theorem \ref{thm 1.1}.

\hfill $\Box$ \\

\hspace{0.4cm}

\hspace{0.3cm}

\end{document}